\documentclass[preprint,11pt,authoryear]{elsarticle}
\usepackage{lineno}
\modulolinenumbers[5]
\usepackage{pgfplots}
\usepackage{multirow}
\usepackage[utf8]{inputenc}
\usepackage{blindtext}
\usepackage[T1]{fontenc}
\usepackage{amsmath}
\usepackage{algorithmic}
\makeatletter
\newcommand{\AlgoResetCount}{\renewcommand{\@ResetCounterIfNeeded}{\setcounter{AlgoLine}{0}}}
\newcommand{\AlgoNoResetCount}{\renewcommand{\@ResetCounterIfNeeded}{}}
\newcounter{AlgoSavedLineCount}

\makeatother

\usepackage{lscape}
\usepackage[dvipsnames]{xcolor}
\usepackage{array}
\newcolumntype{L}[1]{>{\raggedright\let\newline\\\arraybackslash\hspace{0pt}}m{#1}}
\newcolumntype{C}[1]{>{\centering\let\newline\\\arraybackslash\hspace{0pt}}m{#1}}
\newcolumntype{R}[1]{>{\raggedleft\let\newline\\\arraybackslash\hspace{0pt}}m{#1}}

\usepackage{changepage}
\usepackage{subcaption}
\usepackage{graphicx}

\usepackage[normalem]{ulem}
\useunder{\uline}{\ul}{}

\makeatletter
\DeclareFontEncoding{LS1}{}{}
\DeclareFontSubstitution{LS1}{stix}{m}{n}
\DeclareMathAlphabet{\mathscr}{LS1}{stixscr}{m}{n}
\makeatother

\definecolor{blue(pigment)}{rgb}{0.2, 0.2, 0.6}

\makeatletter
\def\@xfootnote[#1]{%
  \protected@xdef\@thefnmark{#1}%
  \@footnotemark\@footnotetext}
\makeatother

\let\originaltable\table
\let\endoriginaltable\endtable
\renewenvironment{table}[1][ht]{%
  \originaltable[#1]
  \centering}%
  {\endoriginaltable}

\usepackage[margin=2.5cm]{geometry}
\usepackage{framed}

\usepackage{amsfonts}
\allowdisplaybreaks

\usepackage{url}

\usepackage{breakurl}
\usepackage[breaklinks]{hyperref}

\usepackage{amssymb}

\usepackage{float}
\usepackage{caption}
\usepackage{multirow}
\usepackage{booktabs}
\usepackage{array} 
\usepackage{makecell}
\usepackage{graphicx}
\usepackage{amsmath}
\PassOptionsToPackage{hyphens}{url}\usepackage{hyperref}
\usepackage{subcaption}
\usepackage[linesnumbered,commentsnumbered,vlined,ruled]{algorithm2e}
\usepackage{bm}

\usepackage{tikz}
\usetikzlibrary{matrix}
\usetikzlibrary{positioning}
\usetikzlibrary{patterns}
\usepackage{tablefootnote}
\usetikzlibrary{calc}
\usetikzlibrary{shapes}
\usetikzlibrary{arrows}

\tikzset{
    vertex/.style={circle,draw,minimum size=0.7em,inner sep=0pt,font=\microsize},
    edge/.style={->,> = latex'}
}

\usepackage{pgfplots}
\usepackage{pgfplotstable}
\usepgfplotslibrary{dateplot}
\pgfplotsset{compat=newest}
\pgfplotsset{
        show sum on top/.style={
            /pgfplots/scatter/@post marker code/.append code={%
                \node[
                    at={(normalized axis cs:%
                            \pgfkeysvalueof{/data point/x},%
                            \pgfkeysvalueof{/data point/y})%
                    },
                    anchor=south,
                ]
                {\pgfmathprintnumber{\pgfkeysvalueof{/data point/y}}};
            },
        },
    }

\usepackage{algorithmic}
\usepackage{ragged2e}
\usepackage{mathtools}
\usepgfplotslibrary{colorbrewer}
\pgfplotsset{compat = 1.15, cycle list/Set1-8} 
\usetikzlibrary{pgfplots.statistics, pgfplots.colorbrewer} 
\usepackage{pgfplotstable}
\usepackage{filecontents}

\makeatletter
\newcommand\resetstackedplots{%
\pgfplots@stacked@isfirstplottrue
}
\makeatother

\newcommand{\microsize}{\fontsize{4pt}{5pt}\selectfont}

\makeatletter
\newcounter{groupcount}
\pgfplotsset{
    draw group line/.style n args={5}{
        after end axis/.append code={
            \setcounter{groupcount}{0}
            \pgfplotstableforeachcolumnelement{#1}\of\datatable\as\cell{%
                \def\temp{#2}
                \ifx\temp\cell
                    \ifnum\thegroupcount=0
                        \stepcounter{groupcount}
                        \pgfplotstablegetelem{\pgfplotstablerow}{[index]0}\of\datatable
                        \coordinate [yshift=#4] (startgroup) at (axis cs:\pgfplotsretval,0);
                    \else
                        \pgfplotstablegetelem{\pgfplotstablerow}{[index]0}\of\datatable
                        \coordinate [yshift=#4] (endgroup) at (axis cs:\pgfplotsretval,0);
                    \fi
                \else
                    \ifnum\thegroupcount=1
                        \setcounter{groupcount}{0}
                        \draw [
                            shorten >=-#5,
                            shorten <=-#5
                        ] (startgroup) -- node [anchor=north] {#3} (endgroup);
                    \fi
                \fi
            }
            \ifnum\thegroupcount=1
                        \setcounter{groupcount}{0}
                        \draw [
                            shorten >=-#5,
                            shorten <=-#5
                        ] (startgroup) -- node [anchor=north] {#3} (endgroup);
            \fi
        }
    }
}
\makeatother
\journal{OR Spectrum}

\begin{document}

\begin{frontmatter}

\title{A tactical time slot management problem under mixed logit demand}

\author[inst1]{Dorsa Abdolhamidi\corref{cor1}\fnref{label2}}
\fntext[label2]{Corresponding author}
\ead{dorsa.abdolhamidi@unil.ch}
\author[inst1]{Virginie Lurkin}
\ead{virginie.lurkin@unil.ch}

\affiliation[inst1]{organization={University of Lausanne, HEC Lausanne Faculty of Business and Economics},
            addressline={Quartier de Chamberonne}, 
            city={Lausanne},
            postcode={1015}, 
            state={Vaud},
            country={Switzerland}}
\begin{abstract}
We study the tactical time slot management problem under mixed logit demand for attended home delivery in subscription settings. We propose a static mixed-integer linear programming model that integrates delivery slot assortment, price discount decisions, and routing optimization while capturing customer heterogeneity through the mixed logit model. To overcome the computational challenges posed by simulation-based choice probabilities, we develop a simulation-based Adaptive Large Neighborhood Search method aligned with a Sample Average Approximation reformulation. Computational experiments on large-scale instances demonstrate the effectiveness of our approach in capturing stochastic customer behavior and preference heterogeneity, providing a scalable and flexible method for optimizing time slot management under complex demand structures.
\end{abstract}
\end{frontmatter}

\section{Introduction}

Online shopping has experienced remarkable growth over the past decade, rising from 7.4\% of total retail sales in 2015 to 18.9\% in 2022 \citep{Coppola2023GlobalStatista}. This trend is expected to continue, with online retail projected to account for approximately 23\% of total sales by 2027. Consumers increasingly demand fast, reliable, and convenient deliveries. To meet these expectations, many retailers offer multiple delivery modes, the most common being in-store pickup, collection at designated locations, and home delivery, which accounts for 79\% of online purchases \citep{Justen2021OnlinePreferences}. Attended home delivery (AHD) is particularly relevant for perishable or high-value goods, or when orders require payment confirmation or signature.

In this context, subscription-based models (SBM) offer an increasingly popular format for AHD. Companies such as HelloFresh\footnote{\url{https://www.hellofresh.com}} and Gousto\footnote{\url{https://www.gousto.co.uk}} deliver fresh ingredients and recipes on a recurring basis. Similarly, local farms distribute weekly produce baskets. These services enhance convenience and quality assurance but pose logistical challenges, particularly in aligning delivery routes with uncertain customer time slot selections.

Although SBMs benefit from access to historical customer data, such as past orders and locations, the inherent uncertainty over which delivery slots will be selected can lead to suboptimal planning and higher routing costs. Effective models must therefore incorporate stochastic choice behavior to ensure operational efficiency and service quality.

Customer time slot selection is largely driven by trade-offs between price and delivery timing \citep{Kapser2021AutonomousAcceptance, Nguyen2019WhatRetailing}. By anticipating how customers balance these dimensions, retailers can optimize both the delivery offerings and underlying routing plans \citep{Campbell2005DecisionInitiatives, Vinsensius2020DynamicDelivery}. The focus here is on a tactical time slot assortment problem, where an e-retailer offering AHD via a subscription model must determine a fixed set of time slots and associated price discounts to offer over a selling horizon. While dynamic approaches exist, static differentiated offerings remain operationally attractive due to their simplicity and ease of implementation \citep{Soppert2022DifferentiatedEffects}.

To model customer behavior, extensive existing literature relies on the Multinomial Logit (MNL) model, valued for its closed-form choice probabilities and computational tractability. However, MNL imposes restrictive assumptions, notably homogeneous preferences across individuals and the independence of irrelevant alternatives, that limit its ability to reflect realistic substitution patterns. Finite mixtures of MNLs, such as latent class models, offer more flexibility by segmenting the population into discrete customer types. While this partially addresses preference heterogeneity, such models still assume a finite number of customer types and cannot capture smoothly varying individual preferences. The Mixed Logit (ML) model overcomes these limitations by introducing random coefficients, allowing for continuous variation in individual preferences and more flexible substitution patterns \citep{Train2009DiscreteSimulation}.

Integrating ML into the slot assortment problem is nontrivial, as the choice probabilities cannot be expressed in closed-form, resulting in nonlinear models that require simulation-based methods. To address this, we build upon the work of \citet{PachecoPaneque2021IntegratingOptimization} by using a reformulation based on the Sample Average Approximation (SAA), a Monte Carlo simulation-based technique. Each scenario represents a draw from the distribution of random coefficients and error terms in the utility function. Evaluating deterministic customer choices for each utility realization eliminates the need for closed-form choice probabilities and enables support for a wide range of Random Utility Models (RUM), such as ML, probit, and other simulation-based discrete choice models.

Using this reformulation, the resulting Tactical Time Slot Management problem under ML demand (TTSM-ML) is formulated as a static Mixed-Integer Linear Programming (MILP) model that simultaneously captures three interrelated decisions: the e-retailer’s selection of delivery slots and price discounts; customers’ time slot choices driven by their heterogeneous preferences modeled via ML; and the e-retailer’s routing optimization based on these choices.

Monte Carlo–based MILP formulations are typically difficult to solve for large instances, especially when the number of scenarios and customers increases. To efficiently solve TTSM-ML, we develop a simulation-based Adaptive Large Neighborhood Search (sALNS) heuristic. sALNS, along with its general form ALNS, has proven effective in solving a variety of complex optimization problems \citep{Ropke2006AnWindows,MattosRibeiro2012AnProblem,Dang2016DesignALNS,Pisinger2019LargeSearch,Cantu-Funes2023Simulation-basedSystems}. Its flexibility, enabled by multiple adaptive operators, allows efficient exploration of the solution space while adapting to the problem’s specific characteristics. Leveraging a simulation-based objective evaluation aligned with the SAA-based MILP structure, sALNS maintains computational scalability while effectively exploring solutions.

To evaluate the performance of our model and solution approach, we conduct a series of numerical experiments on instances of the TTSM-ML problem with up to 100 customers, sizes that are computationally intractable for exact solvers. The experiments compare the sALNS heuristic against exact methods and assess its robustness and computational efficiency in capturing the stochastic nature of customer choices and the underlying heterogeneity modeled by ML. The results demonstrate the relevance of our model and method.

The remainder of the paper is structured as follows. Section~\ref{Sec2} reviews the relevant literature. Section~\ref{sec3} presents the MILP model and its components. Section~\ref{secALNS} describes the sALNS method. Section~\ref{secRes} reports the numerical results. Conclusions are discussed in Section~\ref{sec6}.

\section{Literature review and main contributions} \label{Sec2}

The field of AHD has witnessed extensive research, with comprehensive literature reviews such as those by \citet{Wamuth2023DemandReview} and \citet{Agatz2010TimeDelivery} providing a foundational understanding of demand management strategies. These studies organize the domain across planning levels (strategic, tactical, and operational) and decision levers (assortment and pricing). Within this broader context, our work focuses on tactical planning for time slot management, with Section~\ref{LR} reviewing the main contributions and Section~\ref{contri} highlighting our own contributions.

\subsection{Research on tactical time slot management}
\label{LR}
Tactical Time Slot Management (TTSM) refers to the planning decisions made before the ordering process begins, where retailers determine which delivery slots or prices to offer based on anticipated customer preferences and behavior.

Table~\ref{tb:LR} summarizes key contributions in the TTSM literature, organized along multiple modeling dimensions. \textit{Delivery options offered} distinguishes whether customers are presented with a single delivery option or a menu of alternatives. \textit{Decision levers} capture the planning variables considered (slot assortment, pricing, or both). \textit{Customer granularity} indicates whether decisions are made for individual customers or based on aggregate demand. \textit{Behavioral assumptions} specify whether customer responses are treated deterministically or modeled as stochastic choices. \textit{Choice model type} describes how customer behavior is captured (\textit{e.g.}, uniform rules, utility-based models). Finally, \textit{context / industry} outlines the operational setting of each contribution.

\begin{table}[h]
\caption{Important contributions in tactical time slot management\tablefootnote{The main content of the table is derived from \citet{Wamuth2023DemandReview}, with some changes for coherence.}. MNL: Multinomial Logit; ML: Mixed Logit; SBM: Subscription-Based Model.}
\label{tb:LR}
\resizebox{\textwidth}{!}{%
\begin{tabular}{lcccccc}
\toprule
Paper                                       & Delivery Options Offered & Decision Levers  & Customer Granularity & Behavioral Assumptions & Choice Model Type & Context / Industry         \\ \midrule
\citet{Agatz2010TimeDelivery}               & Multiple       & Assortment & Aggregate   & Deterministic & Uniform          & E-grocery               \\
\citet{Cleophas2014WhenProfitable}          & Multiple       & Assortment & Aggregate   & Deterministic & -          & E-grocery               \\
\citet{Hernandez2017HeuristicsView}         & Multiple       & Assortment & Aggregate   & Deterministic & -  & E-retail                \\
\citet{Visser2019StrategicRetailing}        & Single         & Assortment & Aggregate   & Deterministic & -          & E-grocery             \\
\citet{Spliet2014TheProblem}                & Single         & Assortment & Individual   & Deterministic & -          & B2B SBM                 \\
\citet{Spliet2015TheProblem}                & Single         & Assortment & Individual   & Deterministic & -          & B2B SBM                 \\
\citet{Spliet2017TheTimes}                  & Single         & Assortment & Individual   & Deterministic & -          & B2B SBM                 \\
\citet{Fallahtafti2021TimeWindows}          & Single         & Assortment & Individual   & Deterministic & -          & B2B                     \\
\citet{Bruck2018AServices}                  & Multiple       & Assortment & Individual  & Deterministic & Rule-based & Service                 \\
\citet{Karaenke2020Non-monetaryDelivery}    & Single         & Assortment & Individual   & Deterministic & Valuation-based & B2B              \\
\citet{Cote2019TacticalDelivery}            & Multiple       & Assortment & Aggregate   & Stochastic & Uniform  & E-retail          \\
\citet{Mackert2019IntegratingLogistics}     & Multiple       & Assortment & Aggregate   & Stochastic    & Finite-mixture MNL       & E-grocery               \\
\citet{Klein2019DifferentiatedDelivery}     & Multiple       & Pricing      & Aggregate   & Deterministic & Rank-based & E-grocery               \\
\citet{Zhang2025IntegratedDelivery}         & Multiple       & Pricing      & Aggregate   & Deterministic & Rank-based & Food               \\ \midrule
Our paper                                   & Multiple       & Assortment \& Pricing & Individual & Stochastic & ML    & SBM          \\
\bottomrule
\end{tabular}%
}
\end{table}

Early contributions to TTSM typically adopted aggregate demand representations and deterministic behavioral assumptions, focusing on operational efficiency rather than preference-driven decision-making. \citet{Agatz2010TimeDelivery}, \citet{Cleophas2014WhenProfitable}, \citet{Hernandez2017HeuristicsView} applied zone-based models without accounting for customer behavior. \citet{Agatz2010TimeDelivery} examined time slot allocation in e-grocery. Their model assumed uniform demand across slots without any behavioral consideration. Similarly, \citet{Cleophas2014WhenProfitable} developed an iterative framework for setting slot capacities independently of customer demand responses. \citet{Hernandez2017HeuristicsView} proposed heuristic methods to solve large-scale problems. \citet{Visser2019StrategicRetailing} offered only a single time slot, optimizing acceptance rates and profitability under deterministic customer responses.

A parallel stream addressed B2B contexts with individual-level granularity but without modeling customer preferences. This includes \citet{Spliet2014TheProblem, Spliet2015TheProblem, Spliet2017TheTimes}, and \citet{Fallahtafti2021TimeWindows}, all of which assumed demand follows the only single time slot optimally offered by the supplier without explicit modeling of customer preferences or choice behavior.

Advancing customer-level modeling, \citet{Bruck2018AServices} introduced rule-based slot preferences derived from historical data, assigning individual customers fixed selection rules. Similarly, \citet{Karaenke2020Non-monetaryDelivery} used valuation-based assumptions for user choices. Yet, both relied on deterministic behavior.

Stochastic demand-side modeling first emerged with \citet{Cote2019TacticalDelivery}, who treated customer selections probabilistically using a uniform choice model. Building on this, \citet{Mackert2019IntegratingLogistics} introduced a finite-mixture MNL model to explicitly represent variation and diversity in customers’ slot choices. Their approach, relying on the closed-form MNL probability formulation, enables segmentation based on observable traits; however, it assumes homogeneous preferences within segments and does not model individual-level customer decisions.

Other contributions focused on pricing as a key tactical control lever. \citet{Klein2019DifferentiatedDelivery} and \citet{Zhang2025IntegratedDelivery} explicitly modeled customer price sensitivity using rank-based utility models, allowing retailers to influence slot selection through monetary incentives. While highlighting the importance of price differentiation, both approaches relied on deterministic assumptions and aggregate-level modeling, limiting their ability to capture nuanced, individual-level choice behavior under uncertainty.

In summary, although the literature has gradually evolved towards more realistic demand modeling, most approaches still depend on deterministic assumptions, aggregate-level planning, or discrete customer segments. To the best of our knowledge, no prior TTSM work integrates continuous individual heterogeneity via the ML model, which our paper addresses. 

\subsection{Our contribution}
\label{contri}

Our contributions to the TTSM literature are as follows:
\begin{itemize}
    \item[-] \textit{Demand-side modeling with continuous individual heterogeneity: }
    We develop a TTSM problem that incorporates individual-level customer preferences. This extends existing models that either omit choice behavior or rely on closed‑form specifications such as finite‑mixture MNL, which capture at most segment‑level (discrete) heterogeneity. In contrast, ML enables continuous heterogeneity in preferences, thus supporting customer‑specific choice modeling.
    \item[-] \textit{Scalable solution method integrating simulation-based RUMs:} 
    We rely on SAA to embed ML choice probabilities within a stochastic MILP formulation. To solve this model at scale, we develop a two-phase heuristic designed for simulation-based RUMs. First, the \textit{Route-First Time-Second (RFTS)} heuristic constructs initial solutions by clustering customers via K-means, sequencing them with a nearest neighbor search, and assigning slot–price menus using utility-aware rules and bidirectional time-window checks. While RFTS adopts the familiar two-phase “route‑then‑time” logic of \citet{Klein2023DynamicDelivery}, it differs in context and mechanics: RFTS operates offline, planning all routes, offerings, and prices before any bookings, leverages ML parameters rather than closed-form MNL, and expands feasible slot sets via bidirectional window computations. Second, an sALNS iteratively refines the solution using Monte Carlo demand evaluation, custom destroy-repair operators, and a record-to-record travel acceptance mechanism.
    \item[-] \textit{Numerical validation:} 
    We perform comprehensive numerical experiments to assess the performance of our solution approach, and the value of modeling individual-level choice heterogeneity within tactical slot management.
\end{itemize}

\section{Mathematical formulation}
\label{sec3}

This section presents the mathematical formulation of the TTSM-ML problem as a joint optimization model that captures the interaction between the e-retailer (leader) and customers (followers) in the tactical management of delivery time slots. The model integrates three interdependent decision components:

\begin{enumerate}
    \item \textbf{Slot and Pricing Decision (Retailer)} 
    The retailer selects the set of available delivery slots and corresponding price discounts. This decision is denoted by \( \gamma \in I \), where \( I \) defines the feasible combinations of time slot offerings and associated discounts.

    \item \textbf{Customer Choice (Follower)} 
   Given the retailer's offerings \( \gamma \), each customer chooses the time slot that maximizes their utility. This behavior is represented as \( w^\ast(\gamma) = \arg \max_{w \in \mathbb{W}(\gamma)} U(w) \), where \( \mathbb{W}(\gamma) \) is the set of feasible customer choices conditional on the offerings.

    \item \textbf{Routing Decision (Retailer)}  
    Based on the resulting customers selections \( w^\ast(\gamma) \), the retailer determines optimal delivery routes. This routing decision, denoted by \( x \in \mathbb{X} \), is dependent on both the initial offerings \( \gamma \) and the realized customer choices \( w^\ast(\gamma) \). The feasible routing decisions are represented as \( \mathbb{X}(\gamma, w^\ast(\gamma)) \).
\end{enumerate}

The retailer aims to jointly optimize the delivery slot offerings, pricing, and routing decisions in response to customer choices. The objective function \( \mathcal{Q} \) depends on the retailer’s decisions \( \gamma \) and \( x \), as well as the induced customer selections \( w^\ast(\gamma) \). The retailer’s optimization problem is formulated as:
\[
\max_{\gamma \in I,\, x \in \mathbb{X}(\gamma, w^\ast(\gamma))} \mathcal{Q}(\gamma, w^\ast(\gamma), x).
\]

We formalize the TTSM-ML optimization problem through a MILP formulation, divided into three key components: Section~\ref{problemdef1} addresses the retailer’s decisions regarding slot availability and pricing; Section~\ref{problemdef2} models customers' preferences using a simulation-based RUM; and Section~\ref{secModel} integrates both stakeholders’ decisions into a unified optimization model, where the retailer seeks to maximize profit while adhering to routing constraints and associated costs.\footnote{{
TTSM-ML shares structural features with two-stage stochastic programming and bilevel optimization, but with important differences. Unlike classical two-stage models with exogenous uncertainty, our customer choices are endogenous, shaped by retailer decisions. Routing is evaluated over simulated selections rather than optimized as recourse. Unlike standard bilevel models with deterministic follower responses, we model stochastic customer behavior and interdependent multi-level decisions.}}

\subsection{Assortment and pricing decisions of the retailer}
\label{problemdef1}

An online retailer operating in a subscription-based system must decide, before capturing any orders and statically, which time slots and corresponding price discounts to offer to customers, aiming to maximize profit. These customers are known but uncertain in their final delivery time choices. Let $C$ represent the set of customers, $T$ the set of available time slots, and $H$ the set of possible discrete price discount rates, where each $h \in H$ satisfies $0 \leq h < 1$. The discounts are applied to a base fee $f$, representing the price before any discounts are applied. The set $I$ containing all alternatives offered to customers is then formed by combining each time slot with each discount rate, along with an opt-out alternative (denoted by $0$). Mathematically, we define $I = T \times H \cup \{0\}$.

When considering an element $i \in I$, which is not the opt-out alternative, we define it as an object of the form $(t_i, h_i) \in T \times H$. The retailer's assortment decision is captured by the binary variable $\gamma_{in}$, defined as:
\begin{equation*}
    \gamma_{in} =
    \begin{cases}
      1, & \text{if alternative $i \in I$ is offered to customer $n \in C$}, \\
      0, & \text{otherwise}.
    \end{cases}
\end{equation*}

We also define the set $I_t = \{i \,|\, t_i = t\} \subset I$ as the set of alternatives associated with time slot $t$.

The assortment decision problem for the retailer is then formulated as:
\begin{align}
    \max \mathcal{Q}(\gamma) && \label{eq:master_obj}\\
    \text{subject to} \notag && \\
    \gamma_{0n} &= 1 && \forall n \in C \label{eq:optout}, \\
    \sum_{i \in I} \gamma_{in} &\geq \nu && \forall n \in C, \label{eq:service_quality}\\
    \sum_{i \in I_t} \gamma_{in} &\leq 1 && \forall n \in C, \, \forall t \in T, \label{eq:one_price}\\
    \gamma_{in} &\in \{0, 1\} && \forall n \in C,\, \forall i \in I. \label{eq:bound_4}
\end{align}

Constraints (\ref{eq:optout}) ensure that the opt-out alternative is always available to customers. Constraints (\ref{eq:service_quality}) provide flexibility in determining the minimum number of alternatives, $\nu \in \mathbb{N}$, that the retailer must offer. For instance, by setting $\nu \geq 2$, the retailer guarantees at least one alternative to the opt-out alternative, contributing to customer satisfaction and retention. Constraints (\ref{eq:one_price}) prevent offering the same time slot with different prices, which would lead customers to always choose the cheaper alternative. Constraints (\ref{eq:bound_4}) define the binary decision variables.

The assortment and pricing model is designed to leverage differentiated pricing, allowing the retailer to tailor assortment and discount strategies to drive customers toward more profitable and efficient delivery decisions. To better understand the benefits of such differentiation, we compare it to a case where uniform pricing is applied across customers within the same time slot. This comparison is captured by the additional Constraints (\ref{eq:discrimination}), which ensure that all customers selecting the same time slot receive the same price, thereby limiting differentiation to the time slot level.
Specifically, the constraints enforce that the prices offered to any two customers \( n \) and \( m \) must be equal if both are assigned the same slot. The use of a large constant \( M^f \) allows the constraints to be inactive unless both customers are actually offered slot \( t \):
\begin{equation}
\label{eq:discrimination}
    -M^f \left(2 - \sum_{i \in I_t} (\gamma_{in} + \gamma_{im})\right) \leq \sum_{i \in I_t} h_i f (\gamma_{in} - \gamma_{im}) \leq M^f \left(2 - \sum_{i \in I_t} (\gamma_{in} + \gamma_{im})\right), \quad \forall n, m \in C, \, \forall t \in T,
\end{equation}
where \( M^f = f \), ensuring that the maximum possible price difference does not exceed the base price \( f \). The impact of maintaining price consistency is analyzed in Section~\ref{secFairnessCheck}.

\subsection{Time Slot Decision of the Customers}
\label{problemdef2}

To model how customers respond to the retailer’s slot offerings, we assume they behave as rational utility maximizers. Their decisions are captured by a logit choice model, in which each customer selects the alternative \(i \in I\) that yields the highest utility. Each utility value \(u_{in}\) consists of a systematic component \(V_{in}\), reflecting observable attributes and customer characteristics, and a random error term \(\xi_{in}\), which captures unobservable factors affecting decision-making, including specification and measurement errors:
\begin{equation*}
    u_{in} = V_{in} + \xi_{in}, \quad \forall i \in I,\;\forall n \in C,
\end{equation*}
where \(\xi_{in} \sim \mathrm{EV1}(0,1)\). 

In the case of slot assortment, each alternative \(i\) corresponds to a delivery time slot and is characterized by its time, price, and potentially other attributes. The systematic component of the utility function is given by:
\begin{equation*}
    V_{in} = f_{\mathrm{time}}(\beta^{\mathrm{time}}_{i}, t_i)
           + f_{\mathrm{price}}(\beta^{\mathrm{price}}, h_{in}, f)
           + f_{\mathrm{other}}(\beta^{\mathrm{other}}, \mathbf{o}_i),
\end{equation*}
where:
\begin{itemize}
    \item[-] \(t_i\) is the delivery time associated with alternative \(i\),
    \item[-] \(h_{in}\) represents the discount rate applied to time slot \(i\) for customer \(n\), and \(f\) is the base delivery fee before any discount,
    \item[-] \(\mathbf{o}_i\) is a vector of other observable characteristics of alternative \(i\),
    \item[-] \(\beta^{\mathrm{time}}_{i}\), \(\beta^{\mathrm{price}}\), and \(\beta^{\mathrm{other}}\) are customer-specific sensitivity parameters for time, price, and other attributes, respectively.
\end{itemize}

Under the assumption that customers are utility maximizers, the probability \(P_{in}\) that customer \(n\) selects alternative \(i\) depends both on its availability \(\gamma_{in}\) and on the likelihood that it offers the highest utility among all available options:
\begin{equation*}
    P_{in}
    = \gamma_{in}\,\Pr\bigl[u_{in} \ge u_{jn},\;\forall j \in I \mid \gamma_{jn} = 1\bigr], \quad \forall i \in I,\;\forall n \in C.
\end{equation*}

To account for individual-level heterogeneity in preferences, acknowledging that customers differ in their sensitivities, we allow certain utility parameters \(\beta\) to vary randomly across customers according to a distribution \(f(\beta \mid \Theta)\), where \(\Theta\) denotes the distribution parameters. This approach leads to the well-known ML model, where the choice probability is expressed as:
\begin{equation*}
    P_{in} = \int \frac{\gamma_{in} e^{V_{in}(\beta)}}{\sum_{j \in I} \gamma_{jn} e^{V_{jn}(\beta)}} f(\beta \mid \Theta) \, d\beta, \quad \forall i \in I, \; n \in C,
    \label{eq:choiceprob_mixed}
\end{equation*}
which typically requires numerical integration for evaluation.

To address the nonlinearity of ML choice probabilities, we adopt the SAA reformulation proposed by \citet{PachecoPaneque2021IntegratingOptimization}. This Monte Carlo–based method replaces the integral over the distribution of customer preferences with an empirical average over \(R\) independently sampled scenarios. Each scenario corresponds to a draw of customer-specific preference parameters and error terms, enabling a tractable formulation that embeds flexible choice behavior within a linear optimization model. For each scenario \(r\), the utility of customer \(n\) for alternative \(i\) is computed as:
\begin{equation}
    u_{inr}
    = f_{\mathrm{time}}(\beta^{\mathrm{time}}_{i}, t_i)
    + f_{\mathrm{price}}(\bar{\beta}^{\mathrm{price}}_{r}, h_i, f)
    + f_{\mathrm{other}}(\bar{\beta}^{\mathrm{other}}_{r}, \mathbf{o}_i)
    + \bar{\xi}_{inr}
    - M_{nr}(1-\gamma_{in}), \quad \forall i \in I, \; n \in C,\; r\in R,
    \label{eq:utility_gamma}
\end{equation}
where
\[
M_{nr} = \max_{i \in I} \bigl( V_{inr}(\bar{\beta}_r) + \bar{\xi}_{inr} \bigr), \quad \forall n \in C,\; r\in R.
\]
In each scenario, unavailable alternatives are penalized using the large constant \(M_{nr}\), effectively removing them from consideration.

To capture customer choice within each scenario, we introduce binary decision variables \(w_{inr}\) indicating whether customer \(n\) selects alternative \(i\) in scenario \(r\). Additionally, continuous variables \(U_{nr}\) represent the maximum utility attained by customer \(n\) in scenario \(r\) among all available alternatives. The following constraints enforce logical consistency between these variables:
\begin{align}
    \sum_{i \in I} w_{inr} &= 1, 
    && \forall n \in C, \; r=1,\dots,R, 
    \label{eq:one_choice} \\
    U_{nr} &\ge u_{inr}, 
    && \forall i \in I, \; n \in C, \; r=1,\dots,R, 
    \label{eq:max_utiliy} \\
    U_{nr} &\le u_{inr} + M^U_{nr} (1 - w_{inr}), 
    && \forall i \in I, \; n \in C, \; r=1,\dots,R, 
    \label{eq:UtilityLin_2} \\
    w_{inr} &\in \{0, 1\}, 
    && \forall i \in I, \; n \in C, \; r=1,\dots,R. 
    \label{eq:bound_w2}
\end{align}

Here, \(M_{nr}^U = 2 M_{nr}\) are sufficiently large constants used to deactivate constraints for unavailable alternatives. By doubling the baseline \(M_{nr}\), it effectively cancels out the utility reduction introduced earlier, restoring the original utility value and thus rendering the constraint inactive for those alternatives. Constraints \eqref{eq:one_choice} ensure that each customer selects exactly one delivery alternative in each scenario. Constraints \eqref{eq:max_utiliy} and \eqref{eq:UtilityLin_2} enforce that the selected alternative corresponds to the maximum utility experienced by the customer in that scenario. Finally, Constraints \eqref{eq:bound_w2} define the domain of the binary selection variables.

Using these variables, empirical choice probabilities can be approximated as
\begin{equation*}
    P_{in} \approx \frac{1}{R} \sum_{r=1}^R w_{inr}, \quad \forall i \in I, \; n \in C,
    \label{eq:prob_empirical}
\end{equation*}
which enables the integration of ML and, more generally, any random utility-based choice model, regardless of whether closed-form choice probabilities exist, directly within a linear optimization model for the retailer’s joint assortment, pricing, and routing decisions, as further detailed in the next section.

\subsection{Routing decisions of the retailer}
\label{secModel}
The expected profit combines revenues from customer choices with
routing costs incurred for serving them. As outlined in the previous section, the probability that customer $n \in C$ chooses alternative $i \in I$ is approximated via a sample average over $R$ scenarios. Accordingly, expected revenues are computed as:

\begin{equation}
\frac{1}{R} \sum_{r\in R} \sum_{n\in C}\sum_{i \in N} w_{inr} h_{in}f.\notag
\end{equation}

The retailer's expected routing costs are also influenced by customers' choices, as routes are determined by the promised delivery times. We formulate a Capacitated Vehicle Routing Problem with Time Windows (CVRPTW) to reflect these decisions. Let $K$ denote the set of vehicles, each with capacity $Q$, and let the depot be denoted by $\{O\}$. The full node set is $V=C\cup \{O\}$. Each customer $n\in C$ has demand $d_n$, and travel between any two nodes $m, n \in V$ incurs time $t_{mn}$ and cost $c_{mn}$, proportional to travel time. Each delivery option $i \in I$ has a time window [$\underline{T}_i, \overline{T}_i$].

In each scenario $r$, routing is modeled using binary variables $x_{mnkr}$ and $z_{kr}$, defined as:
\begin{align}
   x_{mnkr} &=
    \begin{cases}
      1, & \text{if vehicle $k$ visits customer $n$ after customer $m$ in scenario $r\in R$}, \\
      0, & \text{otherwise},
    \end{cases} \notag \\
   z_{kr} &=
    \begin{cases}
      1, & \text{if vehicle $k$ is used in scenario $r\in R$}, \\
      0, & \text{otherwise}.
    \end{cases} \notag
\end{align}

The actual delivery time of each customer $n \in C$ in each scenario $r \in R$ is obtained by continuous real variables $\tau_{nr}$. The routing constraints are given by:
\begin{align}
    \sum_{k\in K} \sum_{\substack{m \in V \\ m\neq n}} x_{mnkr} &= \sum_{i\in I \setminus \{0\}} w_{inr}, && \quad \forall \, n \in C,~\forall \, r\in R, \label{eq:customer_once}\\
    \sum_{\substack{m \in V \\ m \neq n}} x_{mnkr} &= \sum_{\substack{m \in V \\ m \neq n}} x_{nmkr}, && \quad \forall \, n \in V,~\forall \, k \in K,~\forall \, r\in R, \label{eq:enter_leave}\\
    \sum_{n\in V} x_{Onkr} &= 1, && \quad \forall \, k \in K ,~\forall \, r\in R, \label{eq:other_origions}\\
    \sum_{i\in I} \overline{T}_{i}w_{inr} &\geq \tau_{nr} ,  &&\quad \forall \, n \in C,~\forall \, r\in R, \label{eq:connection_1}\\
     \sum_{i\in I} \underline{T}_{i}w_{inr} &\leq \tau_{nr} ,  &&\quad \forall \, n \in C,~\forall \, r\in R, \label{eq:connection_2}\\
    \tau_{mr} + t_{mn} &\leq \tau_{nr} +M_{mn}^\tau\left(1-\sum_{k\in K} x_{mnkr}\right), && \quad \forall \, m\in V,~\forall \,n \in C,~\forall \, r\in R, 
    \label{eq:time_continuety} \\
    \sum_{n\in C}d_n\sum_{m\in V}x_{mnkr} &\leq Q z_{kr}, && \quad \forall \, k \in K,~\forall \, r\in R, \label{eq:capacity} \\
    \sum_{m,n\in V} x_{mnkr} &\leq M^K \sum_{m,n\in V} x_{mn(k-1)r} , &&\quad \forall k \in K \setminus\{1\} ,~\forall \, r\in R,\label{eq:symmetry} \\
    x_{mnkr} &\in \{0,1\}, && \quad \forall \, m,n \in V,~\forall \, k \in K,~\forall \, r\in R, \label{eq:bound_x}\\
    \tau_{nr} &\geq 0, && \quad \forall \, n \in V,~\forall \, r\in R, \label{eq:bound_tau}\\
    z_{kr} &\in \{0,1\}, && \quad \forall \, r\in R. \label{eq:bound_z}
\end{align}

Constraints (\ref{eq:customer_once}) ensure that all customers who select a delivery time slot (\textit{i.e.}, the ones who do not opt out) are visited by a vehicle. Constraints (\ref{eq:enter_leave}) are the standard flow conservation constraints. Constraints (\ref{eq:other_origions}) make sure all routes start from the depot ${O}$. Constraints (\ref{eq:connection_1})-(\ref{eq:connection_2}) verify that in each scenario the customer is indeed visited during the promised time slot. Time-continuity is handled by Constraints (\ref{eq:time_continuety}), where $M_{mn}^\tau = \overline{T}_{|T|} + t_{mn}$, make sure that the constraints are not binding when $x_{mnkr}=0$. Constraints (\ref{eq:capacity}) prevent vehicles from carrying loads that exceed their capacities, they also force $z_{kr}=1$ if a vehicle is used to serve customers. Constraints (\ref{eq:symmetry}) break the symmetry of a homogeneous fleet search for each route. $M^K=|V|$ is the maximum number of arcs a vehicle can potentially traverse.
Finally, Constraints (\ref{eq:bound_x})-(\ref{eq:bound_z}) are the domain constraints. 

Assuming the retailer pays a fixed vehicle cost $c^v$, the expected routing costs are given by:
\begin{equation*}
    \frac{1}{R} \left( \sum_{r\in R} \sum_{m,n\in V}\sum_{k \in K} x_{mnkr}  c_{mn} + \sum_{r\in R}\sum_{k \in K} z_{kr} c^v \right),
\end{equation*}

Consequently, the expected profit of the online retailer is given by:
\begin{equation}
\label{eq:Obj1} 
  \frac{1}{R}\left(\sum_{r\in R} \sum_{n\in C}\sum_{i \in N} w_{inr}  h_if - \sum_{r\in R} \sum_{m,n\in V}\sum_{k \in K} x_{mnkr}  c_{mn}  - \sum_{r\in R}\sum_{k \in K} z_{kr} c^v \right),
\end{equation}

The complete TTSM-ML problem maximizes the profit function (Equation~\ref{eq:Obj1}) through a structured integration of three decision components:
\begin{align}
    & \quad \text{Constraints (\ref{eq:optout})-(\ref{eq:one_price}), (\ref{eq:discrimination})} \rightarrow \text{Retailer's assortment \& price discounting decisions,} \notag \\
    & \quad \text{Constraints (\ref{eq:utility_gamma})-(\ref{eq:UtilityLin_2})} \rightarrow \text{Customer choice responses,} \notag \\
    & \quad \text{Constraints (\ref{eq:customer_once})-(\ref{eq:symmetry})} \rightarrow \text{Retailer's routing optimization,} \notag \\
    & \quad \text{Constraints (\ref{eq:bound_4}), (\ref{eq:bound_w2}), (\ref{eq:bound_x})-(\ref{eq:bound_z})} \rightarrow \text{Variable domains.} \notag
\end{align}

The ML choice model allows for rich representation of customer preferences, but embedding its nonlinear choice probabilities via the SAA reformulation results in a large-scale stochastic MILP, which motivates the heuristic solution approach presented in the next section.

\section{Solution approach} \label{secALNS}

In this section, we present the solution approach to address the TTSM-ML problem. As described in Section~\ref{subsec:initial}, a constructive heuristic, referred to as RFTS, is employed to rapidly generate good feasible solutions. In the subsequent improvement phase, detailed in Section~\ref{subsec:ALNS}, we apply an sALNS to refine and enhance the initial solution produced by RFTS. 

\subsection{Constructive phase: RFTS heuristic} \label{subsec:initial}

Our constructive heuristic, called RFTS, first generates a route plan and then allocates time slots. It adapts the Cluster-First Route-Second (CFRS) heuristic from \citet{Gillett1974AProblem}, which consists of two main steps. A detailed description of the heuristic is provided in Algorithm~\ref{alg:init}. 

\begin{algorithm}[]
\DontPrintSemicolon
\let\oldnl\nl%
\newcommand{\nonl}{\renewcommand{\nl}{\let\nl\oldnl}}
\BlankLine
$R \leftarrow$ Cluster customers to establish routes \label{algline:cluster}\;

Check capacity constraints in each $route \in R$ and repair if necessary \label{algline:capacity}\;
\BlankLine
\nonl \For{$route \in R$} {
    \BlankLine
    Solve a Nearest Neighbor Search in $route$ \label{algline:NN}\;
    \BlankLine
    \nonl \For{$c \in Customers$}{
        $\underline{T}_c \leftarrow$ Assign the earliest time $c$ can be met by timing customers in a forward manner starting from the beginning of the horizon \label{algline:forward_t}\;
        
        $\overline{T}_c \leftarrow$ Assign the latest time $c$ can be met by timing customers in a backward manner starting from the end of the horizon \label{algline:backward_t}\;
        
        $T_c \leftarrow$ List of all time slots $t$ such that $\underline{T}_c \leq t \leq \overline{T}_c$ \label{algline:all_t}\;
        \BlankLine
        \nonl \For{$t \in T_c$}{
            $h_t \leftarrow$ The lowest discount rate which on average has a higher probability than the probability of customer $c$ opting out \label{algline:pricing}\;
            
            $H_c \leftarrow H_c \bigcup \{h_t\}$\;
        }
        Assign $(T_c, H_c)$ to customer $c$ \label{algline:option_assignment}\;
    }
}
\caption{RFTS heuristic}
\label{alg:init}
\end{algorithm}

Firstly, the K-means algorithm clusters customers based on their geographical proximity, prioritizing the spatial efficiency of vehicle routes while ignoring temporal delivery constraints. We explore various configurations for the number of clusters, each corresponding to a vehicle's route. Among these, we select the configuration with the lowest delivery cost that keeps the number of vehicles within the interval $[|K|-\zeta, |K|]$ (line~\ref{algline:cluster}), where $K$ is the maximum number of available vehicles.

Afterwards, we check if all clusters meet capacity constraints. If a cluster exceeds the limit, we remove the fewest possible customers to bring the demand down to the allowed level. Then, we see if these customers can be added to other clusters without exceeding their capacities. If this is not possible, we create a new cluster for these customers (line~\ref{algline:capacity}).
Following clustering, 
routes are constructed within each cluster using the nearest neighbor search to reduce travel distance. Routes originate from the depot, servicing all assigned customers within the cluster before returning to the depot. This routing approach aims to reduce travel distance, decreasing routing costs (line~\ref{algline:NN}).

Once routes are constructed, time slots are assigned to each customer using a forward assignment process. Starting from the depot at the beginning of the operational day, the earliest possible time slot ($\underline{T}_c$) for each customer is determined based on cumulative travel time along the route. This forward assignment ensures that each customer's earliest service time is captured as the vehicle progresses through its route (line~\ref{algline:forward_t}).

Following the forward phase, the process is reversed, starting with the last customer on the route. This reverse assignment adjusts the latest possible time slots ($\overline{T}_c$) by backtracking through the route and recalculating time slots against the direction of travel (line~\ref{algline:backward_t}). If discrepancies exist between the earliest ($\underline{T}_c$) and latest ($\overline{T}_c$) time slots, we offer the complete range from $\underline{T}_c$ to $\overline{T}_c$ to the customer, allowing them to select the most convenient time within this range (line~\ref{algline:all_t}). This dual-direction assignment process increases the flexibility offered to customers and improves overall routing efficiency.

For the assortment decision, we first assess the utility of each slot given various discount rates, with the aim of increasing overall revenue. Our strategy involves setting the lowest feasible discount rate for each slot, ensuring it still has greater value than the alternative for the customer to opt out. We achieve this by ranking the utility of all available alternatives against the opt-out choice in each scenario. We calculate averages of these rankings, and the slot that offers the lowest discount rate yet still outperforms the opt-out alternative rank is chosen for the specific slot of that customer (lines~\ref{algline:pricing}-\ref{algline:option_assignment}). 

The RFTS heuristic constructs feasible solutions by first prioritizing routing costs and then revenue, pursuing these objectives sequentially. This approach may overlook important trade-offs in the TTSM-ML problem. To address these limitations, we introduce an improvement phase that refines the solutions for better overall quality.

\subsection{Improvement phase: sALNS method} \label{subsec:ALNS}

This section describes the sALNS method for the improvement phase. Section~\ref{subsec:structure} outlines the overall structure of the sALNS algorithm, while Section~\ref{secNum} focuses on the destroy and repair operators. And Section~\ref{subsec:ls} details the local search component. Each section contributes to the overall solution refinement process.

\subsubsection{Structure of the sALNS algorithm} \label{subsec:structure}

Building on the initial solutions from the RFTS heuristic, we develop an sALNS method that jointly refines routing, slot offerings, and pricing decisions through an integrated, adaptive search. 

In our sALNS, a solution represents a complete plan specifying which time slots are offered to each customer and their corresponding discount rates. Candidate solutions are evaluated using a Monte Carlo simulation that accounts for customer choice behavior across multiple scenarios. For each scenario, customer preferences are sampled from ML model, or more generally, any simulation-based RUM, and delivery slot choices are simulated accordingly. Given these choices, routing costs are estimated using a heuristic to solve the resulting vehicle routing problem with time windows. Several routing heuristics can be applied, including Clarke and Wright’s savings heuristic (CW), its improved variant with local search (ICW), and the CFRS heuristic, which is tailored to clustering-based routing. The overall evaluation procedure, summarized in Algorithm~\ref{alg:obj}, computes the expected profit of a candidate solution by aggregating simulated revenues and routing costs over all scenarios.

\begin{algorithm}
\DontPrintSemicolon

\ForEach{scenario $r \in R$}{
    \ForEach{customer $n \in C$}{
        \ForEach{option $i \in I$}{
            $\bar{\xi}_{inr} \sim \text{EV1}(0,1)$\;
            $\bar{\beta} \sim f(\beta \mid \Theta)$\, where $\bar{\beta} = [\bar{\beta}^{\mathrm{time}}_{i}, \bar{\beta}^{\mathrm{price}}_{nr}, \bar{\beta}^{\mathrm{other}}_{nr}]$\;
            $u_{inr} = f_{\mathrm{time}}(\bar{\beta}^{\mathrm{time}}_{i}, t_i) 
            + f_{\mathrm{price}}(\bar{\beta}^{\mathrm{price}}_{nr}, h_{in}, f) 
            + f_{\mathrm{other}}(\bar{\beta}^{\mathrm{other}}_{nr}, \mathbf{o}_i) + \bar{\xi}_{inr}$\;
        }
        $t_n, h_n \gets \arg\max_{i \in I} u_{inr}$\;
    }
    \textnormal{Cost$_r$ $\gets$ Solve a routing problem with time windows $\{t_n\}_{n \in C}$}\;
    \textnormal{Revenue$_r$ $\gets \sum_{n \in C} h_nf$}\;
}
\textnormal{Obj $\gets \frac{1}{R}\sum_{r \in R} (\text{Revenue}_r - \text{Cost}_r)$}\;

\caption{Scenario‑based expected objective evaluation under individual‑level random utility demand and routing}
\label{alg:obj}
\end{algorithm}

The overall structure of the sALNS is outlined in Algorithm~\ref{alg:pseudocodeALNS}, and its key components, the destroy and repair operators and the local search phase, are detailed in the sections that follow.

\begin{algorithm}
\DontPrintSemicolon

solution$_i$, solution$_b \leftarrow$ Construct the initial solution with algorithm~\ref{alg:init} \label{algline:init_sol}\;
obj$_i$, obj$_b \leftarrow$ Evaluate solution$_i$ with algorithm~\ref{alg:obj} \label{algline:f_calc}\;
$\omega_d, \omega_r$ $\leftarrow$ Initialize weights for destroy and repair operators equally \label{algline:param_w}\;
$\epsilon, N, \phi, \phi_{min}, \nu \leftarrow$ Initialize parameters for stopping condition and acceptance mechanism\;
\BlankLine
\Repeat{No $\epsilon \%$ improvement in $N$ iterations \label{algline:stop_ALNS}}{
    $\delta, r \leftarrow$ Select destroy and repair operators \label{algline:destroy_repair_1}\;
    solution$_c \leftarrow r(\delta(\text{solution}_i))$, Perform selected destroy and repair operators \label{algline:destroy_repair_2}\;
   obj$_c \leftarrow$ Evaluate solution$_c$ with algorithm~\ref{alg:obj} \label{algline:r_calc}\;
    \BlankLine
    \If{obj$_c \leq$ obj$_i$ with certain probability}{\label{algline:LS-1}
        solution$_c \leftarrow$ Implement local search algorithm (section~\ref{subsec:ls}) \label{algline:LS-2}\;
        obj$_c \leftarrow$ Evaluate solution$_c$ with algorithm~\ref{alg:obj} \label{algline:r2_calc}\;
    }
    \BlankLine
    \If{obj$_c$ - obj$_i \leq \phi$}{\label{algline:RRT}
        solution$_i$, obj$_i$ $\leftarrow$ solution$_c$, obj$_c$\;\label{algline:RRT_assignment}
        \If{obj$_c \leq$ obj$_b$}{
            solution$_b$, obj$_b$ $\leftarrow$ solution$_c$, obj$_c$\; \label{algline:best}
        }
    } 
    \BlankLine 
    $\phi \leftarrow$ max($\phi_{min}$, $\phi - \lambda$)  \label{algline:RRT2}\;
    \BlankLine
    Update destroy and repair selection probabilities ($\omega_d, \omega_r$) with roulette wheel scheme \label{algline:RW}\;
    \BlankLine
} \label{alg:ALNS}
\caption{Simulation-based Adaptive Large Neighborhood Search}
\label{alg:pseudocodeALNS}
\end{algorithm}

Initially, the incumbent solution (solution$_i$) and the best-known solution (solution$_b$) are created using the RFTS heuristic (line~\ref{algline:init_sol}). These are immediately evaluated with Algorithm~\ref{alg:obj}. 

Equal weights are then assigned to the destroy and repair operators, ensuring a uniform probability of operator selection at the start (line~\ref{algline:param_w}). The destroy and repair operators are selected based on their current weights and applied to the incumbent solution, generating a new current solution (solution$_c$) (lines~\ref{algline:destroy_repair_1}-\ref{algline:r_calc}). 
If the new solution improves upon the incumbent solution, then, with a certain probability, a local search (described in Section~\ref{subsec:ls}) is applied to further refine it (lines~\ref{algline:LS-1}--\ref{algline:r2_calc}).

The record-to-record travel acceptance criterion is then applied to determine whether the current solution should replace the incumbent. This criterion accepts solutions whose values are within a threshold ($\phi$) of that of the incumbent solution. $\phi$ decreases over time by $\lambda$, facilitating a transition from exploration to exploitation (lines~\ref{algline:RRT}-\ref{algline:RRT2}) \citep{Santini2018AMetaheuristic}.

The probability of selecting each destroy and repair operator is dynamically adjusted using a roulette wheel selection method (line~\ref{algline:RW}) \citep{Dueck1993NewTravel}. The weights of the operators are updated based on their performance, using the equation \(\omega_y = \theta\omega_y + (1-\theta)s_j\), where $y$ represents each destroy and repair operator and \(s_j\) is the performance score. This mechanism allows for a shift in focus towards more effective operators, improving the search process over time.

The algorithm continues iterating until the improvement in the solution is less than \(\epsilon\)\% over \(N\) iterations, indicating convergence towards a solution (line~\ref{algline:stop_ALNS}).

\subsubsection{Destroy and repair operators} \label{secNum} 

Our sALNS method modifies the current solution by removing delivery alternatives for some customers using four different destroy operators. Customers whose delivery alternatives are removed are added to an array \(D\).

\begin{itemize}
    \item[-] \textbf{Random destroy:} This operator removes delivery alternatives for a random subset of customers, specifically \([\kappa|C|]\) customers, where \(\kappa\) is a random variable uniformly distributed in the range \(0 < \kappa \leq 1\).
    
    \item[-] \textbf{Neighborhood destroy:} This operator selects a random customer and removes the delivery alternatives for that customer as well as all customers in its geographical proximity, determined using the K-means clustering algorithm. This approach disrupts the solution from a vehicle routing perspective, promoting exploration of alternative routing configurations.
    
    \item[-] \textbf{Worst destroy:} Customers are removed based on their impact on profit across all scenarios. Specifically, the profit loss for each customer is calculated as the additional routing cost of visiting the customer minus the delivery reward. The \([\kappa|C|]\) customers with the lowest average profit loss are added to \(D\).
    
    \item[-] \textbf{Adaptive destroy:} Customers are removed based on their historical contribution to solution improvement in past iterations. Selection probabilities are adjusted dynamically, giving higher priority to customers whose removal has led to more solution improvements compared to the others.
\end{itemize}

After altering the solution with destroy operators, the sALNS method employs repair operators to reassign delivery alternatives to all customers in \(D\). To repair the solution, we first determine the number of delivery alternatives to add, uniformly selecting a value at random. Then, we apply the following repair operators:

\begin{itemize}
    \item[-] \textbf{Random repair:} This operator randomly assigns delivery alternatives to customers in \(D\).
    
    \item[-] \textbf{High-utility slot repair:} This operator identifies delivery alternatives that offer higher utility than the opt-out alternative. For each customer in \(D\), it evaluates potential alternatives across all scenarios, ranking them based on their utility. Assignments are then made based on this ranked list.
        
    \item[-] \textbf{Greedy repair:} This strategy prioritizes delivery alternatives within the same spatial neighborhood to reduce vehicle utilization. It identifies underutilized time slots in each customer's neighborhood, selecting the lowest discount rate that exceeds the opt-out threshold for new slots.
    
    \item[-] \textbf{Two-regret repair:} This operator focuses on assigning alternatives ranked second or lower in customers' preference lists, deliberately avoiding the top alternative to explore alternative solutions. 
    
    \item[-] \textbf{Best repair:} This operator reconstructs routes for customers in \(D\) by finding best-performing positions in existing routes to decrease average costs across scenarios. The discounting mechanism aligns with the one used in the Greedy repair operator.
    
    \item[-] \textbf{Discount-rate-swap repair:} This operator retains existing time slot assignments for customers in \(D\) but randomly adjusts their discount rates, making sure the new rates differ from previous ones.
\end{itemize}

Details on the performance of these operators are provided in Section~\ref{secNum2}.

\subsubsection{Local search component} \label{subsec:ls}

While the sALNS destroy operators remove customer delivery alternatives and the repair strategies reconstruct them, the local search phase focuses on uncovering marginal enhancements to the current solution. This phase is applied probabilistically to solutions that show improvement over the incumbent solution generated by the destroy and repair operators. The local search employs various operators sequentially, terminating each operator’s execution as soon as it achieves an improvement. This early termination strategy optimizes computational efficiency while maintaining solution quality. 

The operators target specific customers or delivery alternatives, with the targets being processed in a random order. This randomization diversifies results and prevents excessive focus on any single customer or delivery alternative, which could lead to suboptimal search behavior. The following local search operators have been developed for this phase:

\begin{itemize}
    \item[-] \textbf{Random-alternative-inclusion operator:}  
    This operator identifies customers without complete slot offerings and expands their alternatives by randomly assigning an unused time slot and a discount rate. 

    \item[-] \textbf{Random-alternative-elimination operator:}  
    This operator evaluates the number of alternatives currently assigned to a customer. If the number exceeds the e-retailer's specified threshold (\(\nu\)), it randomly removes one of the existing alternatives to streamline the solution.

    \item[-] \textbf{Discount-rate-swap operator:}  
    This operator adjusts the discount rate of an assigned delivery slot to a customer. It randomly selects a new discount rate that differs from the existing rate.

    \item[-] \textbf{Common-alternative-elimination operator:}  
    This operator scans the delivery alternatives available to customers and removes any alternatives that are universally offered to all customers. By eliminating such common alternatives, the operator encourages diversity in slot assignments.

    \item[-] \textbf{High-utility-alternative-elimination operator:}  
    This operator identifies the highest-utility alternative provided to a customer and removes it from their list of delivery alternatives. This encourages the algorithm to explore alternative configurations and potentially identify better overall solutions.
\end{itemize}

By leveraging these targeted local search operators, the sALNS algorithm efficiently fine-tunes solutions, complementing the broader adjustments made during the destroy and repair phases. This phase ensures that even minor improvements are exploited, contributing to the overall quality of the solutions generated.

\section{Numerical results}
\label{secRes}

This section presents the results of our numerical experiments. We begin by describing the generation of synthetic data and problem instances in Section~\ref{secInstance}, addressing the absence of real-world data. Section~\ref{algoperformance} evaluates the algorithmic performance of our solution approach, focusing on computational efficiency and solution quality. Finally, Section~\ref{managerialinsights} analyzes several aspects of the problem, providing numerical and computational insights and exploring key implications for decision-making. 

\subsection{Instance generation}\label{secInstance}

The problem requires three primary categories of input data. The first category, routing inputs, includes customer and depot locations, travel costs ($c_{mn}$), travel times between nodes ($t_{mn}$), customer demands ($d_n$), and vehicle capacity ($Q$). To generate this data, Solomon’s CVRPTW dataset \citep{Solomon1987AlgorithmsConstraints} was adapted. Customer subsets of sizes ranging from 5 to 100 were sampled, preserving the original proportions of random, clustered, and mixed (\textit{i.e.}, combination of random and clustered points) distributions. Vehicle capacities and customer demands were scaled down by a factor of 10 to align with the scope of the TTSM-ML problem, resulting in vehicle capacities of 10 and integer demands ranging from 1 to 4. Travel times were calculated using the Euclidean distance between locations, and the delivery costs were computed proportionally to the travel times using the formula $c_{mn} = 0.4\; t_{mn}$.

The second category, delivery time slots, refers to the predefined time windows ($T := [\underline{T}_i, \overline{T}_i]$) associated with delivery options. These slots were created by analyzing the time windows in Solomon’s dataset. The maximum upper bound of all time windows was identified, and the range from 0 to this maximum was divided into three equal periods ($t_1,t_2,$ and $t_3$). This approach aligns the dataset with the model's assumptions by creating discrete delivery periods that customers can choose from.

The final category, customer sensitivities, captures preferences for delivery time slots ($\beta^{\mathrm{time}}$) and pricing ($\beta^{\mathrm{price}}$). To estimate these sensitivities, we generated five synthetic datasets, each comprising 50,000 customer choice observations, using the procedure detailed in Algorithm~\ref{alg:discretechoice}. Each observation simulates a customer selecting a delivery option from a randomly chosen subset of available alternatives (line~\ref{alglin:random_option}). The individual-level price sensitivity is sampled from a normal distribution with mean $\mu^{\text{price}}$ and standard deviation $\sigma^{\text{price}}$ (line~\ref{alglin:price_s}). Time-slot-specific preferences are likewise sampled for each available slot (line~\ref{alglin:time_s}). 
Offered prices include either no discount or a 15\% discount, each with 50\% probability (line~\ref{alglin:disc}). Utilities for each option are then computed based on these discounted prices combined with a Gumbel-distributed error term (lines~\ref{alglin:gumble}–\ref{alglin:opt-out_u}).
The option with the highest utility (including the opt-out) is recorded as the customer’s choice (lines~\ref{alglin:choice}–\ref{alglin:record}). The specific input parameters used to generate these datasets are summarized in Table~\ref{tab:beta_inputs_app} in Appendix~\ref{app:beta_inputs}.

\begin{algorithm}[]
\DontPrintSemicolon
\let\oldnl\nl%
\newcommand{\nonl}{\renewcommand{\nl}{\let\nl\oldnl}}
\BlankLine
\SetKwComment{Comment}{$\triangleright$\ }{}
\BlankLine
\nonl \ForEach{observation}{
    Randomly assign time-slot availability pattern from set of all possible combinations\; \label{alglin:random_option}
    $\beta^{\mathrm{price}} \sim \mathcal{N}(\mu^{price},\sigma^{price})$\; \label{alglin:price_s}
    \BlankLine
    \nonl \ForEach{available time-slot}{
        $\beta^{\mathrm{time}}_t  \sim \mathcal{N}(\mu^{time}_t,\sigma^{time}_t)$\; \label{alglin:time_s}
        $h_t \sim$ DiscreteUniform(discountSet) \;\label{alglin:disc}
    }
    \BlankLine

    \nonl \ForEach{time-slot and opt\text{-}out (0)}{
        $\xi_t \sim EV1(1,0)$\;\label{alglin:gumble}
    }
    \BlankLine
    
    \nonl \ForEach{time-slot $t$}{
        \eIf{time-slot $t$ is available}{
            $U_t = \beta^{\mathrm{time}}_t + \beta^{\mathrm{price}} \cdot f \cdot h_t + \xi_t$\;
        }{
            $U_t = -\infty$\;
        }
    }
    \BlankLine
    $U_{opt\text{-}out} = \xi_0$\; \label{alglin:opt-out_u}
    \BlankLine
    $choice = \text{argmax}_t\{U_t, U_{opt\text{-}out}\}$\; \label{alglin:choice}
    \BlankLine
    Store time-slot's availability, prices, and choice\; \label{alglin:record}
}
\caption{Dataset generation for customers' behavioral model estimation}
\label{alg:discretechoice}
\end{algorithm}

The five datasets vary in time-slot preferences and price sensitivity (see Table~\ref{tab:beta_inputs_app} in Appendix~\ref{app:beta_inputs}), enabling a robustness check of model performance across heterogeneous customer populations. Tables~\ref{tab:beta_mnl_ml_all1} and~\ref{tab:beta_mnl_ml_all2} in Appendix~\ref{app:ml_mnl_tables} present the parameter estimates and log-likelihood values for both the MNL and ML models fitted to the datasets. The results show that all datasets exhibit customer preference heterogeneity, with the ML model consistently outperforming the MNL model, as confirmed by statistically significant likelihood-ratio tests. This makes the datasets well-suited for testing the TTSM-ML, which is designed to capture individual-level preference variation.\footnote{For datasets where the ML model does not outperform simpler alternatives such as the MNL, we advocate the use of frameworks specifically designed for those choice structures, offering greater efficiency and parsimony without sacrificing modeling fidelity.}

In our experiments, these estimated coefficients are used linearly in the utility functions for delivery slot choices, consistent with standard discrete choice modeling practice \citep{Sifringer2020EnhancingLearning}. Under both MNL and ML models, the utility is defined as:
\begin{align}
&\textit{\footnotesize MNL:} 
&u_{inr} &= \beta^{\mathrm{time}}_{i} \cdot t_i + \beta^{\mathrm{price}} \cdot h_{in} \cdot f + \bar{\xi}_{inr} - M_{nr} \cdot (1-\gamma_{in}), &\quad \forall i \in I, ~\forall n \in C, ~\forall r \in R. \notag \\
&\textit{\footnotesize ML:} 
&u_{inr} &= \beta^{\mathrm{time}}_{i} \cdot t_i + \bar{\beta}^{\mathrm{price}}_{r} \cdot h_{in} \cdot f + \bar{\xi}_{inr} - M_{nr} \cdot (1-\gamma_{in}), &\quad \forall i \in I, ~\forall n \in C, ~\forall r \in R. \notag
\end{align}

where $\bar{\beta}^{\mathrm{price}}_{r}$ is a scenario-specific draw from a normally distributed parameter calibrated to reflect price sensitivity heterogeneity.

\subsection{Algorithmic performance}
\label{algoperformance}

This section evaluates the performance of the sALNS method. We first calibrate the number of scenarios ($R$) in Section~\ref{secSupSize} to ensure a reliable approximation of the stochastic components captured in our SAA-based MILP formulation. We then assess the computational performance of sALNS, including heuristic configuration and operator analysis, across varying instance sizes in Section~\ref{secNum2}.

All experiments were run on a 10‑core AMD EPYC processor with 32GB of RAM. Models were implemented in Python 3.11.6 and solved in Gurobi 11.0.3 with a 12‑hour time limit per instance, following \cite{Mackert2019Choice-basedDelivery}.

\subsubsection{Number of scenarios} \label{secSupSize}

This section identifies the smallest number of scenarios, $R$, needed to ensure a reliable approximation of the underlying stochastic model within our SAA-based MILP formulation. To that end, we conduct both in-sample and out-of-sample analyses using 5-customer instances, where exact optimization is tractable. Tests are performed across the three map configurations and the five synthetic datasets, providing diverse instances in terms of customer locations and preference configurations.

In the in-sample test, we examined how the objective values evolve as the number of scenarios ($R$) increases. For each instance and scenario size, we generated 20 fresh sets of ($\bar{\beta}^{\mathrm{price}}_{r}, \bar{\xi}_{inr}$) vectors and solved the corresponding TTSM-ML models. To enable consistent comparison, we computed the relative deviation of each solution’s objective value from a benchmark, defined as the best objective value achieved for the same instance using 120 scenarios. Figure~\ref{fig:In_Sample_Test} presents the average percentage deviations across all instances, with error bars indicating the standard deviation of these deviations. Each sub-figure corresponds to a dataset and contains three line charts, each representing the results of a different map configuration. As shown, the deviations decrease with increasing $R$, stabilizing around 2\% at $R=100$ for all instances. Moreover, the error bars also converge to less than 2\% beyond 100 scenarios.

\begin{figure}[h!]
\centering
\scriptsize
\begin{subfigure}{0.32\textwidth}
\captionsetup{justification=centering}
\footnotesize
    \tikz{
    \begin{axis}[
  width=\linewidth,
      ytick={0,10,...,90}, ytick align=outside, ytick style={draw=none},
      xtick={0,20,...,100,120}, xtick align=outside, xtick style={draw=none},
      xlabel=Scenarios,
      ylabel={In-sample deviations (\%)},
      height=\axisdefaultheight,
      grid=both,
  legend to name=sharedlegend
      ]
    \addplot+[
      BlueViolet, 
      mark=triangle,
      mark options={BlueViolet, scale=0.8},
      smooth, 
      error bars/.cd, 
        y fixed,
        y dir=both, 
        y explicit
    ] 
    table [x=x, y=y,y error=error] {
        x  y       error
        5	51.66583324	31.24310175
        10	43.99666301	27.84299022
        15	33.03493405	23.0759637
        20	24.04157194	17.80203288
        30	18.64466344	16.59475433
        40	15.86558757	14.13804263
        50	14.02671298	12.02087659
        60	9.000536967	7.657332086
        80	5.288945856	4.085866868
        100	5.206920377	3.129444622
        110	4.384292543	2.829345925
    };
    \addlegendentry{Random map}
    \addplot+[
      purple, 
      mark=square,
      mark options={purple, scale=0.8},
      smooth, 
      error bars/.cd, 
        y fixed,
        y dir=both, 
        y explicit
    ] 
    table [x=x, y=y,y error=error] {
        x  y       error
        5	20.88504527	7.665189462
        10	15.51401384	10.39703502
        15	10.87505505	8.974544229
        20	8.043672583	7.484088072
        30	7.002135609	5.472276618
        40	7.677046298	5.58786016
        50	7.49078385	4.503620129
        60	4.920206234	3.442534964
        80	2.49790465	2.169872162
        100	1.236548825	1.007605471
        110	1.269458938	0.71678
    };
    \addlegendentry{Clustered map}
    \addplot+[
      teal, 
      mark=*,
      mark options={teal, scale=0.8},
      smooth, 
      error bars/.cd, 
        y fixed,
        y dir=both, 
        y explicit
    ] 
    table [x=x, y=y,y error=error] {
        x  y       error
        5	47.17009051	15.78357616
        10	30.32740055	17.73685329
        15	26.37704505	18.54802067
        20	21.49290959	15.64492865
        30	11.88590151	8.111031233
        40	10.48271228	7.218666173
        50	7.895725539	5.721409709
        60	7.151233172	3.722600849
        80	3.314341988	3.736235906
        100	2.428428978	2.258552367
        110	1.408388471	1.372343957
    };
    \addlegendentry{Mixed map}
    \end{axis}
    }
\caption{Dataset 1}
\label{fig:Out-of-sample-Test-10}
\end{subfigure}
\hfill
\begin{subfigure}{0.32\textwidth}
\captionsetup{justification=centering}
\footnotesize
\tikz{
\begin{axis}[
  width=\linewidth,
  ytick={0,10,...,100}, ytick align=outside, ytick style={draw=none},
  xtick={0,20,...,120}, xtick align=outside, xtick style={draw=none},
  xlabel=Scenarios,
  ylabel={In-sample deviations (\%)},
  height=\axisdefaultheight,
  grid=both
  ]
    \addplot+[
      BlueViolet, 
      mark=triangle,
      mark options={BlueViolet, scale=0.8},
      smooth,
      error bars/.cd, 
        y fixed,
        y dir=both, 
        y explicit
    ] 
    table [x=x, y=y,y error=error] {
        x  y       error
        5	32.62347265	16.66872728
        10	22.99584067	15.39010835
        15	21.90097758	12.95762686
        20	20.26517613	12.66034424
        30	15.16462267	11.07500938
        40	14.96549332	8.528250165
        50	10.99233335	7.753232992
        60	7.864453294	6.487018103
        80	6.019656963	4.16990685
        100	3.113778224	2.329033276
        110	1.893340919	1.213499049
    };
    \addplot+[
      purple, 
      mark=square,
      mark options={purple, scale=0.8},
      smooth, 
      error bars/.cd, 
        y fixed,
        y dir=both, 
        y explicit
    ] 
    table [x=x, y=y,y error=error] {
        x  y       error
        5	34.01996159	19.85012348
        10	13.76904803	9.749780051
        15	12.78371215	9.727577384
        20	9.227403377	8.96906335
        30	8.270514662	5.66612804
        40	6.913353571	4.752889139
        50	5.917531653	4.001608305
        60	4.376031815	2.978129238
        80	1.860722169	1.618345548
        100	1.371629688	1.106764355
        110	1.230002658	1.021109334
    };
    \addplot+[
      teal, 
      mark = *,
      mark options={teal, scale=0.8},
      smooth, 
      error bars/.cd, 
        y fixed,
        y dir=both, 
        y explicit
    ] 
    table [x=x, y=y,y error=error] {
        x  y       error
        5	63.14703246	27.42731243
        10	33.40398541	24.02801944
        15	26.23912942	21.16152001
        20	13.15081964	13.53721079
        30	12.28858984	12.6664624
        40	10.16593025	6.957737858
        50	8.146322055	8.311320485
        60	11.2237023	6.644825904
        80	5.462706006	4.883488691
        100	3.61700807	2.573356607
        110	2.465894427	2.079570187
    };
    \end{axis}
    }
\caption{Dataset 2}
\label{fig:Out-of-sample-Test-20}
\end{subfigure}
\hfill
\begin{subfigure}{0.32\textwidth}
\captionsetup{justification=centering}
\footnotesize
\tikz{
\begin{axis}[
  width=\linewidth,
  ytick={0,5,...,45}, ytick align=outside, ytick style={draw=none},
  xtick={0,20,...,120}, xtick align=outside, xtick style={draw=none},
  xlabel=Scenarios,
  ylabel={In-sample deviations (\%)},
  height=\axisdefaultheight,
  grid=both
  ]
    \addplot+[
      BlueViolet, 
      mark=triangle,
      mark options={BlueViolet, scale=0.8},
      smooth,
      error bars/.cd, 
        y fixed,
        y dir=both, 
        y explicit
    ] 
    table [x=x, y=y,y error=error] {
        x  y       error
        5	8.342119534	3.938403022
        10	4.83351727	3.240548972
        15	3.17325797	1.806509442
        20	2.931933499	1.424348817
        30	2.323193143	1.605310438
        40	1.619029382	1.2412661
        50	1.928012199	0.975460286
        60	1.781731474	0.867992335
        80	1.605221766	1.459827787
        100	1.166252705	0.88460159
        110	1.222505575	1.045603727

    };
    \addplot+[
      purple, 
      mark=square,
      mark options={purple, scale=0.8},
      smooth, 
      error bars/.cd, 
        y fixed,
        y dir=both, 
        y explicit
    ] 
    table [x=x, y=y,y error=error] {
        x  y       error
        5	5.410147262	4.556179547
        10	3.795779224	2.801616172
        15	2.704174145	1.645990579
        20	2.284477134	2.128120199
        30	1.642431223	1.404667104
        40	1.255446239	1.024562935
        50	2.375653447	0.882697119
        60	1.348353731	0.773038626
        80	1.363011444	1.417708663
        100	1.478188601	1.108904903
        110	1.011327017	1.095469096
    };
    \addplot+[
      teal, 
      mark = *,
      mark options={teal, scale=0.8},
      smooth, 
      error bars/.cd, 
        y fixed,
        y dir=both, 
        y explicit
    ] 
    table [x=x, y=y,y error=error] {
        x  y       error
        5	20.81231357	13.96832214
        10	11.93539506	7.966966468
        15	10.10717613	6.129266408
        20	10.8096619	6.359633194
        30	7.106607108	4.863218381
        40	4.112635701	3.521202978
        50	3.665514143	3.186265527
        60	4.279650282	3.551002098
        80	2.281026945	1.962364368
        100	1.656549223	1.001137177
        110	1.515800957	0.765270299 
    };
    \end{axis}
    }
\caption{Dataset 3}
\label{fig:Out-of-sample-Test-30}
\end{subfigure}
\hfill
\begin{subfigure}{0.32\textwidth}
\captionsetup{justification=centering}
\footnotesize
\tikz{
\begin{axis}[
  width=\linewidth,
  ytick={0,5,...,35,40}, ytick align=outside, ytick style={draw=none},
  xtick={0,20,...,120}, xtick align=outside, xtick style={draw=none},
  xlabel=Scenarios,
  ylabel={In-sample deviations (\%)},
  height=\axisdefaultheight,
  grid=both
  ]
    \addplot+[
      BlueViolet, 
      mark=triangle,
      mark options={BlueViolet, scale=0.8},
      smooth,
      error bars/.cd, 
        y fixed,
        y dir=both, 
        y explicit
    ] 
    table [x=x, y=y,y error=error] {
        x  y       error
        5	21.75863301	10.93891497
        10	15.49936836	10.54977769
        15	10.9756661	5.302124318
        20	6.444810545	5.602649905
        30	5.909402584	4.980417623
        40	8.332031362	5.458636896
        50	5.837765019	3.950786889
        60	2.718668969	3.943286179
        80	1.178384524	1.072379727
        100	1.144033746	1.474422097
        110	1.224508102	0.276081551
    };
    \addplot+[
      purple, 
      mark=square,
      mark options={purple, scale=0.8},
      smooth, 
      error bars/.cd, 
        y fixed,
        y dir=both, 
        y explicit
    ] 
    table [x=x, y=y,y error=error] {
        x  y       error
        5	24.2293603	13.47832729
        10	17.1286383	10.97090089
        15	9.190324628	6.134386227
        20	9.766408889	5.006460212
        30	5.958964643	4.084042556
        40	8.676196774	2.599527541
        50	4.953350538	2.610982562
        60	4.837498017	1.543510001
        80	1.679017859	2.347054359
        100	1.446180596	0.013760242
        110	1.690276328	0.712992061
    };
    \addplot+[
      teal, 
      mark = *,
      mark options={teal, scale=0.8},
      smooth, 
      error bars/.cd, 
        y fixed,
        y dir=both, 
        y explicit
    ] 
    table [x=x, y=y,y error=error] {
        x  y       error
        5	13.79157065	7.417851811
        10	8.393938755	4.723194519
        15	6.072875702	2.910120101
        20	6.261434987	3.101212803
        30	4.776185032	2.310763038
        40	3.319299109	3.061215578
        50	2.865550439	3.134026914
        60	1.504052052	2.383943967
        80	0.107473097	1.788518344
        100	0.58546214	1.371354603
        110	0.591503178	0.180902394
    };
    \end{axis}
    }
\caption{Dataset 4}
\label{fig:Out-of-sample-Test-40}
\end{subfigure}
\hfill
\begin{subfigure}{0.32\textwidth}
\captionsetup{justification=centering}
\footnotesize
\tikz{
\begin{axis}[
  width=\linewidth,
  ytick={0,5,...,35}, ytick align=outside, ytick style={draw=none},
  xtick={0,20,...,120}, xtick align=outside, xtick style={draw=none},
  xlabel=Scenarios,
  ylabel={In-sample deviations (\%)},
  height=\axisdefaultheight,
  grid=both
  ]
    \addplot+[
      BlueViolet, 
      mark=triangle,
      mark options={BlueViolet, scale=0.8},
      smooth,
      error bars/.cd, 
        y fixed,
        y dir=both, 
        y explicit
    ] 
    table [x=x, y=y,y error=error] {
        x  y       error
        5	16.95353763	10.46961644
        10	9.170832258	4.082110248
        15	8.960197545	4.375487531
        20	7.419131362	6.863285669
        30	7.398044745	2.125818532
        40	7.08720256	4.511113512
        50	3.696688623	3.625612698
        60	4.153036686	3.165225335
        80	3.00825832	0.766760374
        100	1.888887462	0.239100398
        110	1.335741583	0.042912509
    };
    \addplot+[
      purple, 
      mark=square,
      mark options={purple, scale=0.8},
      smooth, 
      error bars/.cd, 
        y fixed,
        y dir=both, 
        y explicit
    ] 
    table [x=x, y=y,y error=error] {
        x  y       error
        5	21.08992902	12.27598974
        10	12.3847875	5.507816841
        15	12.02205605	6.297067538
        20	8.237117052	5.415400149
        30	7.802017701	3.755355635
        40	5.388822817	2.674821933
        50	7.380690906	4.715671318
        60	6.067203101	3.69754386
        80	2.149889013	3.525156494
        100	1.681940484	1.368495866
        110	0.800832409	2.014446806
    };
    \addplot+[
      teal, 
      mark = *,
      mark options={teal, scale=0.8},
      smooth, 
      error bars/.cd, 
        y fixed,
        y dir=both, 
        y explicit
    ] 
    table [x=x, y=y,y error=error] {
        x  y       error
        5	15.64836645	10.62089163
        10	8.706113748	8.474347271
        15	5.166664172	4.277001649
        20	5.902275562	5.936872596
        30	4.755275211	3.96420892
        40	5.02642709	1.453948535
        50	4.190568935	3.235124854
        60	3.7867185	2.11348733
        80	1.859045947	1.461338912
        100	1.819547665	1.889783075
        110	0.710603357	2.062567308
    };
    \end{axis}
    }
\caption{Dataset 5}
\label{fig:Out-of-sample-Test-50}
\end{subfigure}
\hfill
\begin{subfigure}{0.32\textwidth}
\centering
\begin{tikzpicture}
\node at (0,0) {\ref{sharedlegend}}; 
\end{tikzpicture}
\end{subfigure}
\caption{In-sample test: It evaluates the deviation in objective values across problems with varying scenario sizes, ranging from 5 to 110 scenarios, compared to the problem with a fixed number of 120 scenarios.}
\label{fig:In_Sample_Test}
\end{figure}

To assess the generalization performance of our solutions under stochastic customer preferences, we conduct a weak version of out-of-sample test. For each map and dataset, we generate 20 independent sets of simulated preference scenarios (\textit{i.e.}, realizations of $\bar{\beta}^{\mathrm{price}}_{r}$ and $\bar{\xi}_{inr}$). The TTSM-ML model is then solved using one of these sets, referred to as the training scenario. The resulting solution is subsequently evaluated on the remaining 19 sets, referred to as the test scenarios, with the slot and price decisions held fixed.

The term optimal profit refers to the objective value achieved on the training scenario used for optimization, while realized profits denote the profit values obtained when the same solution is evaluated on the test scenarios. We compute the average percentage deviation between realized and optimal profits, measured relative to the optimal profit, to assess each solution’s robustness to preference variation. Figure~\ref{fig:Out-of-sample-Test} presents the average deviation and standard deviation across these comparisons.

\begin{figure}[h!]
\centering
\scriptsize
\begin{subfigure}{0.32\textwidth}
\captionsetup{justification=centering}
\footnotesize
    \tikz{
    \begin{axis}[
  width=\linewidth,
      ytick={0,2,4.00,...,28}, ytick align=outside, ytick style={draw=none},
      xtick={0,20,...,100,120}, xtick align=outside, xtick style={draw=none},
      xlabel=Scenarios,
      ylabel={Out-of-sample deviations (\%)},
      height=\axisdefaultheight,
      grid=both,
  legend to name=sharedlegend
      ]
    \addplot+[
      BlueViolet, 
      mark=triangle,
      mark options={BlueViolet, scale=0.8},
      smooth,
      error bars/.cd, 
        y fixed,
        y dir=both, 
        y explicit
    ] 
    table [x=x, y=y,y error=error] {
        x  y       error
        5	7.072009029	3.650150038
        10	3.884200305	2.465341121
        15	3.353445603	2.082022102
        20	2.825408014	1.562308061
        30	1.595238435	0.997679326
        40	1.115177365	0.832107018
        50	1.137757577	0.669315244
        60	0.770928781	0.558536938
        80	0.79256139	0.579684265
        100	0.653254312	0.416967528
        110	0.586295941	0.384343443
        120	0.616859429	0.40378436
    };
    \addlegendentry{Random map}
    \addplot+[
      purple, 
      mark=square,
      mark options={purple, scale=0.8},
      smooth, 
      error bars/.cd, 
        y fixed,
        y dir=both, 
        y explicit
    ] 
    table [x=x, y=y,y error=error] {
        x  y       error
        5	8.440372505	5.018555281
        10	3.699300442	3.146289141
        15	3.946158018	2.41019968
        20	3.948074942	2.364910016
        30	3.902910613	2.426000445
        40	3.64444756	2.460854022
        50	2.481927175	1.767971841
        60	2.085285293	1.259131107
        80	2.037166885	1.280600305
        100	1.93534625	1.177913818
        110	1.603709348	0.976230094
        120	1.668100872	1.023966876
    };
    \addlegendentry{Clustered map}
    \addplot+[
      teal, 
      mark = *,
      mark options={teal, scale=0.8},
      smooth, 
      error bars/.cd, 
        y fixed,
        y dir=both, 
        y explicit
    ] 
    table [x=x, y=y,y error=error] {
        x  y       error
        5	7.61967888	6.089435112
        10	6.540947072	3.969511011
        15	4.409152814	2.809999266
        20	4.165657473	2.890636403
        30	3.890312679	2.174218158
        40	3.538242616	1.924723863
        50	2.303060625	1.51271428
        60	2.341793859	1.402799107
        80	2.517159279	1.485360452
        100	1.62534734	1.110406168
        110	1.651075206	1.09910445
        120	1.415300342	0.847967613
    };
    \addlegendentry{Mixed map}
    \end{axis}
    }
\caption{Dataset 1}
\label{fig:Out-of-sample-Test-10}
\end{subfigure}
\hfill
\begin{subfigure}{0.32\textwidth}
\captionsetup{justification=centering}
\footnotesize
\tikz{
\begin{axis}[
  width=\linewidth,
  ytick={0,5,...,50}, ytick align=outside, ytick style={draw=none},
  xtick={0,20,...,120}, xtick align=outside, xtick style={draw=none},
  xlabel=Scenarios,
  ylabel={Out-of-sample deviations (\%)},
  height=\axisdefaultheight,
  grid=both
  ]
    \addplot+[
      BlueViolet, 
      mark=triangle,
      mark options={BlueViolet, scale=0.8},
      smooth,
      error bars/.cd, 
        y fixed,
        y dir=both, 
        y explicit
    ] 
    table [x=x, y=y,y error=error] {
        x  y       error
        5	9.114420759	5.183449092
        10	2.262996122	1.696703307
        15	3.518212808	2.482370685
        20	1.277747152	1.000607006
        30	1.570933022	0.962052665
        40	0.831363877	0.675805317
        50	1.484347595	0.822934014
        60	1.233274836	0.542896197
        80	0.57680221	0.412391108
        100	0.417600448	0.301720095
        110	0.394037997	0.280611527
        120	0.946310858	0.569297824
    };
    \addplot+[
      purple, 
      mark=square,
      mark options={purple, scale=0.8},
      smooth, 
      error bars/.cd, 
        y fixed,
        y dir=both, 
        y explicit
    ] 
    table [x=x, y=y,y error=error] {
        x  y       error
        5	5.501872896	2.994966643
        10	2.501553651	1.486785699
        15	3.371440772	1.901378059
        20	2.963053961	1.419216666
        30	1.534367982	1.019963852
        40	1.319127334	0.815351046
        50	0.975852884	0.684887409
        60	0.917373527	0.581683312
        80	0.709735686	0.524675616
        100	0.486892559	0.359048607
        110	0.571538491	0.369695861
        120	0.570600165	0.346073284
    };
    \addplot+[
      teal, 
      mark = *,
      mark options={teal, scale=0.8},
      smooth, 
      error bars/.cd, 
        y fixed,
        y dir=both, 
        y explicit
    ] 
    table [x=x, y=y,y error=error] {
        x  y       error
        5	10.41557707	6.387542356
        10	7.452290773	3.896795878
        15	5.752700334	2.845708155
        20	3.938980954	2.253917732
        30	3.159932734	2.034784452
        40	2.811896312	1.715257191
        50	2.762554265	1.594680893
        60	3.079178038	1.584170768
        80	2.289487095	1.314970416
        100	1.726946601	1.013836355
        110	1.66373835	0.975985004
        120	1.650445033	1.008997356
    };
    \end{axis}
    }
\caption{Dataset 2}
\label{fig:Out-of-sample-Test-10}
\end{subfigure}
\hfill
\begin{subfigure}{0.32\textwidth}
\captionsetup{justification=centering}
\footnotesize
\tikz{
\begin{axis}[
  width=\linewidth,
  ytick={0,5,...,45}, ytick align=outside, ytick style={draw=none},
  xtick={0,20,...,120}, xtick align=outside, xtick style={draw=none},
  xlabel=Scenarios,
  ylabel={Out-of-sample deviations (\%)},
  height=\axisdefaultheight,
  grid=both
  ]
    \addplot+[
      BlueViolet, 
      mark=triangle,
      mark options={BlueViolet, scale=0.8},
      smooth,
      error bars/.cd, 
        y fixed,
        y dir=both, 
        y explicit
    ] 
    table [x=x, y=y,y error=error] {
        x  y       error
        5	5.66029	3.9210486
        10	3.6513239	3.3680067
        15	2.9792186	1.9021851
        20	3.9874723	2.6770019
        30	3.2560011	1.8715427
        40	2.9317921	1.4506696
        50	2.4479206	1.2981638
        60	2.5227064	1.2254376
        80	1.7260742	1.1938341
        100	1.8677303	1.1061565
        110	1.401982774	0.79753555
        120	1.326622013	0.721042782
    };
    \addplot+[
      purple, 
      mark=square,
      mark options={purple, scale=0.8},
      smooth, 
      error bars/.cd, 
        y fixed,
        y dir=both, 
        y explicit
    ] 
    table [x=x, y=y,y error=error] {
        x  y       error
        5	5.1627668	2.9347126
        10	3.2212042	1.8719903
        15	2.5421512	1.5170376
        20	2.2145765	1.2860289
        30	2.0193173	1.2487333
        40	2.4713522	1.4445329
        50	1.9684625	1.0328478
        60	2.1493202	1.0803944
        80	1.0587914	0.6524256
        100	0.984995	0.7527927
        110	0.903049573	0.652372271
        120	0.699103731	0.577738772
    };
    \addplot+[
      teal, 
      mark = *,
      mark options={teal, scale=0.8},
      smooth, 
      error bars/.cd, 
        y fixed,
        y dir=both, 
        y explicit
    ] 
    table [x=x, y=y,y error=error] {
        x  y       error
        5	6.7871971	3.2158015
        10	5.4230325	2.8142322
        15	2.717929	1.8630837
        20	2.6175946	1.8509486
        30	1.4792142	0.8856115
        40	1.7094583	1.1007223
        50	1.2223926	0.8452679
        60	1.1341807	0.8882537
        80	0.9786349	0.6788966
        100	1.0150122	0.7060217
        110	0.835797541	0.516903674
        120	0.871187709	0.628252975
    };
    \end{axis}
    }
\caption{Dataset 3}
\label{fig:Out-of-sample-Test-10}
\end{subfigure}
\hfill
\begin{subfigure}{0.32\textwidth}
\captionsetup{justification=centering}
\footnotesize
\tikz{
\begin{axis}[
  width=\linewidth,
  ytick={0,5,...,30}, ytick align=outside, ytick style={draw=none},
  xtick={0,20,...,120}, xtick align=outside, xtick style={draw=none},
  xlabel=Scenarios,
  ylabel={Out-of-sample deviations (\%)},
  height=\axisdefaultheight,
  grid=both
  ]
    \addplot+[
      BlueViolet, 
      mark=triangle,
      mark options={BlueViolet, scale=0.8},
      smooth,
      error bars/.cd, 
        y fixed,
        y dir=both, 
        y explicit
    ] 
    table [x=x, y=y,y error=error] {
        x  y       error
        5	9.445931045	6.82587437
        10	4.085677105	2.974421875
        15	3.86680703	2.782176255
        20	4.276201065	3.007355175
        30	3.16883435	2.41212476
        40	2.04963251	1.620106675
        50	1.77868044	1.57752831
        60	2.47193092	1.685749875
        80	2.379754235	1.141922865
        100	2.11487529	1.143661115
        110	2.057431245	1.175885555
        120	1.53204177	0.89055169
    };
    \addplot+[
      purple, 
      mark=square,
      mark options={purple, scale=0.8},
      smooth, 
      error bars/.cd, 
        y fixed,
        y dir=both, 
        y explicit
    ] 
    table [x=x, y=y,y error=error] {
        x  y       error
        5	8.845449	5.227205603
        10	7.330357988	5.896879884
        15	3.18358652	2.150307104
        20	3.344311936	2.047057892
        30	2.241026436	1.484623656
        40	2.360051332	1.51224142
        50	1.399817312	0.824654308
        60	1.417758528	1.069555824
        80	1.088774196	0.660114056
        100	1.181843412	0.873277472
        110	1.131529724	0.796762548
        120	0.825258612	0.57595688
    };
    \addplot+[
      teal, 
      mark = *,
      mark options={teal, scale=0.8},
      smooth, 
      error bars/.cd, 
        y fixed,
        y dir=both, 
        y explicit
    ] 
    table [x=x, y=y,y error=error] {
        x  y       error
        5	10.09134028	7.70706658
        10	5.73434122	3.76021356
        15	6.939662336	5.513552452
        20	4.116095412	2.986554492
        30	4.827449732	3.209989664
        40	2.75816604	2.113339416
        50	1.444662876	1.317070076
        60	0.565216406	0.603928191
        80	1.112849798	0.762102656
        100	2.012487566	1.68361554
        110	1.287683584	1.254592063
        120	1.666735288	1.280770334
    };
    \end{axis}
    }
\caption{Dataset 4}
\label{fig:Out-of-sample-Test-10}
\end{subfigure}
\hfill
\begin{subfigure}{0.32\textwidth}
\captionsetup{justification=centering}
\footnotesize
\tikz{
\begin{axis}[
  width=\linewidth,
  ytick={0,5,...,45}, ytick align=outside, ytick style={draw=none},
  xtick={0,20,...,120}, xtick align=outside, xtick style={draw=none},
  xlabel=Scenarios,
  ylabel={Out-of-sample deviations (\%)},
  height=\axisdefaultheight,
  grid=both
  ]
    \addplot+[
      BlueViolet, 
      mark=triangle,
      mark options={BlueViolet, scale=0.8},
      smooth,
      error bars/.cd, 
        y fixed,
        y dir=both, 
        y explicit
    ] 
    table [x=x, y=y,y error=error] {
        x  y       error
        5	14.22047804	10.64743709
        10	15.58443708	8.203903584
        15	7.598447407	5.328921924
        20	9.951151307	4.852038469
        30	5.44639592	4.092206089
        40	3.00405581	2.390466326
        50	4.604455951	2.95480345
        60	5.208797855	3.068227506
        80	3.071090988	2.236238247
        100	2.793339757	2.013328437
        110	2.798586099	1.776854669
        120	2.65107939	1.690074413
    };
    \addplot+[
      purple, 
      mark=square,
      mark options={purple, scale=0.8},
      smooth, 
      error bars/.cd, 
        y fixed,
        y dir=both, 
        y explicit
    ] 
    table [x=x, y=y,y error=error] {
        x  y       error
        5	22.27818987	16.64591194
        10	11.00868761	8.287289888
        15	6.418877121	4.238878855
        20	6.944136478	6.007042202
        30	2.591192565	1.697511762
        40	3.491222312	2.143729683
        50	2.634175646	1.52809306
        60	2.125395208	1.176152223
        80	1.983774859	1.380620951
        100	1.158198226	0.942980594
        110	1.469411707	1.057738052
        120	1.323459201	1.105659459
    };
    \addplot+[
      teal, 
      mark = *,
      mark options={teal, scale=0.8},
      smooth, 
      error bars/.cd, 
        y fixed,
        y dir=both, 
        y explicit
    ] 
    table [x=x, y=y,y error=error] {
        x  y       error
        5	12.7618001	9.442415
        10	6.162972181	4.344198762
        15	7.105721591	3.715210267
        20	4.499876511	2.29738458
        30	2.513245128	0.938707593
        40	2.103092102	1.257173393
        50	3.222397863	1.614890866
        60	3.465838571	1.076162885
        80	2.375338367	1.053314609
        100	1.080933729	0.437432339
        110	1.67398491	0.330567491
        120	1.38583917	0.448263238
    };
    \end{axis}
    }
\caption{Dataset 5}
\label{fig:Out-of-sample-Test-10}
\end{subfigure}
\hfill
\begin{subfigure}{0.32\textwidth}
\centering
\vspace{-2cm} 
\ref{sharedlegend}
\end{subfigure}

\caption{Out-of-sample test (weak version): It measures the deviation between the stochastic optimal objective value and the optimal objective value obtained by re-evaluating the solution with different $\xi$ vector realizations.}
\label{fig:Out-of-sample-Test}
\end{figure}

Results show that once 100 scenarios are used, deviations remain below 2\%, with standard deviations stabilizing around 2\%. These findings confirm that using 100 scenarios provides reliable and stable objective values across different instances.

We also conducted a similar scenario count analysis on 10-customer instances. However, due to the increased problem size, exact optimization is no longer tractable. Instead, we applied the sALNS heuristic for both in-sample and out-of-sample evaluations. Since the validity and performance of sALNS are demonstrated in the following section, we refer the reader to Appendix~\ref{app:Scenario10} for detailed results, which confirm that the stabilization behavior observed at 100 scenarios holds for larger instances as well. Accordingly, we use $R = 100$ in all subsequent experiments unless otherwise stated.

\subsubsection{Performance evaluation of the sALNS method}\label{secNum2}

This section focuses on evaluating the performance of the sALNS method. We begin by clarifying two key algorithmic design decisions: the choice of routing heuristic and the number of algorithm repetitions.

We evaluated three routing heuristics: CFRS (Cluster First, Route Second), CW (Clarke and Wright’s savings heuristic), and ICW (Improved Clarke and Wright’s savings heuristic). A detailed comparison is provided in Appendix~\ref{subsec:route}. Among the three, the CW heuristic offered the best trade-off between solution quality and computational efficiency, and was therefore selected for integration into our sALNS method.

Regarding the determination of the algorithmic repetition count, a preliminary variability assessment over ten repetitions per instance showed a coefficient of variation of about 2\% and a standard deviation of about 1\%, indicating that ten repetitions are sufficient to obtain stable and representative results without incurring excessive computational costs. {To ensure a fair comparison with the Gurobi solver, all sALNS solutions are held fixed while routing costs are re-evaluated by solving the induced VRPTW exactly in Gurobi. This provides consistent objective values and optimality gaps.

We start with an evaluation of sALNS on five‑customer instances, enabling benchmarking against the optimal solutions produced by the Gurobi solver. Table~\ref{tab:small-size_res} compares the solve times of Gurobi and sALNS, and reports sALNS’s best and mean optimality gaps relative to the optimal solutions found by Gurobi. The results demonstrate the limited scalability of exact methods like Gurobi as the number of scenarios increases. Although Gurobi guarantees optimality, its solve time becomes prohibitive, exceeding four hours for 100 scenarios, even for a small five-customer instance. In contrast, sALNS handles the increasing scenario complexity efficiently, solving the same problems in under 13 seconds while maintaining optimality gaps below 1.5\%. This makes sALNS a more practical choice for larger and more realistic instances.

\begin{table}[h!]
    \centering
    \scriptsize
    \caption{Comparison of computational performance for instances with five customers: solving times for Gurobi and sALNS, and optimality gaps for sALNS.}
    \label{tab:small-size_res}
    \begin{tabular}{cc@{\hskip 10pt}cccc}
    \toprule
    \multirow{2}{*}{Scenarios} & Gurobi & & \multicolumn{3}{c}{sALNS} \\ \cline{2-2} \cline{4-6}
    & Time (s)  & & Time (s) & Best Gap (\%) & Mean Gap (\%) \\ \midrule
    5   & 9     & & 1.16 & 0.81 & 1.10 \\
    10  & 30    & & 1.84 & 0.93 & 1.05 \\
    15  & 84    & & 2.38 & 1.12 & 1.16 \\
    20  & 203   & & 2.88 & 0.98 & 1.05 \\
    30  & 3494  & & 4.11 & 1.04 & 1.09 \\
    40  & 5662  & & 5.35 & 1.15 & 1.19 \\
    50  & 9307  & & 7.30 & 1.10 & 1.24 \\
    60  & 10152 & & 8.38 & 1.11 & 1.15 \\
    80  & 16824 & & 10.70& 0.80 & 0.95 \\
    100 & 18784 & & 13.59& 1.11 & 1.23 \\ \bottomrule
    \end{tabular}
\end{table}

Figure~\ref{fig:exact_comparison_10-20} extends the performance analysis by comparing the computational efficiency and solution accuracy of sALNS with Gurobi on larger instances involving 10, 15, and 20 customers (denoted as C\#) across various scenario counts (denoted as R\#). The bars represent the computational times (left axis), while the overlaid line graph displays the gap between the solutions produced by sALNS and the best result reported by Gurobi within the 12-hour time limit (right axis). This dual representation enables the assessment of sALNS performance in both speed and solution quality.

Notably, for instances with 15 and 20 customers, Gurobi fails to provide any solution for cases involving more than 60 and 30 scenarios, respectively, which is insufficient to accurately capture customer choice behavior. This highlights the impracticality of Gurobi for real-world applications, where a higher number of scenarios is often required. 
In contrast, the sALNS method consistently produces high-quality solutions (within a 10\% gap) in less than 80 seconds for all of the evaluated instances. 

\pgfplotstableread{
1	3.35	8.67
2	6.07	10.26
3	10.43	10.92
4	12.48	10.56
5	19.37	9.81
6	28.18	9.35
7	38.01	8.80
8	42.50	4.62
9	50.51	2.03
10	69.61	3.19
		
12	4.32	7.568984025
13	8.55	7.201120234
14	14.43	5.174624844
15	21.49	2.596973169
16	33.62	2.332392541
17	48.76	3.176782024
18	67.78	1.215338042
19	74.30	0.176433659
		
21	8.30	5.291591673
22	12.25	6.509020546
23	17.24	5.169624577
24	31.47	3.173123496
25	42.95	-9.529712239		
}\datatableFour
\begin{figure}[h!]
\centering
\begin{tikzpicture}[scale=0.75,
  every axis/.style={ 
    ybar,
    x tick label style={rotate=45,font=\scriptsize},
    xticklabels={
    C10-R5,C10-R10,C10-R15,
    C10-R20,C10-R30,C10-R40,C10-R50,C10-R60,C10-R80,C10-R100,
    C15-R5,
    C15-R10,C15-R15,C15-R20,C15-R30,C15-R40,C15-R50,C15-R60,
    C20-R5,C20-R10,C20-R15,C20-R20,C20-R30},
    width = \textwidth,
    height = 7.5cm,
    legend style={draw=none},
    y label style={font=\scriptsize},
    y tick label style={font=\scriptsize} 
  },
]
\begin{axis}[
    ylabel=Time (s),
    xtick=data,
    enlarge y limits=false,
    enlarge x limits=0.05,
    ymin=0,ymax=80,
    legend style={
      font=\scriptsize,
      cells={anchor=west},
      legend columns=5,
      at={(0.27,1.02)},
      anchor=north,
    },
    axis lines*=left,
]
\addplot+[gray] table[x index=0,y index=1] \datatableFour;
\end{axis}
\begin{axis}[
    sharp plot, 
    axis y line*=right,
    axis x line=none,
    ylabel=Gap rate with the best result (\%),
    xtick=data,
    enlarge y limits=false,
    enlarge x limits=0.05,
    ymin=-20,ymax=12,
    legend style={
      font=\scriptsize,
      cells={anchor=west},
      legend columns=5,
      at={(0.62,1.03)},
      anchor=north,
    }
]
\addplot[mark=triangle*, color=Maroon] table[x index=0,y index=2] \datatableFour;
\end{axis}
\end{tikzpicture}
\caption{Comparison of sALNS computation time (bar chart, left axis) and optimality gap (line chart, right axis) relative to Gurobi’s best solution after 12 hours.}
\label{fig:exact_comparison_10-20}
\end{figure}

Figure~\ref{fig:exact_comparison_30-100} highlights the scalability of sALNS, as it continues to deliver efficient solutions for larger instances, including those with up to 100 customers. Despite the increased problem size, solving times remain manageable, with even the most complex instances (C100-R100) being solved within less than 4 hours.
This demonstrates sALNS' ability to handle large-scale, real-world problems efficiently. When considered alongside the results in Table~\ref{tab:small-size_res}, Figure~\ref{fig:exact_comparison_10-20}, and Figure~\ref{fig:exact_comparison_30-100}, the analysis demonstrates that sALNS consistently delivers high-quality solutions across all tested instances, ranging from small to larger problem sizes. For the larger instances (C20–R100 to C100–R100), sALNS solves these within reasonable time frames, demonstrating its scalability. 
In contrast, Gurobi fails to prove optimality in 86\% of the cases, further highlighting the practical value of sALNS for tactical assortment problems, especially in larger-scale settings. 
Moreover, feasibility remains a challenge for Gurobi, as approximately 37\% of the small- to medium-sized instances (5–20 customers) did not yield a feasible solution within the time limit, underscoring the problem’s combinatorial complexity and computational limits.

\pgfplotstableread{
1	161.6910208
2	275.4132907
3	685.427553
4	1233.277314
5	2220.746745
6	3393.542505
7	7141.340776
8	14015.0377
}\datatableThree
\begin{figure}[h!]
\centering
        \begin{tikzpicture}[scale=0.7,
          every axis/.style={ 
            ybar,
            x tick label style={rotate=45,font=\scriptsize},
            xticklabels={C15-R100,C20-R100,C30-R100,C40-R100,C50-R100,C60-R100,C80-R100,C100-R100},%
            ymin=0,ymax=15000,
            height = 7.5cm,
            legend style={draw=none},
            y label style={font=\scriptsize},
            y tick label style={font=\scriptsize} 
          },
        ]
        \begin{axis}[
            ylabel=Time (s),
            xtick=data,
            enlarge y limits=false,
            enlarge x limits=0.1,
            axis lines*=left,
        ]
        \addplot+[gray] table[x index=0,y index=1] \datatableThree;
        \end{axis}
        \end{tikzpicture}       
\caption{Time to solve larger instances (15 to 100 customers) with 100 scenarios using sALNS.}
\label{fig:exact_comparison_30-100}
\end{figure}

To conclude on the performance of sALNS, Figure~\ref{fig:dist&rep} presents a radar chart visualization comparing the contribution of each destroy and repair operator to the search process across all evaluated instances, which include customer sizes ranging from 5 to 100 and varying scenario sizes. Each axis of the radar charts corresponds to one operator, and the three overlaid polygons represent distinct aspects of operator performance: the \textit{Best} polygon (solid blue-violet line) reflects the relative frequency with which each operator contributed to reaching the best-known solution during the search; the \textit{Better} polygon (dotted purple line) indicates how often an operator generated a solution better than the current one; and the \textit{Accepted} polygon (dashed green line) shows the rate at which solutions generated by an operator were accepted into the search process, based on the adaptive acceptance criterion used in sALNS.

Among the destroy operators, all demonstrate balanced performance, but the \textit{Worst destroy} operator stands out, particularly in its ability to target weak solution components and uncover better solutions. For the repair operators, the performance of the random operator appears less effective, highlighting the increased complexity of the problem in such cases and the greater reliance on more structured operators.
the \textit{Two-regret repair} performs best excelling in both routing efficiency and customer coverage. The consistent performance of all operators reflects the robustness and effective integration within the sALNS method, ensuring a balanced search process where no single operator dominates disproportionately. 

\begin{figure}[h!]      
    \captionsetup{justification=centering}
    \begin{subfigure}{0.45\textwidth}
    \begin{tikzpicture}[scale=0.9]
          \path (0:0cm) coordinate (O); 
        
          \def\D{4} 
          \def\R{2.77} 
          \def\A{360/\D} 
          \def\U{35} 
          \def\L{3.42} 
        
          \foreach \X in {1,...,\D}{
            \draw [opacity=0.3] (\X*\A:0) -- (\X*\A:\R);
          }
        
          \foreach \Y in {5,10,15,20,...,\U}{
            \draw [opacity=0.3] (0:\Y*\R/\U) \foreach \X in {1,...,\D}{
                -- (\X*\A:\Y*\R/\U)
            } -- cycle;
            
          }
          \foreach \Y in {0,...,\U}{
            \foreach \X in {1,...,\D}{
              \path [opacity=0.3] (\X*\A:\Y*\R/\U) coordinate (D\X-\Y);
              \fill [opacity=0.1] (D\X-\Y) circle (0.5pt);
            }
          }
          
          \path (1*\A:\L) node (L1) {\microsize Adaptive};
          \path (2*\A:\L) node (L2) {\microsize Random};
          \path (3*\A:\L) node (L3) {\microsize Worst};
          \path (4*\A:\L) node (L4) {\microsize Neighborhood};
        
          \draw [color=BlueViolet,line width=1.0pt]
            (D1-0.8)--
            (D2-0.98)--
            (D3-3.26)--
            (D4-0.96)-- cycle;

          \draw [color=purple,line width=1.0pt, dotted]
            (D1-13.22)--
            (D2-17.89)--
            (D3-24.46)--
            (D4-16.42)-- cycle;
        
          \draw [color=PineGreen,line width=1.0pt,dashed]
            (D1-31.2)--
            (D2-31.03)--
            (D3-31.62)--
            (D4-30.69)-- cycle;
      \draw[very thick, BlueViolet] ($(4.5-0.3,0*.32cm-85)$) -- ++(0.6,0);
      \node[anchor=west] at ($(4.82cm, 0*0.32cm-0.15cm-80)$) {\scriptsize Best};
      \draw[very thick, purple, dotted] ($(4.5-0.3,-1*.32cm-85)$) -- ++(0.6,0);
      \node[anchor=west] at ($(4.82cm, -1*0.32cm-0.15cm-80)$) {\scriptsize Better};
      \draw[very thick, PineGreen, dashed] ($(4.5-0.3,-2*.32cm-85)$) -- ++(0.6,0);
      \node[anchor=west] at ($(4.82cm, -2*0.32cm-0.15cm-80)$) {\scriptsize Accepted};
        \end{tikzpicture}
        \caption{Destroy operators performance.}
        \label{fig:destroy_comparison_obj}
    \end{subfigure}
\hfill
    \begin{subfigure}{0.45\textwidth}
        \begin{tikzpicture}[scale=0.9]
          \path (0:0cm) coordinate (O); 
        
          \def\D{6} 
          \def\R{3.2} 
          \def\A{360/\D} 
          \def\U{25} 
          \def\L{3.85} 
        
          \foreach \X in {1,...,\D}{
            \draw [opacity=0.3] (\X*\A:0) -- (\X*\A:\R);
          }
        
          \foreach \Y in {5,10,15,20,...,\U}{
            \draw [opacity=0.3] (0:\Y*\R/\U) \foreach \X in {1,...,\D}{
                -- (\X*\A:\Y*\R/\U)
            } -- cycle;
            
          }
          \foreach \Y in {0,...,\U}{
            \foreach \X in {1,...,\D}{
              \path [opacity=0.3] (\X*\A:\Y*\R/\U) coordinate (D\X-\Y);
              \fill [opacity=0.1] (D\X-\Y) circle (0.5pt);
            }
          }
          
          \path (1*\A:\L) node (L1) {\microsize Greedy};
          \path (2*\A:\L) node (L2) {\microsize High-utility slot};
          \path (3*\A:\L) node (L3) {\microsize Best};
          \path (4*\A:\L) node (L4) {\microsize Discount-rate-swap};
          \path (5*\A:\L) node (L5) {\microsize Two-regret};
          \path (6*\A:\L) node (L6) {\microsize Random};
        
          \draw [color=BlueViolet,line width=1.0pt]
            (D1-0.96)--
            (D2-1.06)--
            (D3-1.06)--
            (D4-1.02)--
            (D5-1.88)--
            (D6-0.14)-- cycle;
        
          \draw [color=purple,line width=1.0pt, dotted]
            (D1-15.79)--
            (D2-14.39)--
            (D3-16.52)--
            (D4-14.552)--
            (D5-23.62)--
            (D6-2.44)-- cycle;
        
          \draw [color=PineGreen,line width=1.0pt, dashed]
            (D1-22.28)--
            (D2-23.48)--
            (D3-23.04)--
            (D4-21.244)--
            (D5-23.2)--
            (D6-14.22)-- cycle;

        \end{tikzpicture}
        \caption{Repair operators performance.}
        \label{fig:repair_comparison_obj}
    \end{subfigure}
    \caption{Destroy and repair operators performance}
    \label{fig:dist&rep}
\end{figure}

Numerical and computational} insights
\label{managerialinsights} 

This section begins by assessing the benefits of the slot management policy in Section~\ref{sec:policy}, followed by the impact of stochastic customer behaviors in Section~\ref{secUncertainty}. Section~\ref{secML} then explores the importance of capturing customer heterogeneity. In Section~\ref{secDiscount}, we conduct a sensitivity analysis on discount rates. Section~\ref{sec:demographic} tests the solution's robustness across demographic parameters, while Section~\ref{secFairnessCheck} explores the trade-offs between uniform pricing and profitability.

\subsubsection{Impact of tactical assortment and price discount decisions} \label{sec:policy}

To assess the benefits of tactical assortment and price discount decisions, we compare the performance of TTSM-ML against a no-policy baseline using four key metrics: (i) total revenue, (ii) routing cost, (iii) profit (revenue minus cost), and (iv) customer coverage, \textit{i.e.,} the percentage of customers who select a time slot and are thus served. 

Both approaches assume the same ML demand model and are therefore evaluated using sALNS with 100 scenarios and 10 repetitions per instance. The two approaches differ only in their tactical decisions: while the no-policy baseline passively offers all delivery slots at regular prices, TTSM-ML aims to optimize assortment and pricing to steer customer choices and improve system-wide performance.

The evaluation spans customer volumes ranging from 5 to 100. For each size, results are averaged over 15 instances, generated by combining the three spatial configurations and five datasets.

Performance improvements are reported as percentage gaps between TTSM-ML and the no-policy baseline, computed as:
\[
\text{Gap} = \frac{\lambda_{\text{TTSM-ML}} - \lambda_{\text{No-policy}}}{\lambda_{\text{No-policy}}} \times 100,
\]
where $\lambda$ represents revenue, routing cost, profit, or customer coverage.

\begin{table}[h!]
\centering
\caption{Percentage improvement of TTSM-ML over no-policy baseline across instance sizes.}
\label{tab:deterministic_routing}
\resizebox{0.8\textwidth}{!}{%
\begin{tabular}{ccccc}
\toprule
{Customer Size} & {Revenue Gap (\%)} & {Cost Gap (\%)} & {Profit Gap (\%)} & {Coverage Gap (\%)} \\ \midrule
5 & 1.97 & 1.95 & 51.12 & 33.33 \\
10 & 2.26 & 2.22 & 34.00 & 51.43 \\
15 & 4.38 & 4.34 & 11.56 & 18.90 \\
20 & 4.79 & 4.75 & 16.52 & 25.47 \\
30 & 4.44 & 4.39 & 15.55 & 23.16 \\
40 & 0.22 & 0.19 & 8.12 & 15.58 \\
50 & 1.04 & 1.01 & 7.51 & 14.26 \\
60 & 0.34 & 0.31 & 7.78 & 15.72 \\
80 & 2.30 & 2.26 & 9.02 & 14.09 \\
100 & 2.06 & 2.02 & 10.63 & 17.58 \\
\bottomrule
\end{tabular}%
}
\end{table}

The results, summarized in Table~\ref{tab:deterministic_routing}, show that tactical decisions affect both economic and operational outcomes. TTSM-ML consistently enhances profitability across all customer volumes. These improvements stem from a combination of increased revenue, driven by higher slot acceptance (reflected in improved customer coverage) and better price alignment with customer preferences, and reduced routing costs, achieved by improving slot offerings to encourage geographical clustering and more efficient vehicle usage.

\subsubsection{Impact of considering stochastic customer behavior} \label{secUncertainty}

To evaluate the benefits of incorporating customer behavior uncertainty in slot management, we use two key metrics: the Value of Stochastic Solution (VSS) and the Expected Value of Perfect Information (EVPI). 

The VSS measures the profit improvement achieved by explicitly modeling customer behavior uncertainty compared to a deterministic approach. Specifically, we compare the profit from our stochastic model (TTSM-ML) with that from a deterministic counterpart, referred to as the Deterministic TTSM (D-TTSM) model. The D-TTSM eliminates all randomness in customer preferences by setting the random utility components ($\beta^{\mathrm{price}}, \xi$) to their expected values, resulting in deterministic customer behavior. After optimizing the D-TTSM for instances with 60, 80, and 100 customers, the deterministic solution is re-evaluated using the stochastic profit function over 100 scenarios to obtain ($\text{Obj}_{\text{D-TTSM}}$), and then compared with the stochastic TTSM-ML objective value ($\text{Obj}_{\text{TTSM}}$). The VSS is then calculated as:
\begin{equation*}
    \text{VSS} = \frac{\text{Obj}_{\text{TTSM}} - \text{Obj}_{\text{D-TTSM}}}{\text{Obj}_{\text{TTSM}}} \times 100
\end{equation*}
This percentage reflects the benefit of incorporating uncertainty in the decision-making process.

The EVPI quantifies the maximum possible improvement in profit achievable if uncertainty about customer behavior were completely eliminated. For each scenario \(r \in R\), we compute the optimal decisions assuming perfect knowledge of customer choices, then average the resulting profits to get \(\text{Obj}_{\text{PI}}\). The EVPI relative to the stochastic model profit is:
\begin{equation*}
    \text{EVPI} = \frac{\text{Obj}_{\text{PI}} - \text{Obj}_{\text{TTSM}}}{\text{Obj}_{\text{TTSM}}} \times 100
\end{equation*}

Table~\ref{tab:uncertainty} presents VSS and EVPI values across different instance sizes, map configurations and datasets (thus behavioral settings), along with the ratio \(\frac{\text{VSS}}{\text{VSS} + \text{EVPI}} \times 100\%\), indicating the proportion of the total possible improvement captured by the stochastic solution.

\begin{table}[!htbp]
\centering
\scriptsize
\caption{Comparison of the deviation between stochastic and deterministic solutions (VSS), and EVPI across 60-, 80-, and 100-customer instances with 100 scenarios. The last column shows the percentage of the total possible improvement (VSS + EVPI) captured by the stochastic solution.}
\label{tab:uncertainty}
\resizebox{.6\textwidth}{!}{%
\begin{tabular}{cccccc}
\toprule
\multicolumn{1}{l}{Customer} &
  \multicolumn{1}{l}{Map} &
  Dataset &
  \multicolumn{1}{l}{VSS} &
  \multicolumn{1}{l}{EVPI} &
  \multicolumn{1}{l}{$\frac{\text{VSS}}{\text{VSS}+\text{EVPI}}$(\%)} \\ \midrule
\multirow{16}{*}{60}  & \multirow{5}{*}{Random}    & 1 & 17.33  & 70.67  & 19.69 \\
                      &                            & 2 & 59.61  & 40.16  & 59.75 \\
                      &                            & 3 & 24.72  & 49.09  & 33.49 \\
                      &                            & 4 & 87.43  & 88.98  & 49.56 \\
                      &                            & 5 & 9.10   & 76.91  & 10.58 \\ \cline{2-6}
                      & \multirow{5}{*}{Clustered} & 1 & 15.79  & 64.94  & 19.56 \\
                      &                            & 2 & 47.37  & 33.13  & 58.84 \\
                      &                            & 3 & 22.39  & 91.19  & 19.71 \\
                      &                            & 4 & 65.01  & 67.71  & 48.99 \\
                      &                            & 5 & 8.49   & 76.63  & 9.98  \\ \cline{2-6}
                      & \multirow{5}{*}{Mixed}     & 1 & 13.51  & 98.62  & 12.05 \\
                      &                            & 2 & 33.32  & 37.88  & 46.80 \\
                      &                            & 3 & 18.13  & 43.97  & 29.19 \\
                      &                            & 4 & 41.73  & 38.23  & 52.19 \\
                      &                            & 5 & 7.64   & 78.56  & 8.86  \\ \midrule
                      Average &                            &    & \textbf{31.44}  & \textbf{63.78}  & \textbf{31.95}  \\ 
                      \midrule
\multirow{16}{*}{80}  & \multirow{5}{*}{Random}    & 1 & 20.99  & 78.97  & 21.00 \\
                      &                            & 2 & 10.63  & 66.34  & 13.81 \\
                      &                            & 3 & 43.11  & 62.70  & 40.74 \\
                      &                            & 4 & 53.15  & 28.34  & 65.22 \\
                      &                            & 5 & 11.06  & 74.03  & 13.00 \\ \cline{2-6}
                      & \multirow{5}{*}{Clustered} & 1 & 19.14  & 73.32  & 20.70 \\
                      &                            & 2 & 73.51  & 38.43  & 65.67 \\
                      &                            & 3 & 36.13  & 54.54  & 39.85 \\
                      &                            & 4 & 16.66  & 10.78  & 60.73 \\
                      &                            & 5 & 10.28  & 69.40  & 12.91 \\ \cline{2-6}
                      & \multirow{5}{*}{Mixed}     & 1 & 16.00  & 72.06  & 18.17 \\
                      &                            & 2 & 44.93  & 27.03  & 62.44 \\
                      &                            & 3 & 27.35  & 44.32  & 38.16 \\
                      &                            & 4 & 71.54  & 52.37  & 57.74 \\
                      &                            & 5 & 9.17   & 75.90  & 10.78 \\ \midrule
                      Average &                            &    & \textbf{30.91}  & \textbf{55.24}  & \textbf{36.06}  \\ 
                      \midrule
\multirow{16}{*}{100} & \multirow{5}{*}{Random}    & 1 & 26.99  & 76.85  & 25.99 \\
                      &                            & 2 & 24.14  & 47.94  & 33.49 \\
                      &                            & 3 & 58.37  & 71.33  & 45.00 \\
                      &                            & 4 & 16.23  & 106.98 & 13.17 \\
                      &                            & 5 & 12.76  & 100.80 & 11.24 \\ \cline{2-6}
                      & \multirow{5}{*}{Clustered} & 1 & 24.01  & 71.60  & 25.11 \\
                      &                            & 2 & 24.87  & 81.72  & 23.33 \\
                      &                            & 3 & 47.26  & 66.31  & 41.61 \\
                      &                            & 4 & 36.28  & 20.19  & 64.25 \\
                      &                            & 5 & 11.97  & 106.65 & 10.09 \\ \cline{2-6}
                      & \multirow{5}{*}{Mixed}     & 1 & 20.19  & 68.86  & 22.67 \\
                      &                            & 2 & 96.62  & 44.58  & 68.43 \\
                      &                            & 3 & 36.15  & 60.48  & 37.41 \\
                      &                            & 4 & 33.78 & 14.08 & 70.57 \\
                      &                            & 5 & 10.82  & 99.64  & 9.79  \\ \midrule
                      Average &                            &    & \textbf{32.03}  & \textbf{69.20}  & \textbf{33.48}  \\ 
                      \midrule
                      \textbf{Total Average} &                            &    & \textbf{31.46}  & \textbf{62.74}  & \textbf{33.83}  \\ 
                      \bottomrule
\end{tabular}%
}
\end{table}

The results show that explicitly modeling customer behavior uncertainty yields an average profit improvement of approximately 31.5\% over the deterministic baseline (VSS). The EVPI averages around 62.7\%, representing the potential gain with perfect information. Our stochastic model captures over 30\% of this gap, proving its effectiveness in improving decisions under uncertainty. While perfect information is unrealistic, these results clearly demonstrate the practical benefits of accounting for customer behavior uncertainty in tactical slot management, leading to real profit and service improvements.

\subsubsection{Impact of considering customers' heterogeneity} \label{secML}

An important contribution of TTSM-ML is its ability to capture individual-level heterogeneity, providing a more nuanced understanding of customer preferences. As observed in Section~\ref{alg:discretechoice}, the ML model outperforms the simpler MNL model on our datasets. We now evaluate whether incorporating this heterogeneity into the optimization process leads to better tactical decisions and improved profitability, or if the simpler MNL approach achieves similar outcomes despite its lower complexity. To this end, we compare the performance of the TTSM-ML model, which incorporates individual-level preferences through random coefficients ($\bar{\beta}^{\mathrm{price}}_{r}$), with the TTSM-MNL model, which assumes homogeneous preferences represented by fixed coefficients ($\beta^{\mathrm{price}}$).

We solved both models on an extended set of instances featuring the five datasets (\textit{i.e.}, behavioral settings), with five instances from each of the three spatial configurations (random, clustered, and mixed), and across three customer sizes: 60, 80, and 100.

Each instance was solved over 100 stochastic scenarios ($R=100$), and the resulting solutions were re-evaluated with a new set of 100 scenarios (with $R=1$ per evaluation) to obtain robust profit estimates. We then conducted paired t-tests to compare the profits of TTSM-ML and TTSM-MNL for each instance.
This process yields a comprehensive set of paired outcomes comparing the performance of TTSM-ML and TTSM-MNL. To assess their performance, we employed a paired t-test. This test helps determine whether the mean differences in profits between the two models are statistically significant. The t-test produces a t-value, which quantifies the size of the difference relative to the variability of the results, and a p-value, which indicates whether the observed difference is statistically significant. 

Figure~\ref{fig:t_test_summary} provides a 3D summary of the paired t-test results across all experimental settings. Each point represents a specific configuration, defined by a unique combination of dataset (\textit{i.e.} behavioral setting), map, and customer size. For each such configuration, five independent instances were tested as explained. The color of the point indicates the percentage of those five instances in which TTSM-ML significantly outperformed TTSM-MNL at the 95\% confidence level ($p < 0.05$). The exact t-values and corresponding p-values for each tested instance are provided in Appendix~\ref{appSec:t-test}.

As shown in Figure~\ref{fig:t_test_summary}, TTSM-ML significantly outperforms TTSM-MNL in approximately 90\% of evaluated cases with statistical significance ($p < 0.05$). In contrast, in clustered spatial settings (where customer locations are concentrated) the advantage of TTSM-ML is less pronounced, with around 80\% of cases showing statistically significant improvements. This difference may arise because clustered locations tend to lower routing costs, thereby reducing the incremental profit potential gained from improved slot assignments.

These results demonstrate the practical value of explicitly modeling individual-level preference heterogeneity using the ML model. When customer behavior exhibits meaningful heterogeneity, our proposed framework generates more profitable and tailored slot decisions compared to those based on the less fitting MNL model.


\begin{figure}[h!]
\centering
\footnotesize
\begin{tikzpicture}
\begin{axis}[
    view={70}{30},
    xtick={1,2,3,4,5},
    xlabel={Dataset (behavioral setting)},
    ylabel={Map setting},
    ytick={1,2,3},
    yticklabels={Random,Clustered,Mixed},
    zlabel={Number of customers},
    ztick={60, 80, 100},
    colorbar,
    colormap/viridis,
    point meta min=40,
    point meta max=100,
    scatter/use mapped color={
        draw=black, fill=mapped color
    },
    enlargelimits,
    xlabel style={rotate=-50},
    ylabel style={},
]
\addplot3[
    scatter, only marks,
    point meta=explicit, 
    mark=*,
    visualization depends on={0.02*\thisrow{color} \as \perpointmarksize},
    scatter/@pre marker code/.append style={
        /tikz/mark size=\perpointmarksize,
    }
]
table [x=x, y=y, z=z, meta=color] {scatterdata.csv};
\end{axis}
\end{tikzpicture}
\caption{Results of the paired t-test comparing the performance of the TTSM-ML and TTSM-MNL models. The figure indicates percentage of cases where TTSM-ML significantly outperforms TTSM-MNL at the $\alpha = 0.05$ significance level (95\% confidence level).}
\label{fig:t_test_summary}
\end{figure}

\subsubsection{Sensitivity analysis on discount rates} \label{secDiscount}

This section examines how discount rates affect profitability across different customer mapping strategies.
We analyze the trade-off between revenue generation through customer attraction and increased operational costs. As illustrated in Figure~\ref{fig:discount}, our findings reveal that profitability exhibits a non-linear response to discount rates. Here, profit improvement is measured relative to the base case of 0\% discount, \textit{i.e.}, the difference between the profit at a given discount level and the profit obtained using sALNS with no discounts. Each data point represents the average of 20 repetitions, with standard deviation error bars added to capture variability. The figure reports results for the first dataset as a representative example.

\begin{figure}[h!]
\centering
\footnotesize
\begin{tikzpicture}
  \begin{axis}[
      ytick align=outside, ytick style={draw=none},
      xtick align=outside, xtick style={draw=none},
      height=\axisdefaultheight,
      grid=both,
      xlabel={Discount rate (\%)},
      ylabel={Profit improvement (\%)},
    xlabel={Discount rate (\%)},
    ylabel={Profit improvement (\%)},
    xtick={0,10,20,30,40,50,60,70,80},
    ytick={0,5,...,40},
    legend style={
        at={(0.5,1.1)}, 
        anchor=south,  
        draw=none,     
        fill=none      
    },
    legend columns=-1, 
  ]
    \addplot+[
      BlueViolet,
      mark options={BlueViolet, scale=0.0},
      smooth,
      error bars/.cd,
        y dir=both,
        y explicit
    ]
    coordinates {
      (0,0) +- (0,0)
      (5,18.87) +- (0,1.5)
      (10,5.68) +- (0,1.2)
      (15,13.67) +- (0,2.1)
      (20,0.65) +- (0,0.8)
      (30,22.27) +- (0,1.6)
      (40,8.76) +- (0,1.4)
      (50,12.71) +- (0,1.9)
      (60,19.43) +- (0,2.0)
      (70,14.41) +- (0,1.7)
      (80,23.62) +- (0,2.5)
    };
    \addlegendentry{Clustered map}

    \addplot+[
      teal,
      mark options={teal, scale=0.0},
      smooth,
      style = dashdotted,
      error bars/.cd,
        y dir=both,
        y explicit
    ]
    coordinates {
      (0,0)  +- (0,0)
      (5,2.533013938)  +- (0,0.6)
      (10,13.9499002)  +- (0,0.7)
      (15,20.25983358)  +- (0,0.2)
      (20,2.657339693)  +- (0,0.8)
      (30,5.065748863)  +- (0,1.2)
      (40,6.71766497)  +- (0,1.7)
      (50,2.569740036)  +- (0,2.3)
      (60,13.47166456)  +- (0,0.6)
      (70,1.818741357)  +- (0,1.3)
      (80,0.670838975) +- (0,1.9)
    };
    \addlegendentry{Random map}

    \addplot+[
      purple,
      mark options={purple, scale=0.0},
      smooth,
      style = dashed,
      error bars/.cd,
        y dir=both,
        y explicit
    ]
    coordinates {
    (0,0)  +- (0,0)
    (5,34.40019275)  +- (0,2.1)
    (10,39.36391712)  +- (0,1.1)
    (15,29.88518504)  +- (0,0.9)
    (20,10.46140583)  +- (0,0.7)
    (30,8.613546609)  +- (0,0.4)
    (40,1.808119834)  +- (0,0.8)
    (50,8.687339553)  +- (0,1.8)
    (60,3.961902595)  +- (0,0.6)
    (70,2.075104741)  +- (0,1.1)
    (80,23.0771224) +- (0,1.0)
    };
    \addlegendentry{Mixed map}

  \end{axis}
\end{tikzpicture}
\caption{Impact of discount rates on improving profits for 10 customers and different geographic settings solved with sALNS. The results show that the most effective discount rates vary by map pattern, with each exhibiting distinct profit peaks.}
\label{fig:discount}
\end{figure}

Increasing the discount level boosts customer retention, which drives higher revenues that outweigh additional routing costs. Profit improvement from higher discount rates suggests successful customer retention. However, if further increases in discount rates cause profit gains to plateau or return to zero, it indicates that retained customers no longer generate enough revenue to offset the rising discount and routing costs. Conversely, if profit improvement declines but stays positive, it suggests that while the discount reduces revenue, retained customers still contribute enough to maintain profitability above the baseline scenario with no discount.

The peaks in the profit improvement curves represent critical points where customers re-enter the system due to the discount incentives, signaling optimal balance points where retained customer revenue outweighs operational and discount costs. Beyond these points, diminishing returns become evident, as excessive discounts lower revenue per order while routing costs increase disproportionately. The impact of discount rates is also influenced by customer spatial distribution, as seen in the variations across the three mapping strategies. The results indicate that the most effective discount rates depend on both customer price sensitivity and operational costs, with effectiveness varying across different spatial distributions. This suggests the need for tailored discount strategies based on specific settings.

\subsubsection{Assessing the solution across demographic parameters} \label{sec:demographic}

In this part, we analyze how the slot assortment decisions can be adjusted based on consumer sensitivities to time ($\beta^{time}$) and price ($\beta^{\mathrm{price}}$) to either retain or exclude customers, with the goal of improving overall profitability. Figure~\ref{fig:coverage} shows how the coverage rate, defined as the proportion of customers who remain engaged with the system (\textit{i.e.}, who don’t opt out), varies with different sensitivity parameters. This is illustrated in two cases: (a) offering all slots without any discount and (b) applying assortment decisions obtained through the sALNS. The analysis is based on the first dataset, using a representative instance with 10 customers on the mixed map.

Figure~\ref{fig:coverage_without} depicts the coverage rate when no optimization is applied, \textit{i.e.}, when all slots are offered without any discount. The $\beta^{time}$ axis represents sensitivity to slot timing, where values closer to zero indicate less sensitivity, meaning assortment decisions have less impact on retention. The $\beta^{\mathrm{price}}$ axis reflects sensitivity to pricing changes, with lower values indicating reduced responsiveness to discounts.

As expected, as customers become more price-sensitive, they opt out more frequently, causing the coverage rate to decline. Notably, a remarkable portion of the coverage rate stays below 20\% (gray surface) across a wide range of sensitivities, indicating poor retention performance without optimization. This outcome suggests that the utility of opting out outweighs the utility of the offered slots. The initial decline in coverage occurs because customers who are more sensitive to time and price are likely to opt out when their preferences aren’t met. In Figure~\ref{fig:coverage_without}, this drop happens earlier compared to the unoptimized case, reflecting a trade-off made by the sALNS method. By prioritizing profitability over customer retention, the algorithm avoids targeting high-cost customers and focuses on maximizing the profitability of slot arrangements. In contrast, Figure~\ref{fig:coverage_with-without_fairness} shows the coverage rate when assortment decisions are determined using our sALNS method. The trend remains similar: as customers become more sensitive to time or price, retaining them becomes more difficult. However, the coverage rate improves, consistently staying above the 20\% threshold (gray surface) across all tested settings. This improvement shows that adjusting assortment and discount decisions based on sensitivities enables better customer capture and increased profitability.

To evaluate sALNS performance, we compare its profitability improvements against RFTS across varying levels of customer sensitivity. Figure~\ref{fig:improvement_obj_ALNS_RFTS} presents these improvements. It demonstrates that sALNS remarkably outperforms RFTS, particularly in scenarios where customers show lower or moderate sensitivity to both price and time. This is because highly sensitive customers tend to opt out regardless. sALNS excels at identifying opportunities where customer preferences can be effectively leveraged. For customers highly responsive to both price ($\beta^{\mathrm{price}}$) and time ($\beta^{time}$), sALNS achieves up to a 45\% profitability improvement over RFTS, due to its ability to strategically adjust slot offerings and pricing. For customers with moderate sensitivities, sALNS consistently outperforms RFTS by 15-30\%, demonstrating its adaptability across different customer types. This robust performance across varying sensitivities underscores the algorithm’s versatility in handling diverse customer preferences. 

\pgfplotstableread{
-2	-0.2	20.8
-1.8	-0.2	21
-1.4	-0.2	22.8
-1.2	-0.2	23.6
-1	-0.2	24
-0.9	-0.2	24.2
-0.8	-0.2	24.4
-0.7	-0.2	24.6
-0.6	-0.2	24.2
-0.5	-0.2	24.8
-0.4	-0.2	25.6
-0.3	-0.2	26.8
-0.2	-0.2	26.8
-0.1	-0.2	29.8
0	-0.2	29.4
0.1	-0.2	30.4
0.2	-0.2	29
0.3	-0.2	30.2
0.4	-0.2	31
0.5	-0.2	33.8
0.6	-0.2	33.8
0.7	-0.2	34.2
0.8	-0.2	35.8
0.9	-0.2	34.6
1	-0.2	35.8
1.2	-0.2	35.8
1.4	-0.2	38.4
1.6	-0.2	38.4
1.8	-0.2	40.4
2	-0.2	42.6
-2	-0.18	23.4
-1.8	-0.18	23.8
-1.4	-0.18	24.8
-1.2	-0.18	25
-1	-0.18	26.8
-0.9	-0.18	27.6
-0.8	-0.18	29.4
-0.7	-0.18	28.4
-0.6	-0.18	28.6
-0.5	-0.18	28.6
-0.4	-0.18	30.6
-0.3	-0.18	29.8
-0.2	-0.18	31.4
-0.1	-0.18	32
0	-0.18	32
0.1	-0.18	34
0.2	-0.18	34.8
0.3	-0.18	36
0.4	-0.18	36.2
0.5	-0.18	35.8
0.6	-0.18	38.4
0.7	-0.18	38.8
0.8	-0.18	37.8
0.9	-0.18	39.2
1	-0.18	39.2
1.2	-0.18	41.8
1.4	-0.18	42.2
1.6	-0.18	41.2
1.8	-0.18	43.8
2	-0.18	44.6
-2	-0.14	28.8
-1.8	-0.14	29.8
-1.4	-0.14	33.6
-1.2	-0.14	35
-1	-0.14	35.4
-0.9	-0.14	36
-0.8	-0.14	38.6
-0.7	-0.14	37.6
-0.6	-0.14	39
-0.5	-0.14	38.2
-0.4	-0.14	39.8
-0.3	-0.14	40.8
-0.2	-0.14	42.6
-0.1	-0.14	40.8
0	-0.14	41.4
0.1	-0.14	42.4
0.2	-0.14	42
0.3	-0.14	44
0.4	-0.14	44.4
0.5	-0.14	45.2
0.6	-0.14	43.8
0.7	-0.14	46.4
0.8	-0.14	46.8
0.9	-0.14	45.2
1	-0.14	45.6
1.2	-0.14	47
1.4	-0.14	48
1.6	-0.14	48.2
1.8	-0.14	48.8
2	-0.14	50.2
-2	-0.12	33.2
-1.8	-0.12	35
-1.4	-0.12	39.6
-1.2	-0.12	40.8
-1	-0.12	41.6
-0.9	-0.12	42.4
-0.8	-0.12	42
-0.7	-0.12	42.4
-0.6	-0.12	42.8
-0.5	-0.12	42.8
-0.4	-0.12	43.4
-0.3	-0.12	44
-0.2	-0.12	43.8
-0.1	-0.12	45.2
0	-0.12	44.8
0.1	-0.12	45.2
0.2	-0.12	45.6
0.3	-0.12	46.8
0.4	-0.12	47.2
0.5	-0.12	47.2
0.6	-0.12	47.2
0.7	-0.12	48
0.8	-0.12	48.4
0.9	-0.12	48.2
1	-0.12	49
1.2	-0.12	49.6
1.4	-0.12	50.2
1.6	-0.12	52.6
1.8	-0.12	52.8
2	-0.12	53.2
-2	-0.1	40.4
-1.8	-0.1	41.2
-1.4	-0.1	42.2
-1.2	-0.1	43.2
-1	-0.1	44
-0.9	-0.1	44.8
-0.8	-0.1	45.2
-0.7	-0.1	45.2
-0.6	-0.1	45.6
-0.5	-0.1	46.4
-0.4	-0.1	47.6
-0.3	-0.1	47.6
-0.2	-0.1	47.2
-0.1	-0.1	47.8
0	-0.1	47.8
0.1	-0.1	48.4
0.2	-0.1	49.2
0.3	-0.1	49.8
0.4	-0.1	49.8
0.5	-0.1	49.8
0.6	-0.1	50.4
0.7	-0.1	51
0.8	-0.1	51.8
0.9	-0.1	52.6
1	-0.1	52.8
1.2	-0.1	53.4
1.4	-0.1	54.6
1.6	-0.1	55.2
1.8	-0.1	56.4
2	-0.1	57
-2	-0.09	41.4
-1.8	-0.09	42.4
-1.4	-0.09	44
-1.2	-0.09	46.4
-1	-0.09	45.6
-0.9	-0.09	46.4
-0.8	-0.09	46.6
-0.7	-0.09	47
-0.6	-0.09	47.2
-0.5	-0.09	47.6
-0.4	-0.09	48.4
-0.3	-0.09	48.8
-0.2	-0.09	49.6
-0.1	-0.09	49.8
0	-0.09	49.8
0.1	-0.09	50.6
0.2	-0.09	50.2
0.3	-0.09	51.4
0.4	-0.09	52.4
0.5	-0.09	52.4
0.6	-0.09	53
0.7	-0.09	53.2
0.8	-0.09	53.2
0.9	-0.09	53
1	-0.09	54.2
1.2	-0.09	55.2
1.4	-0.09	56
1.6	-0.09	56.6
1.8	-0.09	56.8
2	-0.09	58.2
-2	-0.08	42.8
-1.8	-0.08	44
-1.4	-0.08	45.6
-1.2	-0.08	46.8
-1	-0.08	47.6
-0.9	-0.08	47.8
-0.8	-0.08	48.4
-0.7	-0.08	48
-0.6	-0.08	49
-0.5	-0.08	49.2
-0.4	-0.08	49.2
-0.3	-0.08	50.4
-0.2	-0.08	50
-0.1	-0.08	51
0	-0.08	52.2
0.1	-0.08	52.6
0.2	-0.08	52.6
0.3	-0.08	53.2
0.4	-0.08	53.4
0.5	-0.08	53
0.6	-0.08	53.8
0.7	-0.08	54.8
0.8	-0.08	55.2
0.9	-0.08	55.8
1	-0.08	56
1.2	-0.08	56.6
1.4	-0.08	56.6
1.6	-0.08	58.2
1.8	-0.08	59.2
2	-0.08	60.6
-2	-0.07	46.4
-1.8	-0.07	45.6
-1.4	-0.07	47.6
-1.2	-0.07	48
-1	-0.07	49.2
-0.9	-0.07	49.4
-0.8	-0.07	50.2
-0.7	-0.07	50.2
-0.6	-0.07	50
-0.5	-0.07	51
-0.4	-0.07	52
-0.3	-0.07	52.8
-0.2	-0.07	52.6
-0.1	-0.07	53
0	-0.07	53.4
0.1	-0.07	53
0.2	-0.07	53.8
0.3	-0.07	54.6
0.4	-0.07	55.2
0.5	-0.07	55.8
0.6	-0.07	56
0.7	-0.07	57
0.8	-0.07	57
0.9	-0.07	57
1	-0.07	57.6
1.2	-0.07	58
1.4	-0.07	58.8
1.6	-0.07	60.6
1.8	-0.07	62.2
2	-0.07	63
-2	-0.06	46.8
-1.8	-0.06	47.4
-1.4	-0.06	49.6
-1.2	-0.06	50
-1	-0.06	50
-0.9	-0.06	51.2
-0.8	-0.06	51.8
-0.7	-0.06	52.8
-0.6	-0.06	52.8
-0.5	-0.06	52.6
-0.4	-0.06	53.4
-0.3	-0.06	53
-0.2	-0.06	53.8
-0.1	-0.06	54.4
0	-0.06	55
0.1	-0.06	55.8
0.2	-0.06	56
0.3	-0.06	56.4
0.4	-0.06	57.4
0.5	-0.06	57.6
0.6	-0.06	57.6
0.7	-0.06	57.8
0.8	-0.06	57.2
0.9	-0.06	58.2
1	-0.06	58.8
1.2	-0.06	59.6
1.4	-0.06	62.2
1.6	-0.06	62.8
1.8	-0.06	63
2	-0.06	63.2
-2	-0.05	48
-1.8	-0.05	49
-1.4	-0.05	50.4
-1.2	-0.05	51.8
-1	-0.05	53
-0.9	-0.05	52.6
-0.8	-0.05	53.4
-0.7	-0.05	53
-0.6	-0.05	53.8
-0.5	-0.05	54.4
-0.4	-0.05	55
-0.3	-0.05	55.8
-0.2	-0.05	56
-0.1	-0.05	56.2
0	-0.05	57
0.1	-0.05	57.2
0.2	-0.05	57.6
0.3	-0.05	57.2
0.4	-0.05	57.2
0.5	-0.05	57.4
0.6	-0.05	58.6
0.7	-0.05	58.6
0.8	-0.05	60.2
0.9	-0.05	61.2
1	-0.05	61.6
1.2	-0.05	62.2
1.4	-0.05	63.2
1.6	-0.05	63.2
1.8	-0.05	64
2	-0.05	65.6
-2	-0.04	49.6
-1.8	-0.04	50.4
-1.4	-0.04	52.4
-1.2	-0.04	52.8
-1	-0.04	53.8
-0.9	-0.04	54.2
-0.8	-0.04	55
-0.7	-0.04	55.6
-0.6	-0.04	55.8
-0.5	-0.04	56.8
-0.4	-0.04	57.4
-0.3	-0.04	57.6
-0.2	-0.04	57.6
-0.1	-0.04	58
0	-0.04	57.6
0.1	-0.04	57.8
0.2	-0.04	58.6
0.3	-0.04	59
0.4	-0.04	59.6
0.5	-0.04	60.8
0.6	-0.04	61.6
0.7	-0.04	61.8
0.8	-0.04	62.2
0.9	-0.04	62.6
1	-0.04	63.2
1.2	-0.04	63.2
1.4	-0.04	64
1.6	-0.04	65.6
1.8	-0.04	66.2
2	-0.04	68.2
-2	-0.03	51.8
-1.8	-0.03	53
-1.4	-0.03	53.8
-1.2	-0.03	55
-1	-0.03	56
-0.9	-0.03	56.6
-0.8	-0.03	57
-0.7	-0.03	57.6
-0.6	-0.03	57.6
-0.5	-0.03	57.8
-0.4	-0.03	57.6
-0.3	-0.03	57.6
-0.2	-0.03	58.2
-0.1	-0.03	59
0	-0.03	59.4
0.1	-0.03	60.8
0.2	-0.03	61.2
0.3	-0.03	61.8
0.4	-0.03	62.6
0.5	-0.03	62.6
0.6	-0.03	63
0.7	-0.03	63.6
0.8	-0.03	63.2
0.9	-0.03	63.6
1	-0.03	63.8
1.2	-0.03	65.4
1.4	-0.03	66.2
1.6	-0.03	68
1.8	-0.03	69
2	-0.03	69.8
-2	-0.02	52.8
-1.8	-0.02	53.8
-1.4	-0.02	56.2
-1.2	-0.02	57.2
-1	-0.02	57.6
-0.9	-0.02	58
-0.8	-0.02	57.8
-0.7	-0.02	58
-0.6	-0.02	58.4
-0.5	-0.02	59
-0.4	-0.02	59.4
-0.3	-0.02	60.2
-0.2	-0.02	61.4
-0.1	-0.02	61.8
0	-0.02	62.2
0.1	-0.02	62.6
0.2	-0.02	63
0.3	-0.02	63.6
0.4	-0.02	63.2
0.5	-0.02	63.6
0.6	-0.02	63.8
0.7	-0.02	64.4
0.8	-0.02	65.4
0.9	-0.02	65.8
1	-0.02	66.2
1.2	-0.02	68
1.4	-0.02	69
1.6	-0.02	69.8
1.8	-0.02	71.4
2	-0.02	71.8
-2	-0.01	55
-1.8	-0.01	56.2
-1.4	-0.01	57.6
-1.2	-0.01	57.6
-1	-0.01	58.2
-0.9	-0.01	59
-0.8	-0.01	59.6
-0.7	-0.01	60
-0.6	-0.01	60.8
-0.5	-0.01	61.4
-0.4	-0.01	62.2
-0.3	-0.01	62.6
-0.2	-0.01	63
-0.1	-0.01	63.4
0	-0.01	63.2
0.1	-0.01	63.6
0.2	-0.01	63.8
0.3	-0.01	64.2
0.4	-0.01	65.4
0.5	-0.01	65.6
0.6	-0.01	66.2
0.7	-0.01	67.6
0.8	-0.01	68
0.9	-0.01	68.4
1	-0.01	69.2
1.2	-0.01	69.6
1.4	-0.01	71.4
1.6	-0.01	71.8
1.8	-0.01	72.2
2	-0.01	72.8
-2	0	56.6
-1.8	0	56.8
-1.4	0	58.4
-1.2	0	59.4
-1	0	60.6
-0.9	0	61.6
-0.8	0	62.2
-0.7	0	62.6
-0.6	0	63
-0.5	0	63
-0.4	0	63.2
-0.3	0	63.4
-0.2	0	63.8
-0.1	0	64.4
0	0	65.2
0.1	0	65.6
0.2	0	66.2
0.3	0	67.6
0.4	0	68
0.5	0	68.6
0.6	0	69
0.7	0	69.2
0.8	0	69.6
0.9	0	70.8
1	0	71.4
1.2	0	71.8
1.4	0	72.2
1.6	0	72.8
1.8	0	74
2	0	74.4
-2	0.01	57.6
-1.8	0.01	58.2
-1.4	0.01	60.2
-1.2	0.01	61.8
-1	0.01	63
-0.9	0.01	63.2
-0.8	0.01	63
-0.7	0.01	63.6
-0.6	0.01	63.8
-0.5	0.01	64.2
-0.4	0.01	65.2
-0.3	0.01	65.4
-0.2	0.01	66
-0.1	0.01	67.6
0	0.01	68
0.1	0.01	68.4
0.2	0.01	68.6
0.3	0.01	69.4
0.4	0.01	69.6
0.5	0.01	70.6
0.6	0.01	71.4
0.7	0.01	71.6
0.8	0.01	71.8
0.9	0.01	72.2
1	0.01	72.2
1.2	0.01	72.8
1.4	0.01	73.4
1.6	0.01	74.6
1.8	0.01	75
2	0.01	75.8
-2	0.02	59.4
-1.8	0.02	60.8
-1.4	0.02	63
-1.2	0.02	63
-1	0.02	63.8
-0.9	0.02	64.2
-0.8	0.02	65.2
-0.7	0.02	65.4
-0.6	0.02	66
-0.5	0.02	67
-0.4	0.02	68
-0.3	0.02	68.4
-0.2	0.02	69
-0.1	0.02	69.2
0	0.02	69.6
0.1	0.02	70.6
0.2	0.02	71.4
0.3	0.02	71.6
0.4	0.02	71.8
0.5	0.02	72.2
0.6	0.02	72.2
0.7	0.02	72.8
0.8	0.02	72.8
0.9	0.02	73.2
1	0.02	73.4
1.2	0.02	74.4
1.4	0.02	75.4
1.6	0.02	75.8
1.8	0.02	76.2
2	0.02	76.8
-2	0.03	61.4
-1.8	0.03	63
-1.4	0.03	63.8
-1.2	0.03	65.2
-1	0.03	66
-0.9	0.03	67
-0.8	0.03	68
-0.7	0.03	68.4
-0.6	0.03	69
-0.5	0.03	69.2
-0.4	0.03	69.6
-0.3	0.03	70.6
-0.2	0.03	71.2
-0.1	0.03	71.6
0	0.03	71.8
0.1	0.03	72.2
0.2	0.03	72.2
0.3	0.03	72.8
0.4	0.03	72.8
0.5	0.03	73.2
0.6	0.03	73.4
0.7	0.03	73.6
0.8	0.03	74.4
0.9	0.03	75.2
1	0.03	75
1.2	0.03	75.6
1.4	0.03	76.2
1.6	0.03	76.8
1.8	0.03	78.6
2	0.03	79.6
-2	0.04	63
-1.8	0.04	63.6
-1.4	0.04	65.6
-1.2	0.04	67.4
-1	0.04	68.4
-0.9	0.04	69.2
-0.8	0.04	69.6
-0.7	0.04	70.6
-0.6	0.04	71.2
-0.5	0.04	71.4
-0.4	0.04	71.8
-0.3	0.04	72.2
-0.2	0.04	72
-0.1	0.04	73
0	0.04	72.8
0.1	0.04	73.2
0.2	0.04	73.4
0.3	0.04	73.6
0.4	0.04	74.4
0.5	0.04	75
0.6	0.04	75
0.7	0.04	75.4
0.8	0.04	75.8
0.9	0.04	75.8
1	0.04	76
1.2	0.04	76.8
1.4	0.04	78.6
1.6	0.04	79.6
1.8	0.04	80.2
2	0.04	81.2
-2	0.05	65
-1.8	0.05	65.6
-1.4	0.05	68.8
-1.2	0.05	69.6
-1	0.05	71.2
-0.9	0.05	71.4
-0.8	0.05	71.8
-0.7	0.05	72.2
-0.6	0.05	72.2
-0.5	0.05	72.8
-0.4	0.05	72.8
-0.3	0.05	73.2
-0.2	0.05	73.4
-0.1	0.05	73.6
0	0.05	73.8
0.1	0.05	75
0.2	0.05	74.8
0.3	0.05	75
0.4	0.05	75.8
0.5	0.05	75.8
0.6	0.05	76.2
0.7	0.05	76.8
0.8	0.05	76.8
0.9	0.05	77.8
1	0.05	78.6
1.2	0.05	79.6
1.4	0.05	80.2
1.6	0.05	81.2
1.8	0.05	81.2
2	0.05	81.8
-2	0.06	67.4
-1.8	0.06	68.8
-1.4	0.06	71.2
-1.2	0.06	71.6
-1	0.06	72.2
-0.9	0.06	72.8
-0.8	0.06	72.8
-0.7	0.06	73.2
-0.6	0.06	73.4
-0.5	0.06	73.6
-0.4	0.06	73.8
-0.3	0.06	74.4
-0.2	0.06	74.8
-0.1	0.06	75
0	0.06	75.4
0.1	0.06	75.8
0.2	0.06	76.2
0.3	0.06	76.8
0.4	0.06	76.8
0.5	0.06	77.8
0.6	0.06	78.6
0.7	0.06	79
0.8	0.06	79.8
0.9	0.06	79.8
1	0.06	80.4
1.2	0.06	81.2
1.4	0.06	81.2
1.6	0.06	81.8
1.8	0.06	82.6
2	0.06	83
-2	0.07	69.4
-1.8	0.07	71.2
-1.4	0.07	72.2
-1.2	0.07	72.8
-1	0.07	73.4
-0.9	0.07	73.6
-0.8	0.07	73.8
-0.7	0.07	74.4
-0.6	0.07	74.8
-0.5	0.07	74.8
-0.4	0.07	75.4
-0.3	0.07	75.8
-0.2	0.07	76.2
-0.1	0.07	76.8
0	0.07	76.8
0.1	0.07	77.8
0.2	0.07	78.6
0.3	0.07	79
0.4	0.07	79.6
0.5	0.07	79.8
0.6	0.07	80.2
0.7	0.07	80.8
0.8	0.07	81.2
0.9	0.07	81.2
1	0.07	81.2
1.2	0.07	81.8
1.4	0.07	82.2
1.6	0.07	82.8
1.8	0.07	84.2
2	0.07	84.4
-2	0.08	71.6
-1.8	0.08	72.2
-1.4	0.08	73.4
-1.2	0.08	73.8
-1	0.08	74.8
-0.9	0.08	74.8
-0.8	0.08	75.4
-0.7	0.08	75.8
-0.6	0.08	76.2
-0.5	0.08	76.8
-0.4	0.08	76.8
-0.3	0.08	77.8
-0.2	0.08	78.6
-0.1	0.08	79
0	0.08	79.8
0.1	0.08	80
0.2	0.08	80.2
0.3	0.08	80.8
0.4	0.08	81.2
0.5	0.08	81.4
0.6	0.08	81.4
0.7	0.08	81.6
0.8	0.08	82
0.9	0.08	82.4
1	0.08	82.4
1.2	0.08	82.8
1.4	0.08	83.8
1.6	0.08	84.4
1.8	0.08	84.6
2	0.08	85
-2	0.09	72.8
-1.8	0.09	74
-1.4	0.09	74.8
-1.2	0.09	75.2
-1	0.09	76.2
-0.9	0.09	76.8
-0.8	0.09	76.8
-0.7	0.09	77.8
-0.6	0.09	78.6
-0.5	0.09	79
-0.4	0.09	79.6
-0.3	0.09	79.8
-0.2	0.09	80.2
-0.1	0.09	81
0	0.09	81.4
0.1	0.09	81.4
0.2	0.09	81.4
0.3	0.09	81.6
0.4	0.09	82.2
0.5	0.09	82.6
0.6	0.09	82.6
0.7	0.09	82.6
0.8	0.09	83
0.9	0.09	83.6
1	0.09	84
1.2	0.09	84.4
1.4	0.09	81.8
1.6	0.09	85
1.8	0.09	85.6
2	0.09	86.6
-2	0.1	73.8
-1.8	0.1	74.8
-1.4	0.1	76.6
-1.2	0.1	77.2
-1	0.1	78.6
-0.9	0.1	79
-0.8	0.1	79.6
-0.7	0.1	80
-0.6	0.1	80.4
-0.5	0.1	81
-0.4	0.1	81.2
-0.3	0.1	81.4
-0.2	0.1	81.4
-0.1	0.1	81.6
0	0.1	82
0.1	0.1	82.6
0.2	0.1	82.6
0.3	0.1	82.8
0.4	0.1	83
0.5	0.1	83.6
0.6	0.1	83.6
0.7	0.1	84.2
0.8	0.1	84.4
0.9	0.1	84.8
1	0.1	85
1.2	0.1	85.4
1.4	0.1	85.6
1.6	0.1	87
1.8	0.1	84.4
2	0.1	85
-2	0.12	77.2
-1.8	0.12	78.6
-1.4	0.12	80.4
-1.2	0.12	81.2
-1	0.12	81.4
-0.9	0.12	81.6
-0.8	0.12	82
-0.7	0.12	82.4
-0.6	0.12	82.4
-0.5	0.12	82.8
-0.4	0.12	83
-0.3	0.12	83.6
-0.2	0.12	83.6
-0.1	0.12	84.2
0	0.12	84.4
0.1	0.12	84.6
0.2	0.12	85
0.3	0.12	85.2
0.4	0.12	85.4
0.5	0.12	82.6
0.6	0.12	85.6
0.7	0.12	83.4
0.8	0.12	83.6
0.9	0.12	84.4
1	0.12	84.4
1.2	0.12	85
1.4	0.12	86
1.6	0.12	86.6
1.8	0.12	87.2
2	0.12	89.6
-2	0.14	81.4
-1.8	0.14	81.4
-1.4	0.14	82.4
-1.2	0.14	83
-1	0.14	83.6
-0.9	0.14	84.2
-0.8	0.14	84.4
-0.7	0.14	84.4
-0.6	0.14	84.6
-0.5	0.14	85.2
-0.4	0.14	85.4
-0.3	0.14	82.6
-0.2	0.14	85.8
-0.1	0.14	83.4
0	0.14	83.6
0.1	0.14	84
0.2	0.14	84.4
0.3	0.14	84.6
0.4	0.14	85.2
0.5	0.14	86
0.6	0.14	86
0.7	0.14	86.2
0.8	0.14	86.6
0.9	0.14	87.2
1	0.14	87.2
1.2	0.14	88.6
1.4	0.14	89.2
1.6	0.14	89.6
1.8	0.14	90.2
2	0.14	90
-2	0.16	83
-1.8	0.16	83.6
-1.4	0.16	84.6
-1.2	0.16	85.4
-1	0.16	85.8
-0.9	0.16	83.6
-0.8	0.16	86.8
-0.7	0.16	83.8
-0.6	0.16	84.6
-0.5	0.16	84.8
-0.4	0.16	85
-0.3	0.16	85.4
-0.2	0.16	86
-0.1	0.16	86.2
0	0.16	86.6
0.1	0.16	87.2
0.2	0.16	87.2
0.3	0.16	87.4
0.4	0.16	87.6
0.5	0.16	88
0.6	0.16	89.2
0.7	0.16	89.4
0.8	0.16	89.4
0.9	0.16	90
1	0.16	89.2
1.2	0.16	90.4
1.4	0.16	91.4
1.6	0.16	91.6
1.8	0.16	92.6
2	0.16	92.6
-2	0.18	85
-1.8	0.18	85.8
-1.4	0.18	84.4
-1.2	0.18	85.2
-1	0.18	86
-0.9	0.18	86.2
-0.8	0.18	86.6
-0.7	0.18	87.2
-0.6	0.18	87.2
-0.5	0.18	87.4
-0.4	0.18	87.6
-0.3	0.18	88
-0.2	0.18	88.2
-0.1	0.18	88.4
0	0.18	90.2
0.1	0.18	89
0.2	0.18	89.2
0.3	0.18	90.8
0.4	0.18	90.4
0.5	0.18	90.8
0.6	0.18	91.4
0.7	0.18	91.4
0.8	0.18	91.2
0.9	0.18	92
1	0.18	92.6
1.2	0.18	92.4
1.4	0.18	93.4
1.6	0.18	94
1.8	0.18	94.4
2	0.18	94.6
-2	0.2	85.2
-1.8	0.2	86
-1.4	0.2	87.2
-1.2	0.2	91
-1	0.2	88.2
-0.9	0.2	88.4
-0.8	0.2	89.4
-0.7	0.2	89
-0.6	0.2	89.2
-0.5	0.2	90.8
-0.4	0.2	90.4
-0.3	0.2	90.8
-0.2	0.2	91.4
-0.1	0.2	91.4
0	0.2	91.2
0.1	0.2	91.8
0.2	0.2	92
0.3	0.2	92.8
0.4	0.2	92.6
0.5	0.2	93
0.6	0.2	93.4
0.7	0.2	93.8
0.8	0.2	93.8
0.9	0.2	94.2
1	0.2	94.6
1.2	0.2	94.6
1.4	0.2	94.6
1.6	0.2	95
1.8	0.2	95
2	0.2	95.6
-2	0.22	87.6
-1.8	0.22	88.2
-1.4	0.22	89.2
-1.2	0.22	90.4
-1	0.22	91.4
-0.9	0.22	91.4
-0.8	0.22	91.6
-0.7	0.22	91.8
-0.6	0.22	91.8
-0.5	0.22	92.8
-0.4	0.22	92.4
-0.3	0.22	93.2
-0.2	0.22	93.4
-0.1	0.22	93.6
0	0.22	93.8
0.1	0.22	94.2
0.2	0.22	94.4
0.3	0.22	94.6
0.4	0.22	94.6
0.5	0.22	94.6
0.6	0.22	94.6
0.7	0.22	94.8
0.8	0.22	95
0.9	0.22	95
1	0.22	95.2
1.2	0.22	95.6
1.4	0.22	96.2
1.6	0.22	97.2
1.8	0.22	97.4
2	0.22	97.8
-2	0.24	90.8
-1.8	0.24	91.4
-1.4	0.24	92.4
-1.2	0.24	92.6
-1	0.24	93.4
-0.9	0.24	93.8
-0.8	0.24	93.8
-0.7	0.24	94.2
-0.6	0.24	94.4
-0.5	0.24	94.6
-0.4	0.24	94.6
-0.3	0.24	94.6
-0.2	0.24	94.6
-0.1	0.24	94.8
0	0.24	95
0.1	0.24	95
0.2	0.24	95
0.3	0.24	95.4
0.4	0.24	95.6
0.5	0.24	95.6
0.6	0.24	96.2
0.7	0.24	96.4
0.8	0.24	97.2
0.9	0.24	97.4
1	0.24	97.4
1.2	0.24	97.8
1.4	0.24	97.8
1.6	0.24	97.8
1.8	0.24	97.8
2	0.24	97.8
-2	0.26	92.6
-1.8	0.26	93.4
-1.4	0.26	94.4
-1.2	0.26	94.6
-1	0.26	94.6
-0.9	0.26	94.8
-0.8	0.26	95
-0.7	0.26	95
-0.6	0.26	95
-0.5	0.26	95.4
-0.4	0.26	95.6
-0.3	0.26	95.6
-0.2	0.26	96.2
-0.1	0.26	96.4
0	0.26	97.2
0.1	0.26	97.4
0.2	0.26	97.4
0.3	0.26	97.6
0.4	0.26	97.8
0.5	0.26	97.8
0.6	0.26	97.8
0.7	0.26	97.8
0.8	0.26	97.8
0.9	0.26	97.8
1	0.26	97.8
1.2	0.26	97.8
1.4	0.26	98.4
1.6	0.26	98.4
1.8	0.26	98.4
2	0.26	98.4
-2	0.28	94.6
-1.8	0.28	94.6
-1.4	0.28	95
-1.2	0.28	95.6
-1	0.28	96.2
-0.9	0.28	96.4
-0.8	0.28	97.2
-0.7	0.28	97.4
-0.6	0.28	97.4
-0.5	0.28	97.6
-0.4	0.28	97.8
-0.3	0.28	97.8
-0.2	0.28	97.8
-0.1	0.28	97.8
0	0.28	97.8
0.1	0.28	97.8
0.2	0.28	97.8
0.3	0.28	97.8
0.4	0.28	97.8
0.5	0.28	98
0.6	0.28	98.4
0.7	0.28	98.4
0.8	0.28	98.4
0.9	0.28	98.4
1	0.28	98.4
1.2	0.28	98.4
1.4	0.28	98.4
1.6	0.28	98.4
1.8	0.28	98.6
2	0.28	98.8
-2	0.3	95.6
-1.8	0.3	96.2
-1.4	0.3	97.4
-1.2	0.3	97.8
-1	0.3	97.8
-0.9	0.3	97.8
-0.8	0.3	97.8
-0.7	0.3	97.8
-0.6	0.3	97.8
-0.5	0.3	97.8
-0.4	0.3	97.8
-0.3	0.3	98
-0.2	0.3	98.4
-0.1	0.3	98.4
0	0.3	98.4
0.1	0.3	98.4
0.2	0.3	98.4
0.3	0.3	98.4
0.4	0.3	98.4
0.5	0.3	98.4
0.6	0.3	98.4
0.7	0.3	98.4
0.8	0.3	98.4
0.9	0.3	98.6
1	0.3	98.6
1.2	0.3	98.8
1.4	0.3	98.8
1.6	0.3	98.8
1.8	0.3	98.8
2	0.3	98.8
-2	0.32	97.8
-1.8	0.32	97.8
-1.4	0.32	97.8
-1.2	0.32	97.8
-1	0.32	98.4
-0.9	0.32	98.2
-0.8	0.32	98.4
-0.7	0.32	98.4
-0.6	0.32	98.4
-0.5	0.32	98.4
-0.4	0.32	98.4
-0.3	0.32	98.4
-0.2	0.32	98.4
-0.1	0.32	98.4
0	0.32	98.4
0.1	0.32	98.6
0.2	0.32	98.6
0.3	0.32	98.8
0.4	0.32	98.8
0.5	0.32	98.8
0.6	0.32	98.8
0.7	0.32	98.8
0.8	0.32	98.8
0.9	0.32	98.8
1	0.32	98.8
1.2	0.32	98.8
1.4	0.32	99
1.6	0.32	99.2
1.8	0.32	99.4
2	0.32	99.6
-2	0.34	97.8
-1.8	0.34	98.4
-1.4	0.34	98.4
-1.2	0.34	98.4
-1	0.34	98.4
-0.9	0.34	98.4
-0.8	0.34	98.4
-0.7	0.34	98.6
-0.6	0.34	98.6
-0.5	0.34	98.8
-0.4	0.34	98.8
-0.3	0.34	98.8
-0.2	0.34	98.8
-0.1	0.34	98.8
0	0.34	98.8
0.1	0.34	98.8
0.2	0.34	98.8
0.3	0.34	98.8
0.4	0.34	98.8
0.5	0.34	98.8
0.6	0.34	99
0.7	0.34	99.2
0.8	0.34	99.2
0.9	0.34	99.2
1	0.34	99.4
1.2	0.34	99.6
1.4	0.34	99.6
1.6	0.34	99.6
1.8	0.34	99.6
2	0.34	99.8
-2	0.36	98.4
-1.8	0.36	98.4
-1.4	0.36	98.6
-1.2	0.36	98.8
-1	0.36	98.8
-0.9	0.36	98.8
-0.8	0.36	98.8
-0.7	0.36	98.8
-0.6	0.36	98.8
-0.5	0.36	98.8
-0.4	0.36	98.8
-0.3	0.36	98.8
-0.2	0.36	99
-0.1	0.36	99.2
0	0.36	99.2
0.1	0.36	99.2
0.2	0.36	99.4
0.3	0.36	99.6
0.4	0.36	99.6
0.5	0.36	99.6
0.6	0.36	99.6
0.7	0.36	99.6
0.8	0.36	99.6
0.9	0.36	99.6
1	0.36	99.6
1.2	0.36	99.8
1.4	0.36	99.8
1.6	0.36	99.8
1.8	0.36	99.8
2	0.36	99.8
-2	0.38	98.8
-1.8	0.38	98.8
-1.4	0.38	98.8
-1.2	0.38	98.8
-1	0.38	99
-0.9	0.38	99.2
-0.8	0.38	99.2
-0.7	0.38	99.2
-0.6	0.38	99.4
-0.5	0.38	99.6
-0.4	0.38	99.6
-0.3	0.38	99.6
-0.2	0.38	99.6
-0.1	0.38	99.6
0	0.38	99.6
0.1	0.38	99.6
0.2	0.38	99.6
0.3	0.38	99.6
0.4	0.38	99.8
0.5	0.38	99.8
0.6	0.38	99.8
0.7	0.38	99.8
0.8	0.38	99.8
0.9	0.38	99.8
1	0.38	99.8
1.2	0.38	99.8
1.4	0.38	99.8
1.6	0.38	99.8
1.8	0.38	99.8
2	0.38	99.8
-2	0.4	99.2
-1.8	0.4	99
-1.4	0.4	99.4
-1.2	0.4	99.6
-1	0.4	99.6
-0.9	0.4	99.6
-0.8	0.4	99.6
-0.7	0.4	99.6
-0.6	0.4	99.6
-0.5	0.4	99.6
-0.4	0.4	99.8
-0.3	0.4	99.8
-0.2	0.4	99.8
-0.1	0.4	99.8
0	0.4	99.8
0.1	0.4	99.8
0.2	0.4	99.8
0.3	0.4	99.8
0.4	0.4	99.8
0.5	0.4	99.8
0.6	0.4	99.8
0.7	0.4	99.8
0.8	0.4	99.8
0.9	0.4	99.8
1	0.4	99.8
1.2	0.4	99.8
1.4	0.4	99.8
1.6	0.4	99.8
1.8	0.4	100
2	0.4	100
}\txtwithout
\noindent\begin{figure}[h!]
\captionsetup{justification=centering}
\centering
\scriptsize
\begin{subfigure}{0.45\textwidth}
    \centering
    \begin{tikzpicture}[scale=.75]
    \begin{axis}[
        domain=-2:2,          
        y domain=-0.2:0.2,   
        samples=30,         
        colormap/RdYlBu,
        xtick={-1,0,+1}, 
        xticklabels={
            $-1$,
            $\beta^{time}$,
            $+1$},
        ytick={-.1,0,.1},
        yticklabels={
            $-0.1$,
            $\beta^{\mathrm{price}}$,
            $+0.1$},
        zmin=0,
        zmax=100,
        zlabel = {Coverage rate(\%)}
    ]
    
    \addplot3[surf, fill= white]
        {1 + 6 * exp((x + 6.05 + 40 * (y - 0.1))) > 0 
        ? 100 * (1 - 1 / (1 + 6 * exp((x + 6.05 + 40 * (y - 0.1))))) 
        : 0};
    
    \addplot3[
        thick,
        color=black!70,
        samples=20,
        domain=-2:2
    ]
    ({x}, 
     {(0.1 - (x + 7.43) / 40)}, 
     {20});
    
    \addplot3[
        surf,
        opacity=0.5,    
        shader=interp,  
        draw=none,      
        colormap={customgrey}{
            rgb255(0pt)=(156,156,156)
            rgb255(100pt)=(156,156,156)
        }
    ]
    ({x}, {y}, {20});
    \end{axis}

    \end{tikzpicture}

    \caption{Without assortment policy.}
    \label{fig:coverage_without}
    \end{subfigure}
\noindent\begin{subfigure}{0.45\textwidth}
    \centering
        \begin{tikzpicture}[scale=.75]
        \begin{axis}[
            colormap/RdYlBu,
            xtick={-1,0,+1}, 
            xticklabels={
                $-1$,
                $\beta^{time}$,
                $+1$},
            ytick={-.1,0,.1,.2}, 
            yticklabels={
                $-0.1$,
                $\beta^{\mathrm{price}}$,
                $+0.1$,
                $+0.2$},
            zlabel = {Coverage rate(\%)},
            zmin=0,
            zmax=100,
        ]
        \addplot3
        [surf, fill= white, mesh/rows=40, mesh/cols=30]
        table [x index=0, y index=1, z index=2]\txtwithout;

        \addplot3[
            thick,
            color=black!70,
            mark options={scale=0.5, fill=black}
        ]
        table [
            x index=0,
            y index=1,
            z expr={(abs(\thisrow{2}-20) <= 0.9) ? \thisrow{2} : nan}  
        ]{\txtwithout};      
        
        \addplot3[
            surf,
            opacity=0.5,     
            draw=none,       
            shader=interp,   
            domain=-2:2,
            y domain=-0.2:0.4,
            samples=30,
            colormap={customgrey}{
                rgb255(0pt)=(156,156,156)
                rgb255(100pt)=(156,156,156)
                }
        ]
        ({x}, {y}, {20});
        \end{axis} 
        \end{tikzpicture}
    \caption{With assortment policy.}
    \label{fig:coverage_with-without_fairness}
    \end{subfigure}

\caption{Coverage rate of customers with and without assortment decisions generated by the sALNS method.}
\label{fig:coverage}
\end{figure}

\pgfplotstableread{
-2	-0.2	11.00893388
-1.8	-0.2	10.16448569
-1.4	-0.2	10.9258557
-1.2	-0.2	11.205956
-1	-0.2	10.93549243
-0.9	-0.2	10.67309836
-0.8	-0.2	11.58944113
-0.7	-0.2	11.57578812
-0.6	-0.2	10.85952922
-0.5	-0.2	12.80534289
-0.4	-0.2	12.87081759
-0.3	-0.2	13.60346018
-0.2	-0.2	14.22540561
-0.1	-0.2	15.05728762
0	-0.2	13.39948484
0.1	-0.2	14.23849388
0.2	-0.2	15.4556076
0.3	-0.2	15.0303003
0.4	-0.2	12.52874852
0.5	-0.2	13.67882133
0.6	-0.2	12.11849156
0.7	-0.2	14.40112564
0.8	-0.2	14.93738364
0.9	-0.2	15.7950989
1	-0.2	15.3340267
1.2	-0.2	16.32545789
1.4	-0.2	14.99563928
1.6	-0.2	14.08778858
1.8	-0.2	13.9747691
2	-0.2	14.06065057
-2	-0.18	10.80237376
-1.8	-0.18	10.96844065
-1.4	-0.18	11.56722499
-1.2	-0.18	12.11435136
-1	-0.18	13.54200377
-0.9	-0.18	13.51734168
-0.8	-0.18	14.01394129
-0.7	-0.18	13.74945514
-0.6	-0.18	12.89844096
-0.5	-0.18	14.80759269
-0.4	-0.18	13.71347617
-0.3	-0.18	13.32009207
-0.2	-0.18	12.80765921
-0.1	-0.18	11.77075111
0	-0.18	13.19158454
0.1	-0.18	13.85970627
0.2	-0.18	13.59130474
0.3	-0.18	14.0166833
0.4	-0.18	14.28545798
0.5	-0.18	15.19695963
0.6	-0.18	13.78810367
0.7	-0.18	14.85145549
0.8	-0.18	14.90217717
0.9	-0.18	13.98644211
1	-0.18	15.23687899
1.2	-0.18	14.93311568
1.4	-0.18	13.58367907
1.6	-0.18	12.13086861
1.8	-0.18	11.54399547
2	-0.18	11.70123448
-2	-0.16	11.41994132
-1.8	-0.16	12.03860562
-1.4	-0.16	13.32321429
-1.2	-0.16	12.23635537
-1	-0.16	11.78968686
-0.9	-0.16	12.72203098
-0.8	-0.16	13.31022894
-0.7	-0.16	13.97201552
-0.6	-0.16	13.66422597
-0.5	-0.16	13.70809231
-0.4	-0.16	15.47857578
-0.3	-0.16	15.35424797
-0.2	-0.16	15.16024808
-0.1	-0.16	14.90456862
0	-0.16	13.15473197
0.1	-0.16	13.85264949
0.2	-0.16	14.64663794
0.3	-0.16	14.66860585
0.4	-0.16	14.05054583
0.5	-0.16	14.11894047
0.6	-0.16	12.81026514
0.7	-0.16	13.10635829
0.8	-0.16	11.47750956
0.9	-0.16	11.55806988
1	-0.16	11.52096561
1.2	-0.16	11.61487363
1.4	-0.16	11.79100291
1.6	-0.16	12.76325471
1.8	-0.16	13.05257425
2	-0.16	13.82408502
-2	-0.14	13.49216055
-1.8	-0.14	13.65532851
-1.4	-0.14	12.60729185
-1.2	-0.14	14.30728869
-1	-0.14	15.26846774
-0.9	-0.14	15.63288991
-0.8	-0.14	14.76734135
-0.7	-0.14	14.1229376
-0.6	-0.14	14.0199829
-0.5	-0.14	14.68139758
-0.4	-0.14	15.35547123
-0.3	-0.14	14.88617602
-0.2	-0.14	13.73183159
-0.1	-0.14	11.80939589
0	-0.14	11.86436981
0.1	-0.14	10.77880436
0.2	-0.14	11.89985788
0.3	-0.14	11.43720648
0.4	-0.14	11.7387472
0.5	-0.14	11.8823728
0.6	-0.14	12.09098967
0.7	-0.14	12.10119146
0.8	-0.14	12.49999924
0.9	-0.14	12.64977346
1	-0.14	12.71037829
1.2	-0.14	13.17747966
1.4	-0.14	13.59039939
1.6	-0.14	14.81416275
1.8	-0.14	14.85772918
2	-0.14	15.22826369
-2	-0.12	13.97609398
-1.8	-0.12	15.12425897
-1.4	-0.12	13.47578549
-1.2	-0.12	14.33605902
-1	-0.12	14.25135531
-0.9	-0.12	13.02740985
-0.8	-0.12	12.26762774
-0.7	-0.12	11.95829642
-0.6	-0.12	12.05516034
-0.5	-0.12	13.31613974
-0.4	-0.12	11.71275579
-0.3	-0.12	12.27465473
-0.2	-0.12	12.17888804
-0.1	-0.12	11.43650872
0	-0.12	12.23361873
0.1	-0.12	12.40129691
0.2	-0.12	12.45936687
0.3	-0.12	12.04270404
0.4	-0.12	12.43101907
0.5	-0.12	13.24275299
0.6	-0.12	13.87590092
0.7	-0.12	14.30179263
0.8	-0.12	13.04297417
0.9	-0.12	14.10528146
1	-0.12	14.21871086
1.2	-0.12	14.48864169
1.4	-0.12	5.338825292
1.6	-0.12	4.946440769
1.8	-0.12	4.928820307
2	-0.12	5.664478718
-2	-0.1	14.40867979
-1.8	-0.1	13.65364084
-1.4	-0.1	11.67969335
-1.2	-0.1	11.10250235
-1	-0.1	11.6811353
-0.9	-0.1	14.0353371
-0.8	-0.1	11.8701894
-0.7	-0.1	12.0764123
-0.6	-0.1	6.1337695
-0.5	-0.1	5.9308744
-0.4	-0.1	13.2412816
-0.3	-0.1	14.0270537
-0.2	-0.1	14.1411766
-0.1	-0.1	6.4666417
0	-0.1	6.6935108
0.1	-0.1	6.2747569
0.2	-0.1	6.4496918
0.3	-0.1	11.8693626
0.4	-0.1	11.4109776
0.5	-0.1	12.1469446
0.6	-0.1	12.4737869
0.7	-0.1	13.1634207
0.8	-0.1	13.2978634
0.9	-0.1	13.5316918
1	-0.1	6.9768524
1.2	-0.1	5.304871636
1.4	-0.1	6.391224876
1.6	-0.1	7.617638197
1.8	-0.1	6.776475926
2	-0.1	7.066211897
-2	-0.09	12.2437141
-1.8	-0.09	11.27276177
-1.4	-0.09	11.99571539
-1.2	-0.09	11.25853036
-1	-0.09	12.947747
-0.9	-0.09	6.6866447
-0.8	-0.09	5.5490381
-0.7	-0.09	12.4017895
-0.6	-0.09	12.2313073
-0.5	-0.09	13.7606578
-0.4	-0.09	13.2374005
-0.3	-0.09	14.04756
-0.2	-0.09	13.1382112
-0.1	-0.09	14.178053
0	-0.09	6.074108
0.1	-0.09	13.0308984
0.2	-0.09	13.5923431
0.3	-0.09	6.0735134
0.4	-0.09	6.8989685
0.5	-0.09	6.2992771
0.6	-0.09	6.0148712
0.7	-0.09	6.9872685
0.8	-0.09	6.1102962
0.9	-0.09	6.2912481
1	-0.09	5.4916957
1.2	-0.09	7.205169309
1.4	-0.09	7.480120386
1.6	-0.09	6.388384218
1.8	-0.09	6.842762496
2	-0.09	7.282082639
-2	-0.08	11.45506969
-1.8	-0.08	11.80585903
-1.4	-0.08	12.34179671
-1.2	-0.08	11.9079469
-1	-0.08	13.0408668
-0.9	-0.08	12.8927357
-0.8	-0.08	14.1996898
-0.7	-0.08	6.8689571
-0.6	-0.08	7.9128249
-0.5	-0.08	7.7142527
-0.4	-0.08	6.5770288
-0.3	-0.08	6.4707941
-0.2	-0.08	7.186171
-0.1	-0.08	5.9610753
0	-0.08	6.1648649
0.1	-0.08	7.901989
0.2	-0.08	6.4564845
0.3	-0.08	6.8461374
0.4	-0.08	6.6373698
0.5	-0.08	6.8968794
0.6	-0.08	7.474247
0.7	-0.08	13.4360883
0.8	-0.08	13.9104061
0.9	-0.08	7.4682699
1	-0.08	7.2925127
1.2	-0.08	6.743306871
1.4	-0.08	7.542976806
1.6	-0.08	7.251177542
1.8	-0.08	8.400482368
2	-0.08	7.857157931
-2	-0.07	11.50811149
-1.8	-0.07	11.10032368
-1.4	-0.07	13.21844635
-1.2	-0.07	13.23072979
-1	-0.07	6.4293647
-0.9	-0.07	8.2585437
-0.8	-0.07	7.7602937
-0.7	-0.07	7.7802047
-0.6	-0.07	8.1154817
-0.5	-0.07	7.3640453
-0.4	-0.07	7.2802975
-0.3	-0.07	7.0313284
-0.2	-0.07	7.6608794
-0.1	-0.07	7.4772161
0	-0.07	7.4712381
0.1	-0.07	7.2159857
0.2	-0.07	7.0939403
0.3	-0.07	7.7116893
0.4	-0.07	7.4972744
0.5	-0.07	6.9990437
0.6	-0.07	7.0251563
0.7	-0.07	6.9561739
0.8	-0.07	7.6162899
0.9	-0.07	7.8329768
1	-0.07	8.1567127
1.2	-0.07	8.024716721
1.4	-0.07	8.982380485
1.6	-0.07	8.515792913
1.8	-0.07	9.979989643
2	-0.07	11.2355776
-2	-0.06	11.66402095
-1.8	-0.06	12.94612884
-1.4	-0.06	6.144646963
-1.2	-0.06	7.057511032
-1	-0.06	7.5905814
-0.9	-0.06	7.6110537
-0.8	-0.06	7.3767443
-0.7	-0.06	7.5852637
-0.6	-0.06	6.8508311
-0.5	-0.06	8.0704788
-0.4	-0.06	7.8792082
-0.3	-0.06	5.5604839
-0.2	-0.06	7.5771574
-0.1	-0.06	7.7056181
0	-0.06	7.6610743
0.1	-0.06	8.0853223
0.2	-0.06	7.8298115
0.3	-0.06	7.435412
0.4	-0.06	7.7682315
0.5	-0.06	8.6813163
0.6	-0.06	7.5135398
0.7	-0.06	7.9431637
0.8	-0.06	8.067447
0.9	-0.06	8.0967851
1	-0.06	8.314383
1.2	-0.06	8.583520365
1.4	-0.06	9.765583046
1.6	-0.06	11.20872649
1.8	-0.06	10.90460873
2	-0.06	11.69840092
-2	-0.05	6.305150005
-1.8	-0.05	6.213454551
-1.4	-0.05	6.951072567
-1.2	-0.05	7.596906235
-1	-0.05	7.8300995
-0.9	-0.05	7.456259
-0.8	-0.05	7.1929436
-0.7	-0.05	9.5259694
-0.6	-0.05	8.4218983
-0.5	-0.05	7.7993495
-0.4	-0.05	7.977863
-0.3	-0.05	7.5847257
-0.2	-0.05	8.2618014
-0.1	-0.05	8.0196847
0	-0.05	8.4365895
0.1	-0.05	7.5883924
0.2	-0.05	7.7940157
0.3	-0.05	8.0334562
0.4	-0.05	8.6893256
0.5	-0.05	8.006497
0.6	-0.05	7.0436791
0.7	-0.05	7.2163131
0.8	-0.05	7.0151285
0.9	-0.05	8.7144002
1	-0.05	9.6843058
1.2	-0.05	10.66345257
1.4	-0.05	11.35318963
1.6	-0.05	11.65592136
1.8	-0.05	12.07749259
2	-0.05	14.07860658
-2	-0.04	7.304045317
-1.8	-0.04	7.49164114
-1.4	-0.04	6.947520125
-1.2	-0.04	7.42744155
-1	-0.04	7.1974176
-0.9	-0.04	8.0231362
-0.8	-0.04	6.8846784
-0.7	-0.04	11.0302978
-0.6	-0.04	10.8329003
-0.5	-0.04	7.4959732
-0.4	-0.04	8.0558114
-0.3	-0.04	7.495962
-0.2	-0.04	7.7210572
-0.1	-0.04	7.9737017
0	-0.04	8.8872551
0.1	-0.04	7.5628638
0.2	-0.04	7.9797052
0.3	-0.04	10.2669325
0.4	-0.04	9.0518891
0.5	-0.04	8.3086577
0.6	-0.04	8.0314537
0.7	-0.04	8.745283
0.8	-0.04	7.8863857
0.9	-0.04	8.0941709
1	-0.04	11.3531896
1.2	-0.04	11.18208361
1.4	-0.04	11.90765276
1.6	-0.04	14.07860658
1.8	-0.04	16.14602585
2	-0.04	17.68270701
-2	-0.03	7.215594383
-1.8	-0.03	8.4395533
-1.4	-0.03	7.382783861
-1.2	-0.03	8.542485784
-1	-0.03	8.9173982
-0.9	-0.03	7.679983
-0.8	-0.03	12.134737
-0.7	-0.03	12.4819562
-0.6	-0.03	8.5599987
-0.5	-0.03	8.1542257
-0.4	-0.03	11.8218081
-0.3	-0.03	11.7055683
-0.2	-0.03	8.2686296
-0.1	-0.03	8.2919451
0	-0.03	11.2179852
0.1	-0.03	11.3820136
0.2	-0.03	8.3704414
0.3	-0.03	9.0329364
0.4	-0.03	9.3805372
0.5	-0.03	9.1161338
0.6	-0.03	8.8788954
0.7	-0.03	9.7197044
0.8	-0.03	9.7178742
0.9	-0.03	8.2577115
1	-0.03	13.2402628
1.2	-0.03	13.86323108
1.4	-0.03	16.96528629
1.6	-0.03	17.6519722
1.8	-0.03	18.99023545
2	-0.03	17.97402076
-2	-0.02	7.818613168
-1.8	-0.02	7.812868649
-1.4	-0.02	8.421923935
-1.2	-0.02	8.046924446
-1	-0.02	8.1378503
-0.9	-0.02	15.1535143
-0.8	-0.02	14.4292863
-0.7	-0.02	8.5256592
-0.6	-0.02	8.4910517
-0.5	-0.02	10.3984092
-0.4	-0.02	8.8146563
-0.3	-0.02	8.50589
-0.2	-0.02	9.3047641
-0.1	-0.02	8.5125082
0	-0.02	10.4191983
0.1	-0.02	10.3561847
0.2	-0.02	8.4273598
0.3	-0.02	12.5789662
0.4	-0.02	12.7095089
0.5	-0.02	12.6982751
0.6	-0.02	11.8163921
0.7	-0.02	11.665273
0.8	-0.02	11.5744789
0.9	-0.02	11.4438329
1	-0.02	16.8032538
1.2	-0.02	17.63592176
1.4	-0.02	18.8258826
1.6	-0.02	19.23431696
1.8	-0.02	17.88470299
2	-0.02	19.22842132
-2	-0.01	6.896271442
-1.8	-0.01	8.847618456
-1.4	-0.01	8.454135907
-1.2	-0.01	8.711445011
-1	-0.01	9.0345719
-0.9	-0.01	9.0730271
-0.8	-0.01	13.0927023
-0.7	-0.01	12.4379195
-0.6	-0.01	9.0050695
-0.5	-0.01	9.4613001
-0.4	-0.01	12.4112099
-0.3	-0.01	11.2200016
-0.2	-0.01	15.6364024
-0.1	-0.01	14.1169965
0	-0.01	12.671228
0.1	-0.01	11.484179
0.2	-0.01	16.5816357
0.3	-0.01	18.2863104
0.4	-0.01	11.8346099
0.5	-0.01	18.1981305
0.6	-0.01	18.1208187
0.7	-0.01	8.6456449
0.8	-0.01	9.1615309
0.9	-0.01	11.9113054
1	-0.01	18.4308458
1.2	-0.01	18.6292141
1.4	-0.01	19.68015379
1.6	-0.01	18.29535622
1.8	-0.01	20.43874966
2	-0.01	21.61740665
-2	0	7.813101462
-1.8	0	8.399838612
-1.4	0	9.546820621
-1.2	0	9.673702452
-1	0	10.4960599
-0.9	0	10.5724274
-0.8	0	11.5937575
-0.7	0	19.9653418
-0.6	0	20.0189618
-0.5	0	19.1381288
-0.4	0	19.2607037
-0.3	0	14.6837805
-0.2	0	13.4498809
-0.1	0	18.6002357
0	0	19.4336359
0.1	0	12.8104091
0.2	0	12.9528306
0.3	0	16.5816357
0.4	0	14.8624116
0.5	0	18.3892623
0.6	0	12.6840645
0.7	0	12.841648
0.8	0	16.3465411
0.9	0	13.0949147
1	0	18.3761804
1.2	0	18.37679397
1.4	0	20.43874966
1.6	0	21.62622244
1.8	0	21.3701936
2	0	21.86075798
-2	0.01	9.397608979
-1.8	0.01	9.846666908
-1.4	0.01	9.679173661
-1.2	0.01	12.01482401
-1	0.01	13.4662745
-0.9	0.01	20.6609602
-0.8	0.01	19.0839555
-0.7	0.01	20.1749082
-0.6	0.01	20.5845738
-0.5	0.01	13.2647325
-0.4	0.01	20.5265836
-0.3	0.01	18.5175917
-0.2	0.01	13.5451984
-0.1	0.01	18.7564742
0	0.01	13.7668283
0.1	0.01	13.7823613
0.2	0.01	13.8036722
0.3	0.01	16.0379387
0.4	0.01	19.5349321
0.5	0.01	19.6648297
0.6	0.01	16.1146323
0.7	0.01	20.8665353
0.8	0.01	20.3808261
0.9	0.01	17.0665999
1	0.01	20.936397
1.2	0.01	22.12253741
1.4	0.01	21.88063448
1.6	0.01	22.36492635
1.8	0.01	24.34147553
2	0.01	23.23502096
-2	0.02	10.52892857
-1.8	0.02	11.50062833
-1.4	0.02	13.0299892
-1.2	0.02	14.61068475
-1	0.02	13.5755155
-0.9	0.02	20.6185801
-0.8	0.02	20.3526949
-0.7	0.02	20.276706
-0.6	0.02	16.2091792
-0.5	0.02	16.8742636
-0.4	0.02	18.4619786
-0.3	0.02	19.6332236
-0.2	0.02	19.0664581
-0.1	0.02	19.9340584
0	0.02	18.7652005
0.1	0.02	15.3061218
0.2	0.02	14.5665431
0.3	0.02	19.1820636
0.4	0.02	18.9537123
0.5	0.02	20.3032923
0.6	0.02	20.9198924
0.7	0.02	22.0995319
0.8	0.02	20.8351704
0.9	0.02	20.3017832
1	0.02	21.8401911
1.2	0.02	22.61757104
1.4	0.02	24.39925472
1.6	0.02	24.82443274
1.8	0.02	25.19445771
2	0.02	24.78584151
-2	0.03	12.98882155
-1.8	0.03	12.85855201
-1.4	0.03	14.20040417
-1.2	0.03	15.52364485
-1	0.03	17.0360159
-0.9	0.03	23.6729376
-0.8	0.03	18.7652005
-0.7	0.03	19.8309932
-0.6	0.03	22.617571
-0.5	0.03	19.01133
-0.4	0.03	19.4752402
-0.3	0.03	20.393682
-0.2	0.03	20.0573608
-0.1	0.03	20.5471576
0	0.03	18.864246
0.1	0.03	19.4483647
0.2	0.03	22.0995319
0.3	0.03	19.4806255
0.4	0.03	20.0351086
0.5	0.03	19.5034281
0.6	0.03	18.5020053
0.7	0.03	22.343068
0.8	0.03	21.7311028
0.9	0.03	22.0906818
1	0.03	21.6527118
1.2	0.03	24.82443274
1.4	0.03	23.87381189
1.6	0.03	24.78584151
1.8	0.03	25.5761629
2	0.03	24.09609527
-2	0.04	13.38431386
-1.8	0.04	13.71116241
-1.4	0.04	16.93527856
-1.2	0.04	18.83076356
-1	0.04	20.0786934
-0.9	0.04	22.0995319
-0.8	0.04	19.3140526
-0.7	0.04	22.0995319
-0.6	0.04	20.905014
-0.5	0.04	24.8260186
-0.4	0.04	20.3685199
-0.3	0.04	20.6294716
-0.2	0.04	24.8244327
-0.1	0.04	20.5368031
0	0.04	20.3017832
0.1	0.04	23.661278
0.2	0.04	20.5262716
0.3	0.04	19.1251933
0.4	0.04	24.5301417
0.5	0.04	23.8209387
0.6	0.04	22.617571
0.7	0.04	20.6368753
0.8	0.04	21.5013802
0.9	0.04	21.8401911
1	0.04	23.7114879
1.2	0.04	23.71378894
1.4	0.04	25.5761629
1.6	0.04	25.27150427
1.8	0.04	27.13184779
2	0.04	27.32284246
-2	0.05	15.2972403
-1.8	0.05	16.59361469
-1.4	0.05	18.62926637
-1.2	0.05	20.38401361
-1	0.05	18.84485
-0.9	0.05	22.343068
-0.8	0.05	21.7311028
-0.7	0.05	20.753391
-0.6	0.05	19.4657472
-0.5	0.05	20.8351704
-0.4	0.05	21.1522146
-0.3	0.05	20.9932221
-0.2	0.05	24.983193
-0.1	0.05	24.9304479
0	0.05	25.1944577
0.1	0.05	24.8318699
0.2	0.05	24.8260186
0.3	0.05	24.8244327
0.4	0.05	23.8209387
0.5	0.05	24.5301417
0.6	0.05	21.7989788
0.7	0.05	20.6187915
0.8	0.05	22.478674
0.9	0.05	19.0541716
1	0.05	26.1056383
1.2	0.05	24.61631956
1.4	0.05	27.13184779
1.6	0.05	27.32284246
1.8	0.05	27.41274124
2	0.05	27.99981691
-2	0.06	18.34933002
-1.8	0.06	18.81222224
-1.4	0.06	18.84485002
-1.2	0.06	18.97555368
-1	0.06	20.6402458
-0.9	0.06	24.7992343
-0.8	0.06	22.0588428
-0.7	0.06	22.0588428
-0.6	0.06	25.161594
-0.5	0.06	24.7933138
-0.4	0.06	24.9303226
-0.3	0.06	22.0796886
-0.2	0.06	26.0725984
-0.1	0.06	24.4670003
0	0.06	22.1313698
0.1	0.06	24.4974507
0.2	0.06	24.5160887
0.3	0.06	25.3971484
0.4	0.06	25.8617766
0.5	0.06	23.7885468
0.6	0.06	21.7342741
0.7	0.06	21.8085732
0.8	0.06	22.4468204
0.9	0.06	21.4697498
1	0.06	27.1318478
1.2	0.06	27.32284246
1.4	0.06	24.64338975
1.6	0.06	28.14496206
1.8	0.06	26.68029293
2	0.06	29.05385447
-2	0.07	19.18554899
-1.8	0.07	20.28653152
-1.4	0.07	19.9379745
-1.2	0.07	22.09376877
-1	0.07	21.8085732
-0.9	0.07	24.5534402
-0.8	0.07	26.0876202
-0.7	0.07	23.3216578
-0.6	0.07	25.3773246
-0.5	0.07	25.8617766
-0.4	0.07	24.5023134
-0.3	0.07	25.9894551
-0.2	0.07	21.1203992
-0.1	0.07	22.3112185
0	0.07	24.7992343
0.1	0.07	25.1766444
0.2	0.07	24.8839366
0.3	0.07	23.7885468
0.4	0.07	23.0257372
0.5	0.07	24.8825339
0.6	0.07	24.6470986
0.7	0.07	27.4525642
0.8	0.07	24.6947397
0.9	0.07	24.8977865
1	0.07	27.4127412
1.2	0.07	27.99981691
1.4	0.07	28.12114503
1.6	0.07	29.05385447
1.8	0.07	28.42600211
2	0.07	29.44408373
-2	0.08	19.42058244
-1.8	0.08	19.40377655
-1.4	0.08	19.97282583
-1.2	0.08	22.41498346
-1	0.08	22.9934107
-0.9	0.08	23.2891889
-0.8	0.08	26.0545936
-0.7	0.08	27.9484789
-0.6	0.08	27.9998169
-0.5	0.08	24.4146918
-0.4	0.08	24.7666158
-0.3	0.08	24.8651423
-0.2	0.08	22.3988339
-0.1	0.08	24.8512014
0	0.08	24.6462726
0.1	0.08	24.3486226
0.2	0.08	27.1318478
0.3	0.08	24.6143461
0.4	0.08	24.4025676
0.5	0.08	24.8496952
0.6	0.08	27.4525642
0.7	0.08	27.3228425
0.8	0.08	27.4127412
0.9	0.08	27.6327448
1	0.08	28.121145
1.2	0.08	27.55915663
1.4	0.08	30.33213527
1.6	0.08	29.21077668
1.8	0.08	30.75166203
2	0.08	30.71440328
-2	0.09	21.99521468
-1.8	0.09	20.71203437
-1.4	0.09	24.45066975
-1.2	0.09	24.96711881
-1	0.09	25.186886
-0.9	0.09	29.9866848
-0.8	0.09	27.8619764
-0.7	0.09	26.8193352
-0.6	0.09	27.1468005
-0.5	0.09	25.3506461
-0.4	0.09	25.8416187
-0.3	0.09	24.7680515
-0.2	0.09	24.7371038
-0.1	0.09	24.8844997
0	0.09	24.7489959
0.1	0.09	24.8474966
0.2	0.09	25.8168303
0.3	0.09	24.6947397
0.4	0.09	24.5733372
0.5	0.09	25.8431684
0.6	0.09	27.9484789
0.7	0.09	27.9998169
0.8	0.09	28.121145
0.9	0.09	26.7805554
1	0.09	30.3321353
1.2	0.09	30.49405537
1.4	0.09	29.52932702
1.6	0.09	30.71440328
1.8	0.09	31.42074451
2	0.09	32.40568749
-2	0.1	20.8845795
-1.8	0.1	22.96110127
-1.4	0.1	25.18688596
-1.2	0.1	23.61680989
-1	0.1	26.0065704
-0.9	0.1	25.8468313
-0.8	0.1	23.208406
-0.7	0.1	28.074869
-0.6	0.1	30.6210378
-0.5	0.1	30.6400396
-0.4	0.1	27.9653085
-0.3	0.1	28.121145
-0.2	0.1	27.9484789
-0.1	0.1	29.9866848
0	0.1	27.5591566
0.1	0.1	30.3321353
0.2	0.1	29.8671608
0.3	0.1	25.7896192
0.4	0.1	27.2676801
0.5	0.1	27.4525642
0.6	0.1	27.0554197
0.7	0.1	27.1318478
0.8	0.1	27.3228425
0.9	0.1	25.3309298
1	0.1	29.7838543
1.2	0.1	29.51702875
1.4	0.1	31.42074451
1.6	0.1	32.40568749
1.8	0.1	34.10531639
2	0.1	34.20073054
-2	0.12	24.21934412
-1.8	0.12	26.52987125
-1.4	0.12	27.13184779
-1.2	0.12	27.337787
-1	0.12	27.34268874
-0.9	0.12	27.63274483
-0.8	0.12	28.14496206
-0.7	0.12	24.88199059
-0.6	0.12	28.26593787
-0.5	0.12	27.96530849
-0.4	0.12	27.55915663
-0.3	0.12	28.44108899
-0.2	0.12	30.33213527
-0.1	0.12	28.15878434
0	0.12	28.65979807
0.1	0.12	30.62103784
0.2	0.12	30.75166203
0.3	0.12	29.38169958
0.4	0.12	29.32684185
0.5	0.12	31.42074451
0.6	0.12	29.95340673
0.7	0.12	31.38913789
0.8	0.12	31.39913185
0.9	0.12	31.55343987
1	0.12	34.10531639
1.2	0.12	33.32811215
1.4	0.12	34.41916913
1.6	0.12	35.66667432
1.8	0.12	36.24514298
2	0.12	36.77667642
-2	0.14	25.71344659
-1.8	0.14	27.26768015
-1.4	0.14	28.12114503
-1.2	0.14	27.53937691
-1	0.14	30.33213527
-0.9	0.14	29.86716075
-0.8	0.14	30.49405537
-0.7	0.14	28.8005139
-0.6	0.14	30.75166203
-0.5	0.14	29.68484221
-0.4	0.14	30.71440328
-0.3	0.14	31.42074451
-0.2	0.14	31.42074451
-0.1	0.14	32.6805438
0	0.14	32.40568749
0.1	0.14	31.3936781
0.2	0.14	33.13082178
0.3	0.14	33.57267159
0.4	0.14	34.20073054
0.5	0.14	34.02668136
0.6	0.14	34.39910947
0.7	0.14	34.71252546
0.8	0.14	34.61422841
0.9	0.14	35.17668402
1	0.14	35.17668402
1.2	0.14	36.77667642
1.4	0.14	37.02622177
1.6	0.14	36.97743563
1.8	0.14	35.51489684
2	0.14	38.02441641
-2	0.16	29.05385447
-1.8	0.16	30.33213527
-1.4	0.16	28.74319748
-1.2	0.16	30.74724673
-1	0.16	31.42074451
-0.9	0.16	32.71361909
-0.8	0.16	32.40568749
-0.7	0.16	31.57329609
-0.6	0.16	34.10531639
-0.5	0.16	34.29610701
-0.4	0.16	34.2667848
-0.3	0.16	34.02668136
-0.2	0.16	35.41190167
-0.1	0.16	35.7304645
0	0.16	34.04233031
0.1	0.16	36.24514298
0.2	0.16	36.24514298
0.3	0.16	35.74140186
0.4	0.16	36.77667642
0.5	0.16	35.64446129
0.6	0.16	35.41162947
0.7	0.16	35.41162947
0.8	0.16	35.2532707
0.9	0.16	37.41717638
1	0.16	37.36795139
1.2	0.16	37.16114584
1.4	0.16	37.6723167
1.6	0.16	39.04377007
1.8	0.16	39.21537615
2	0.16	39.31092037
-2	0.18	30.71440328
-1.8	0.18	31.45366159
-1.4	0.18	34.10531639
-1.2	0.18	34.20073054
-1	0.18	34.41916913
-0.9	0.18	34.6924447
-0.8	0.18	35.66667432
-0.7	0.18	35.1967609
-0.6	0.18	36.24514298
-0.5	0.18	36.79651817
-0.4	0.18	36.77667642
-0.3	0.18	37.12966803
-0.2	0.18	37.02622177
-0.1	0.18	37.15815668
0	0.18	35.2532707
0.1	0.18	37.41717638
0.2	0.18	35.97542461
0.3	0.18	35.81968794
0.4	0.18	36.2334279
0.5	0.18	38.11801666
0.6	0.18	38.48709372
0.7	0.18	36.92868353
0.8	0.18	39.04377007
0.9	0.18	38.38750796
1	0.18	38.75276136
1.2	0.18	39.31092037
1.4	0.18	39.31114993
1.6	0.18	39.64460841
1.8	0.18	39.63788228
2	0.18	38.7904807
-2	0.2	34.27295471
-1.8	0.2	35.41190167
-1.4	0.2	36.24514298
-1.2	0.2	36.77667642
-1	0.2	37.02622177
-0.9	0.2	37.15815668
-0.8	0.2	36.97743563
-0.7	0.2	37.41717638
-0.6	0.2	37.36795139
-0.5	0.2	37.32934116
-0.4	0.2	38.02441641
-0.3	0.2	36.95221892
-0.2	0.2	37.80399917
-0.1	0.2	38.48709372
0	0.2	37.59446386
0.1	0.2	38.76501359
0.2	0.2	38.3258782
0.3	0.2	38.69573527
0.4	0.2	37.95893077
0.5	0.2	38.78946828
0.6	0.2	38.2973204
0.7	0.2	38.8443578
0.8	0.2	39.59256397
0.9	0.2	39.5083968
1	0.2	38.56313655
1.2	0.2	39.73344924
1.4	0.2	40.59497739
1.6	0.2	41.54869724
1.8	0.2	41.54869724
2	0.2	41.79440098
-2	0.22	35.64446129
-1.8	0.22	37.02622177
-1.4	0.22	35.97542461
-1.2	0.22	38.02441641
-1	0.22	38.39579425
-0.9	0.22	38.07736434
-0.8	0.22	38.7652885
-0.7	0.22	38.81730359
-0.6	0.22	39.16299787
-0.5	0.22	38.69573527
-0.4	0.22	38.0423317
-0.3	0.22	39.50717108
-0.2	0.22	39.31114993
-0.1	0.22	38.8443578
0	0.22	39.59256397
0.1	0.22	38.70458056
0.2	0.22	39.63788228
0.3	0.22	39.46696256
0.4	0.22	39.73344924
0.5	0.22	40.59497739
0.6	0.22	39.77714238
0.7	0.22	41.19792504
0.8	0.22	40.85851933
0.9	0.22	41.05426377
1	0.22	41.54869724
1.2	0.22	41.00003971
1.4	0.22	42.33608872
1.6	0.22	44.06289576
1.8	0.22	44.67005178
2	0.22	44.72854004
-2	0.24	38.02441641
-1.8	0.24	38.48709372
-1.4	0.24	38.75276136
-1.2	0.24	38.17322181
-1	0.24	38.2973204
-0.9	0.24	38.81991181
-0.8	0.24	39.59256397
-0.7	0.24	38.72894914
-0.6	0.24	38.56313655
-0.5	0.24	38.51826598
-0.4	0.24	39.73344924
-0.3	0.24	39.79707954
-0.2	0.24	40.59497739
-0.1	0.24	40.39726609
0	0.24	41.54869724
0.1	0.24	41.54869724
0.2	0.24	41.07429643
0.3	0.24	41.02004289
0.4	0.24	41.79440098
0.5	0.24	41.21549934
0.6	0.24	42.01954826
0.7	0.24	42.01954826
0.8	0.24	44.06289576
0.9	0.24	44.67005178
1	0.24	44.67005178
1.2	0.24	44.72854004
1.4	0.24	45.20191453
1.6	0.24	45.20191453
1.8	0.24	45.20191453
2	0.24	45.20191453
-2	0.26	39.25867116
-1.8	0.26	39.31114993
-1.4	0.26	39.63788228
-1.2	0.26	38.81029568
-1	0.26	40.59497739
-0.9	0.26	41.19792504
-0.8	0.26	41.54869724
-0.7	0.26	41.05426377
-0.6	0.26	41.54869724
-0.5	0.26	41.61894771
-0.4	0.26	41.00003971
-0.3	0.26	41.79440098
-0.2	0.26	41.99951971
-0.1	0.26	42.51176454
0	0.26	44.06289576
0.1	0.26	44.67005178
0.2	0.26	44.23764528
0.3	0.26	44.59588729
0.4	0.26	44.18442874
0.5	0.26	44.72854004
0.6	0.26	45.20191453
0.7	0.26	44.7891102
0.8	0.26	45.20191453
0.9	0.26	45.20191453
1	0.26	45.16138882
1.2	0.26	45.20191453
1.4	0.26	45.42476813
1.6	0.26	45.42476813
1.8	0.26	45.42476813
2	0.26	45.42476813
-2	0.28	38.81029568
-1.8	0.28	40.59497739
-1.4	0.28	41.07429643
-1.2	0.28	41.00003971
-1	0.28	42.01954826
-0.9	0.28	42.51176454
-0.8	0.28	44.06289576
-0.7	0.28	44.23764528
-0.6	0.28	44.67005178
-0.5	0.28	44.59588729
-0.4	0.28	44.18442874
-0.3	0.28	44.72854004
-0.2	0.28	44.84989876
-0.1	0.28	45.20191453
0	0.28	45.02762968
0.1	0.28	45.02762968
0.2	0.28	45.20191453
0.3	0.28	45.20191453
0.4	0.28	45.16138882
0.5	0.28	45.13464584
0.6	0.28	45.20951396
0.7	0.28	45.42476813
0.8	0.28	45.42476813
0.9	0.28	45.42476813
1	0.28	45.42476813
1.2	0.28	45.42476813
1.4	0.28	45.92462087
1.6	0.28	45.92462087
1.8	0.28	46.07826844
2	0.28	46.21161378
-2	0.3	41.79440098
-1.8	0.3	42.01954826
-1.4	0.3	44.67005178
-1.2	0.3	44.72854004
-1	0.3	45.20191453
-0.9	0.3	44.7891102
-0.8	0.3	45.20191453
-0.7	0.3	44.98710397
-0.6	0.3	45.20191453
-0.5	0.3	45.20191453
-0.4	0.3	45.16138882
-0.3	0.3	45.09418044
-0.2	0.3	45.42476813
-0.1	0.3	45.42476813
0	0.3	45.42476813
0.1	0.3	45.42476813
0.2	0.3	45.42476813
0.3	0.3	45.42476813
0.4	0.3	45.42476813
0.5	0.3	45.42476813
0.6	0.3	45.92462087
0.7	0.3	45.92462087
0.8	0.3	45.94492311
0.9	0.3	46.07826844
1	0.3	46.07826844
1.2	0.3	46.21161378
1.4	0.3	46.21161378
1.6	0.3	45.92893869
1.8	0.3	46.21161378
2	0.3	46.21161378
-2	0.32	44.18442874
-1.8	0.32	44.84989876
-1.4	0.32	45.20191453
-1.2	0.32	45.20191453
-1	0.32	45.42476813
-0.9	0.32	45.42476813
-0.8	0.32	45.42476813
-0.7	0.32	45.42476813
-0.6	0.32	45.42476813
-0.5	0.32	45.42476813
-0.4	0.32	45.42476813
-0.3	0.32	45.42476813
-0.2	0.32	45.92462087
-0.1	0.32	45.92462087
0	0.32	45.94492311
0.1	0.32	46.07826844
0.2	0.32	46.07826844
0.3	0.32	46.21161378
0.4	0.32	46.21161378
0.5	0.32	46.21161378
0.6	0.32	46.0579662
0.7	0.32	46.0579662
0.8	0.32	46.0579662
0.9	0.32	46.0579662
1	0.32	46.0579662
1.2	0.32	46.21161378
1.4	0.32	46.16477002
1.6	0.32	45.91034631
1.8	0.32	46.22077422
2	0.32	46.17406647
-2	0.34	45.20191453
-1.8	0.34	45.42476813
-1.4	0.34	45.42476813
-1.2	0.34	45.42476813
-1	0.34	45.94492311
-0.9	0.34	45.92462087
-0.8	0.34	45.92462087
-0.7	0.34	45.92462087
-0.6	0.34	45.92462087
-0.5	0.34	46.06228403
-0.4	0.34	46.21161378
-0.3	0.34	46.21161378
-0.2	0.34	46.21161378
-0.1	0.34	46.21161378
0	0.34	46.21161378
0.1	0.34	46.0579662
0.2	0.34	46.21161378
0.3	0.34	46.21161378
0.4	0.34	46.21161378
0.5	0.34	46.21161378
0.6	0.34	46.16477002
0.7	0.34	46.08785443
0.8	0.34	45.91034631
0.9	0.34	46.08785443
1	0.34	46.22077422
1.2	0.34	46.17406647
1.4	0.34	45.99681548
1.6	0.34	45.99681548
1.8	0.34	45.99681548
2	0.34	46.78885931
-2	0.36	45.42476813
-1.8	0.36	45.94492311
-1.4	0.36	46.07826844
-1.2	0.36	46.21161378
-1	0.36	46.0579662
-0.9	0.36	46.21161378
-0.8	0.36	46.21161378
-0.7	0.36	46.0579662
-0.6	0.36	46.0579662
-0.5	0.36	46.0579662
-0.4	0.36	46.21161378
-0.3	0.36	46.21161378
-0.2	0.36	46.16477002
-0.1	0.36	46.08785443
0	0.36	46.08785443
0.1	0.36	46.08785443
0.2	0.36	46.22077422
0.3	0.36	45.86408823
0.4	0.36	46.17406647
0.5	0.36	46.17406647
0.6	0.36	46.17406647
0.7	0.36	46.17406647
0.8	0.36	45.99681548
0.9	0.36	46.17406647
1	0.36	46.02542889
1.2	0.36	46.61106363
1.4	0.36	46.61106363
1.6	0.36	46.61106363
1.8	0.36	46.78885931
2	0.36	46.78885931
-2	0.38	45.92462087
-1.8	0.38	46.21161378
-1.4	0.38	46.21161378
-1.2	0.38	46.19131153
-1	0.38	46.14449728
-0.9	0.38	45.89010886
-0.8	0.38	46.08785443
-0.7	0.38	45.92682929
-0.6	0.38	46.22077422
-0.5	0.38	45.86408823
-0.4	0.38	45.86408823
-0.3	0.38	46.17406647
-0.2	0.38	45.99681548
-0.1	0.38	46.06179843
0	0.38	45.99681548
0.1	0.38	46.17406647
0.2	0.38	46.02542889
0.3	0.38	45.99681548
0.4	0.38	46.62325529
0.5	0.38	46.63976498
0.6	0.38	46.61106363
0.7	0.38	46.78885931
0.8	0.38	46.78885931
0.9	0.38	46.61106363
1	0.38	46.62325529
1.2	0.38	46.78885931
1.4	0.38	46.78885931
1.6	0.38	46.78885931
1.8	0.38	46.78885931
2	0.38	46.78885931
-2	0.4	46.0952725
-1.8	0.4	46.14449728
-1.4	0.4	45.91034631
-1.2	0.4	46.17406647
-1	0.4	46.17406647
-0.9	0.4	46.17406647
-0.8	0.4	45.99681548
-0.7	0.4	46.17406647
-0.6	0.4	46.17406647
-0.5	0.4	45.97660734
-0.4	0.4	46.61106363
-0.3	0.4	46.78885931
-0.2	0.4	46.61106363
-0.1	0.4	46.61106363
0	0.4	46.78885931
0.1	0.4	46.62325529
0.2	0.4	46.62325529
0.3	0.4	46.78885931
0.4	0.4	46.78885931
0.5	0.4	46.78885931
0.6	0.4	46.78885931
0.7	0.4	46.78885931
0.8	0.4	46.78885931
0.9	0.4	46.78885931
1	0.4	46.78885931
1.2	0.4	46.78885931
1.4	0.4	46.78885931
1.6	0.4	46.78885931
1.8	0.4	46.75882514
2	0.4	46.75882514
}\tblwithout
\begin{figure}[h!]
\centering
\scriptsize
        \begin{tikzpicture}[scale=.75]
    \begin{axis}[
        colormap/RdYlBu,
        xtick={-1,0,+1}, 
        xticklabels={
            $-1$,
            $\beta^{time}$,
            $+1$},
            ytick={-.1,0,.1}, 
        yticklabels={
                $-0.1$,
                $\beta^{\mathrm{price}}$,
                $+0.1$},
        label style={font=\footnotesize},
        tick label style={font=\footnotesize}  
    ]
    \addplot3 
    [surf, fill= white, mesh/rows=41, mesh/cols=30]
    table [x index=0, y index=1, z index=2] \tblwithout;
    \end{axis}
    \end{tikzpicture}
\caption{Comparison of the objective values from sALNS and RFTS heuristic.}
\label{fig:improvement_obj_ALNS_RFTS}
\end{figure}

\subsubsection{Benefits of personalized discounts over uniform pricing} 
\label{secFairnessCheck} 
In this final section, we analyze the effects of personalized discounts compared to uniform pricing, as captured by the introduction of Constraints~\ref{eq:discrimination}, described in Section~\ref{problemdef1}. To accommodate this additional constraint, the sALNS algorithm was slightly adjusted. This required modifications to the repair operators, local search operators, and the RFTS heuristic. Specifically, the repair and local search operators were updated to preserve existing price levels for already allocated slots, while newly assigned slots followed structured discounting rules. Additionally, the RFTS heuristic was modified to include a final step that enforces uniform pricing by applying the maximum discount rate to each slot, ensuring that all customers remain in the system.

The results of comparing personalized discounts to uniform pricing are shown in Table~\ref{tab:with-without}, which presents both the time gap (the computational time required for each approach) and the objective gap (the percentage of profit lost when using uniform pricing compared to personalized discounts). These results are based on the first dataset (\textit{i.e.} behavioral setting). As shown in the table, the time gap for the personalized discount strategy is negative, indicating that computing personalized discounts requires more time than the uniform pricing approach due to the added complexity of the optimization process. This increase in computational complexity is expected, as the larger solution search space associated with differentiated pricing across customers leads to longer computation times. Additionally, the time gap increases with the number of customers, reflecting the growing complexity of the search space as the customer base expands.

Despite the additional computational cost, the objective gap highlights the profitability gains that come from adopting personalized discounts. The personalized discount strategy consistently delivers higher profit compared to uniform pricing. While uniform pricing is often praised for its simplicity and perceived fairness, it comes with notable limitations in terms of profit optimization and operational efficiency. The core issue with uniform pricing is its inability to fully capitalize on the diversity of customer preferences and price sensitivities. Personalized discount strategies, on the other hand, recognize and leverage this diversity. By adjusting discounts based on individual preferences, as well as operational costs, personalized pricing better aligns with customer needs and enhances resource utilization. In doing so, personalized strategies capture more value and improve operational efficiency.

\begin{table}[h]
\scriptsize
\caption{Comparison of computational time and profit gap between uniform pricing and personalized discount strategy across instances of different sizes and 100 scenarios.}
\begin{tabular}{cccccc}
\toprule
 \multirow{2}{*}{Customers} & \multicolumn{2}{c}{Objective Gap (\%)}  & & \multicolumn{2}{c}{Time Gap (\%)} \\ \cline{2-3} \cline{5-6}
 & Average          & SD         & & Average          & SD           \\ \midrule
5               & 43.93            & 0.36          && -16.12           & 45.86          \\
10              & 54.73            & 0.84          && -114.76          & 50.68          \\
20              & 49.41            & 0.83          && -202.64          & 87.69          \\
30              & 46.37            & 1.69          && -242.28          & 64.19          \\
40              & 48.93            & 1.21          && -277.81          & 84.43          \\
50              & 44.03            & 0.96          && -296.06          & 50.79          \\
60              & 42.49            & 0.59          && -456.84          & 84.32          \\
80              & 41.19            & 0.93          && -370.57          & 120.36         \\
100             & 39.55            & 0.90          & & -603.09          & 217.83         \\ \bottomrule
\end{tabular}
\label{tab:with-without}
\end{table}

\section{Conclusion}\label{sec6}

In this paper, we addressed the challenges of TTSM in AHDs within SBM, focusing on capturing customer heterogeneity in delivery preferences. To this end, we developed an optimization framework based on the ML model, which captures individual-level heterogeneity through random coefficients, enabling a richer and more realistic representation of customer choice behavior compared to traditional closed-form models.

Integrating ML into slot assortment optimization is nontrivial because choice probabilities lack closed-form expressions, resulting in nonlinear and simulation-based formulations. To overcome this, we employed a scenario-based reformulation leveraging SAA. This approach eliminates the need for closed-form choice probabilities and allows the model to accommodate any simulation-based RUM into slot assortment optimization.

The resulting stochastic TTSM-ML was formulated as a MILP, which simultaneously optimizes delivery slot offerings, customer time slot choices driven by heterogeneous preferences, and routing decisions. Given the computational complexity of large-scale instances, we proposed a tailored sALNS that exploits the MILP’s scenario-based structure and employs multiple neighborhood operators to efficiently explore the solution space.

Through extensive numerical experiments on a diverse set of instances varying in customer preferences, spatial configuration, and customer size, we validated the effectiveness of our framework. Results demonstrate that the sALNS heuristic efficiently handles large-scale tactical slot optimization problems and that the ML-based model outperforms traditional approaches based on the simpler MNL model. This improved modeling of customer heterogeneity allows e-retailers to better tailor their service offerings, leading to enhanced profitability and customer retention.

Our research contributes to logistics and operations management by providing a general and flexible framework that integrates advanced choice modeling with tactical decision-making in time slot management. Its ability to support any simulation-based RUM makes it a valuable tool for practitioners requiring customized choice models across diverse contexts. Beyond e-retail, this framework can be adapted to other domains where supply-demand interactions are interdependent, such as healthcare scheduling, or field service operations.

While this study focuses on tactical decisions, an important direction for future work is extending the framework to dynamic and operational contexts. Incorporating real-time demand updates, dynamic fleet routing, and adaptive pricing would increase its practical relevance. Additionally, integrating online decision-making techniques that combine tactical assortment planning with operational dispatching could bridge the gap between strategic objectives and real-time execution, fostering more resilient and efficient attended delivery systems.

Finally, access to real-world data through industry partnerships would be invaluable for refining utility specifications and validating the framework in practice. Collaborations with e-retailers and logistics providers could help quantify the extent to which ML models offer tangible advantages over simpler alternatives, such as finite-mixture MNL, which typically require richer data for accurate group identification, in high-stakes delivery environments.

\bibliographystyle{elsarticle-harv} 
\bibliography{references}
\newpage
\setcounter{section}{0}
\setcounter{page}{1}
\setcounter{figure}{0}
\setcounter{equation}{0}
\renewcommand{\thesection}{A\arabic{section}}
\renewcommand{\thepage}{a\arabic{page}}
\renewcommand{\thetable}{A\arabic{table}}
\renewcommand{\thefigure}{A\arabic{figure}}

\setcounter{ead}{0}
\setcounter{tnote}{0}
\setcounter{fnote}{0}
\setcounter{cnote}{0}
\setcounter{author}{0}
\setcounter{affn}{0}
\resetTitleCounters

\section*{Appendices}

\section{Synthetic dataset configurations} \label{app:beta_inputs}

\begin{table}[h]
\centering
\caption{Input parameters used for generating the five synthetic datasets.}
\label{tab:beta_inputs_app}
\scriptsize
\begin{tabular}{ccccccc}
\toprule
Dataset & $\mu^{time}_{t1}$ & $\mu^{time}_{t2}$ & $\mu^{time}_{t3}$ & $\sigma(\beta^{time})$ & $\mathbb{E}[\beta^{\mathrm{price}}]$ & $\sigma(\beta^{\mathrm{price}})$ \\ \hline
1 & 8.0  & 10.0 & 6.8  & 1& -0.15 & 0.25 \\
2 & 8.0  & 6.8  & 10.0 & 1& -0.05 & 0.25 \\
3 & 8.0  & 6.8  & 10.0 & 1& -0.10 & 0.25 \\
4 & 11.0 & 7.8  & 15.0 & 1& -0.15 & 0.25 \\
5 & 6.0  & 5.8  & 7.5  & 1& -0.15 & 0.25 \\ \bottomrule
\end{tabular}
\end{table}

\section{Full estimation results for all datasets}\label{app:ml_mnl_tables}

\begin{landscape}
\begin{table}[]
\centering
\caption{Comparison of parameter estimates and log-likelihood values for MNL and ML models across five synthetic datasets.}
\label{tab:beta_mnl_ml_all1}
\resizebox{\textwidth}{!}{%
\begin{tabular}{cccccccccccccccc}
\toprule
Dataset &
   &
  \multicolumn{4}{c}{1} &
   &
  \multicolumn{4}{c}{2} &
   &
  \multicolumn{4}{c}{3} \\ \cline{1-1} \cline{3-6} \cline{8-11} \cline{13-16} 
Choice Model &
   &
  \multicolumn{2}{c}{MNL} &
  \multicolumn{2}{c}{ML} &
   &
  \multicolumn{2}{c}{MNL} &
  \multicolumn{2}{c}{ML} &
   &
  \multicolumn{2}{c}{MNL} &
  \multicolumn{2}{c}{ML} \\ \cline{1-1} \cline{3-6} \cline{8-11} \cline{13-16} 
Parameters &
   &
  Value &
  P-value &
  Value &
  P-value &
   &
  Value &
  P-value &
  Value &
  P-value &
   &
  Value &
  P-value &
  Value &
  P-value \\\cline{1-1} \cline{3-6} \cline{8-11} \cline{13-16} 
$\beta^{\mathrm{time}}_{t1}$ &
   &
  1.0690 &
  0.0000 &
  5.8460 &
  0.0000 &
   &
  0.5449 &
  0.0000 &
  6.5141 &
  0.0000 &
   &
  1.1404 &
  0.0000 &
  6.2276 &
  0.0000 \\
$\beta^{\mathrm{time}}_{t2}$ &
   &
  2.0618 &
  0.0000 &
  7.4001 &
  0.0000 &
   &
  2.1961 &
  0.0000 &
  5.5942 &
  0.0000 &
   &
  0.5543 &
  0.0000 &
  5.2542 &
  0.0000 \\
$\beta^{\mathrm{time}}_{t3}$ &
   &
  0.5236 &
  0.0000 &
  4.9178 &
  0.0000 &
   &
  1.1151 &
  0.0000 &
  8.0575 &
  0.0000 &
   &
  2.1434 &
  0.0000 &
  7.7375 &
  0.0000 \\
$\mathbb{E}(\beta^{\mathrm{price}})$ &
   &
  -0.0257 &
  0.0000 &
  -0.0982 &
  0.0000 &
   &
  -0.0079 &
  0.0018 &
  -0.0329 &
  0.0000 &
   &
  -0.0182 &
  0.0000 &
  -0.0665 &
  0.0000 \\
$\sigma(\beta^{\mathrm{price}})$ &
   &
   &
   &
  0.1772 &
  0.0000 &
   &
   &
   &
  0.1910 &
  0.0000 &
   &
   &
   &
  0.1885 &
  0.0000 \\ \cline{1-1} \cline{3-6} \cline{8-11} \cline{13-16} 
Log-Likelihood (LL) &
   &
  \multicolumn{2}{c}{-43151.2603} &
  \multicolumn{2}{c}{-41437.6091} &
   &
  \multicolumn{2}{c}{-38605.6077} &
  \multicolumn{2}{c}{-36959.0563} &
   &
  \multicolumn{2}{c}{-41395.9314} &
  \multicolumn{2}{c}{-39715.9921} \\ \cline{1-1} \cline{3-6} \cline{8-11} \cline{13-16} 
\multirow{2}{*}{LL-Ratio Test} &
   &
  \multicolumn{2}{c}{Statistics} &
  \multicolumn{2}{c}{P-value} &
   &
  \multicolumn{2}{c}{Statistics} &
  \multicolumn{2}{c}{P-value} &
   &
  \multicolumn{2}{c}{Statistics} &
  \multicolumn{2}{c}{P-value} \\ \cline{3-6} \cline{8-11} \cline{13-16} 
 &
   &
  \multicolumn{2}{c}{3427.3024} &
  \multicolumn{2}{c}{0.0000} &
   &
  \multicolumn{2}{c}{3293.1026} &
  \multicolumn{2}{c}{0.0000} &
   &
  \multicolumn{2}{c}{3359.8785} &
  \multicolumn{2}{c}{0.0000} \\ \bottomrule
\end{tabular}%
}
\end{table}

\begin{table}[]
\centering
\tiny
\caption{Comparison of parameter estimates and log-likelihood values for MNL and ML models across five synthetic datasets.}
\label{tab:beta_mnl_ml_all2}
\resizebox{\textwidth}{!}{%
\begin{tabular}{clcccclcccc}
\toprule
Dataset            &  & \multicolumn{4}{c}{4}                 &  & \multicolumn{4}{c}{5}                 \\ \cline{1-1} \cline{3-6} \cline{8-11} 
Choice Model &
   &
  \multicolumn{2}{c}{MNL} &
  \multicolumn{2}{c}{ML} &
   &
  \multicolumn{2}{c}{MNL} &
  \multicolumn{2}{c}{ML} \\ \cline{1-1} \cline{3-6} \cline{8-11} 
Parameters         &  & Value   & P-value & Value   & P-value &  & Value   & P-value & Value   & P-value \\ \cline{1-1} \cline{3-6} \cline{8-11} 
$\beta^{\mathrm{time}}_{t1}$      &  & 2.4289  & 0.0000  & 9.0465  & 0.0000  &  & 0.4492  & 0.0000  & 4.2659  & 0.0000  \\
$\beta^{\mathrm{time}}_{t2}$      &  & 1.3039  & 0.0000  & 6.4592  & 0.0000  &  & 0.3500  & 0.0001  & 4.1106  & 0.0000  \\
$\beta^{\mathrm{time}}_{t3}$      &  & 4.1649  & 0.0000  & 12.2099 & 0.0000  &  & 1.1932  & 0.0000  & 5.3833  & 0.0000  \\
$\mathbb{E}(\beta^{\mathrm{price}})$     &  & -0.0447 & 0.0000  & -0.1105 & 0.0000  &  & -0.0200 & 0.0000  & -0.0985 & 0.0000  \\
$\sigma(\beta^{\mathrm{price}})$ &  &         &         & 0.1976  & 0.0000  &  &         &         & 0.1609  & 0.0000  \\ \cline{1-1} \cline{3-6} \cline{8-11} 
LL &
   &
  \multicolumn{2}{c}{-31535.3630} &
  \multicolumn{2}{c}{-28333.1699} &
   &
  \multicolumn{2}{c}{-46599.8333} &
  \multicolumn{2}{c}{-45612.2656} \\ \cline{1-1} \cline{3-6} \cline{8-11} 
\multirow{2}{*}{LL-Ratio Test} &
   &
  \multicolumn{2}{c}{Statistics} &
  \multicolumn{2}{c}{P-value} &
   &
  \multicolumn{2}{c}{Statistics} &
  \multicolumn{2}{c}{P-value} \\ \cline{3-6} \cline{8-11}
 &
   &
  \multicolumn{2}{c}{6404.3862} &
  \multicolumn{2}{c}{0.0000} &
   &
  \multicolumn{2}{c}{1975.1355} &
  \multicolumn{2}{c}{0.0000} \\ \bottomrule
\end{tabular}%
}
\end{table}
\end{landscape}

\section{Number of scenarios for 10 customers} \label{app:Scenario10}
To verify the consistency of scenario size requirements in larger instances, we replicate the in-sample and out-of-sample analyses described in Section~\ref{secSupSize} for 10-customer instances. Due to the increased computational complexity, exact optimization is no longer tractable. We therefore rely on the sALNS method to solve the TTSM-ML model under varying numbers of scenarios ($R$).

Figures~\ref{fig:InSample10} and~\ref{fig:OutSample10} report the results. Both figures show trends similar to those observed in smaller instances: deviations decrease and stabilize around 2\% when using 100 scenarios, and the standard deviations remain consistently low beyond this threshold. These results confirm that $R = 100$ is a sufficient scenario size for larger problem instances.

\begin{figure}[!htbp]
\centering
\scriptsize
\begin{subfigure}{0.32\textwidth}
\captionsetup{justification=centering}
\footnotesize
    \tikz{
    \begin{axis}[
  width=\linewidth,
      ytick={0,5,...,35}, ytick align=outside, ytick style={draw=none},
      xtick={0,20,...,100,120}, xtick align=outside, xtick style={draw=none},
      xlabel=Scenarios,
      ylabel={In-sample deviations (\%)},
      height=\axisdefaultheight,
      grid=both,
  legend to name=sharedlegend
      ]
    \addplot+[
      BlueViolet, 
      mark=triangle,
      mark options={BlueViolet, scale=0.8},
      smooth, 
      error bars/.cd, 
        y fixed,
        y dir=both, 
        y explicit
    ] 
    table [x=x, y=y,y error=error] {
        x  y       error
        5	15.90	10.37
        10	8.07	7.07
        15	5.02	5.31
        20	5.24	3.77
        30	4.42	2.78
        40	3.34	2.77
        50	2.72	1.55
        60	2.42	2.02
        80	2.09	1.22
        100	1.13	0.88
        110	0.83	0.61   

    };
    \addlegendentry{Random map}
    \addplot+[
      purple, 
      mark=square,
      mark options={purple, scale=0.8},
      smooth, 
      error bars/.cd, 
        y fixed,
        y dir=both, 
        y explicit
    ] 
    table [x=x, y=y,y error=error] {
        x  y       error
        5	10.45	5.89
        10	5.53	4.22
        15	4.77	3.30
        20	4.53	3.85
        30	3.79	2.89
        40	3.04	1.72
        50	3.22	1.69
        60	2.26	1.49
        80	2.06	1.31
        100	1.01	0.79
        110	0.90	0.60
    };
    \addlegendentry{Clustered map}
    \addplot+[
      teal, 
      mark=*,
      mark options={teal, scale=0.8},
      smooth, 
      error bars/.cd, 
        y fixed,
        y dir=both, 
        y explicit
    ] 
    table [x=x, y=y,y error=error] {
        x  y       error
        5	13.37	8.47
        10	6.06	6.04
        15	5.12	3.78
        20	5.24	3.96
        30	4.36	3.20
        40	3.34	2.87
        50	2.82	2.41
        60	2.01	1.53
        80	1.71	1.16
        100	1.26	0.56
        110	0.74	0.60
    };
    \addlegendentry{Mixed map}
    \end{axis}
    }
\caption{Dataset 1}
\label{fig:Out-of-sample-Test-10}
\end{subfigure}
\hfill
\begin{subfigure}{0.32\textwidth}
\captionsetup{justification=centering}
\footnotesize
\tikz{
\begin{axis}[
  width=\linewidth,
  ytick={0,5,...,50}, ytick align=outside, ytick style={draw=none},
  xtick={0,20,...,120}, xtick align=outside, xtick style={draw=none},
  xlabel=Scenarios,
  ylabel={In-sample deviations (\%)},
  height=\axisdefaultheight,
  grid=both
  ]
    \addplot+[
      BlueViolet, 
      mark=triangle,
      mark options={BlueViolet, scale=0.8},
      smooth,
      error bars/.cd, 
        y fixed,
        y dir=both, 
        y explicit
    ] 
    table [x=x, y=y,y error=error] {
        x  y       error
        5	26.69	17.37
        10	16.55	11.88
        15	12.54	9.19
        20	10.04	8.58
        30	9.62	6.01
        40	7.73	4.59
        50	6.73	3.87
        60	3.98	2.75
        80	3.19	2.83
        100	2.10	1.11
        110	1.02	0.95
    };
    \addplot+[
      purple, 
      mark=square,
      mark options={purple, scale=0.8},
      smooth, 
      error bars/.cd, 
        y fixed,
        y dir=both, 
        y explicit
    ] 
    table [x=x, y=y,y error=error] {
        x  y       error
        5	18.09	15.59
        10	12.91	10.11
        15	9.40	6.01
        20	7.46	5.48
        30	5.50	5.77
        40	5.20	4.27
        50	3.51	3.08
        60	3.98	2.81
        80	2.80	1.72
        100	2.40	1.53
        110	1.40	1.03
    };
    \addplot+[
      teal, 
      mark = *,
      mark options={teal, scale=0.8},
      smooth, 
      error bars/.cd, 
        y fixed,
        y dir=both, 
        y explicit
    ] 
    table [x=x, y=y,y error=error] {
        x  y       error
        5	27.99	15.68
        10	14.88	10.01
        15	9.42	9.33
        20	9.88	7.54
        30	8.30	5.97
        40	6.39	4.25
        50	4.24	3.27
        60	3.28	3.09
        80	2.73	2.10
        100	1.53	1.06
        110	1.25	0.90
    };
    \end{axis}
    }
\caption{Dataset 2}
\label{fig:Out-of-sample-Test-20}
\end{subfigure}
\hfill
\begin{subfigure}{0.32\textwidth}
\captionsetup{justification=centering}
\footnotesize
\tikz{
\begin{axis}[
  width=\linewidth,
  ytick={0,5,...,45}, ytick align=outside, ytick style={draw=none},
  xtick={0,20,...,120}, xtick align=outside, xtick style={draw=none},
  xlabel=Scenarios,
  ylabel={In-sample deviations (\%)},
  height=\axisdefaultheight,
  grid=both
  ]
    \addplot+[
      BlueViolet, 
      mark=triangle,
      mark options={BlueViolet, scale=0.8},
      smooth,
      error bars/.cd, 
        y fixed,
        y dir=both, 
        y explicit
    ] 
    table [x=x, y=y,y error=error] {
        x  y       error
        5	12.44	7.22
        10	6.48	5.18
        15	6.46	4.61
        20	5.56	4.41
        30	4.46	3.32
        40	4.13	2.24
        50	3.51	2.59
        60	3.10	2.11
        80	2.55	1.42
        100	1.24	0.79
        110	0.87	0.66
    };
    \addplot+[
      purple, 
      mark=square,
      mark options={purple, scale=0.8},
      smooth, 
      error bars/.cd, 
        y fixed,
        y dir=both, 
        y explicit
    ] 
    table [x=x, y=y,y error=error] {
        x  y       error
        5	16.72	11.89
        10	7.83	7.75
        15	5.61	5.54
        20	6.14	5.59
        30	5.30	4.53
        40	4.51	4.29
        50	3.60	3.32
        60	3.39	2.06
        80	2.40	1.78
        100	1.95	1.93
        110	0.89	0.66
    };
    \addplot+[
      teal, 
      mark = *,
      mark options={teal, scale=0.8},
      smooth, 
      error bars/.cd, 
        y fixed,
        y dir=both, 
        y explicit
    ] 
    table [x=x, y=y,y error=error] {
        x  y       error
        5	23.76	14.95
        10	13.68	11.56
        15	10.39	7.42
        20	9.29	7.56
        30	8.05	6.12
        40	6.91	5.05
        50	5.98	3.89
        60	3.88	2.70
        80	2.74	2.45
        100	1.54	1.15
        110	1.25	1.01
    };
    \end{axis}
    }
\caption{Dataset 3}
\label{fig:Out-of-sample-Test-30}
\end{subfigure}
\hfill
\begin{subfigure}{0.32\textwidth}
\captionsetup{justification=centering}
\footnotesize
\tikz{
\begin{axis}[
  width=\linewidth,
  ytick={0,5,...,40}, ytick align=outside, ytick style={draw=none},
  xtick={0,20,...,120}, xtick align=outside, xtick style={draw=none},
  xlabel=Scenarios,
  ylabel={In-sample deviations (\%)},
  height=\axisdefaultheight,
  grid=both
  ]
    \addplot+[
      BlueViolet, 
      mark=triangle,
      mark options={BlueViolet, scale=0.8},
      smooth,
      error bars/.cd, 
        y fixed,
        y dir=both, 
        y explicit
    ] 
    table [x=x, y=y,y error=error] {
        x  y       error
        5	15.11	13.23
        10	12.58	9.13
        15	9.77	7.02
        20	6.74	5.72
        30	5.63	5.67
        40	5.16	4.17
        50	4.45	3.23
        60	3.80	3.00
        80	2.47	2.00
        100	1.96	1.80
        110	1.23	0.97
    };
    \addplot+[
      purple, 
      mark=square,
      mark options={purple, scale=0.8},
      smooth, 
      error bars/.cd, 
        y fixed,
        y dir=both, 
        y explicit
    ] 
    table [x=x, y=y,y error=error] {
        x  y       error
        5	21.63	14.20
        10	11.78	10.46
        15	7.74	6.94
        20	7.01	4.82
        30	7.62	3.95
        40	5.28	2.67
        50	4.13	2.91
        60	3.10	2.78
        80	2.24	1.81
        100	1.22	0.82
        110	1.17	0.93
    };
    \addplot+[
      teal, 
      mark = *,
      mark options={teal, scale=0.8},
      smooth, 
      error bars/.cd, 
        y fixed,
        y dir=both, 
        y explicit
    ] 
    table [x=x, y=y,y error=error] {
        x  y       error
        5	12.17	7.89
        10	7.24	4.75
        15	6.87	4.55
        20	5.39	4.26
        30	4.04	3.33
        40	3.79	2.27
        50	3.53	2.57
        60	2.65	2.37
        80	2.12	1.47
        100	1.20	0.73
        110	0.77	0.58
    };
    \end{axis}
    }
\caption{Dataset 4}
\label{fig:Out-of-sample-Test-40}
\end{subfigure}
\hfill
\begin{subfigure}{0.32\textwidth}
\captionsetup{justification=centering}
\footnotesize
\tikz{
\begin{axis}[
  width=\linewidth,
  ytick={0,5,...,35}, ytick align=outside, ytick style={draw=none},
  xtick={0,20,...,120}, xtick align=outside, xtick style={draw=none},
  xlabel=Scenarios,
  ylabel={In-sample deviations (\%)},
  height=\axisdefaultheight,
  grid=both
  ]
    \addplot+[
      BlueViolet, 
      mark=triangle,
      mark options={BlueViolet, scale=0.8},
      smooth,
      error bars/.cd, 
        y fixed,
        y dir=both, 
        y explicit
    ] 
    table [x=x, y=y,y error=error] {
        x  y       error
        5	15.31	9.55
        10	8.35	5.47
        15	6.74	4.82
        20	6.35	5.06
        30	4.89	4.44
        40	3.97	3.71
        50	3.32	2.97
        60	2.92	1.88
        80	2.42	1.76
        100	1.20	1.33
        110	1.04	0.52
    };
    \addplot+[
      purple, 
      mark=square,
      mark options={purple, scale=0.8},
      smooth, 
      error bars/.cd, 
        y fixed,
        y dir=both, 
        y explicit
    ] 
    table [x=x, y=y,y error=error] {
        x  y       error
        5	18.45	11.27
        10	11.74	7.69
        15	7.39	6.40
        20	6.17	5.74
        30	6.46	4.04
        40	5.04	3.45
        50	4.13	2.61
        60	3.03	1.78
        80	2.02	1.72
        100	1.54	0.94
        110	0.74	0.67
    };
    \addplot+[
      teal, 
      mark = *,
      mark options={teal, scale=0.8},
      smooth, 
      error bars/.cd, 
        y fixed,
        y dir=both, 
        y explicit
    ] 
    table [x=x, y=y,y error=error] {
        x  y       error
        5	12.63	10.33
        10	8.86	7.08
        15	7.39	5.46
        20	5.67	5.36
        30	4.39	4.90
        40	4.42	3.19
        50	3.30	2.70
        60	3.12	2.48
        80	2.17	1.48
        100	1.65	0.89
        110	0.99	0.67
    };
    \end{axis}
    }
\caption{Dataset 5}
\label{fig:Out-of-sample-Test-50}
\end{subfigure}
\hfill
\begin{subfigure}{0.32\textwidth}
\centering
\begin{tikzpicture}
\node at (0,0) {\ref{sharedlegend}}; 
\end{tikzpicture}
\end{subfigure}
\caption{In-sample test: It evaluates the deviation in objective values across problems with varying scenario sizes, ranging from 5 to 110 scenarios, and 10 customers, compared to the problem with a fixed number of 120 scenarios.}
\label{fig:InSample10}
\end{figure}

\begin{figure}[!htbp]
\centering
\scriptsize
\begin{subfigure}{0.32\textwidth}
\captionsetup{justification=centering}
\footnotesize
    \tikz{
    \begin{axis}[
  width=\linewidth,
      ytick={0,2,4.00,...,28}, ytick align=outside, ytick style={draw=none},
      xtick={0,20,...,100,120}, xtick align=outside, xtick style={draw=none},
      xlabel=Scenarios,
      ylabel={Out-of-sample deviations (\%)},
      height=\axisdefaultheight,
      grid=both,
  legend to name=sharedlegend
      ]
    \addplot+[
      BlueViolet, 
      mark=triangle,
      mark options={BlueViolet, scale=0.8},
      smooth,
      error bars/.cd, 
        y fixed,
        y dir=both, 
        y explicit
    ] 
    table [x=x, y=y,y error=error] {
        x  y       error
        5	0.5417265	0.4224068
        10	0.3010654	0.2391505
        15	0.2167887	0.1811161
        20	0.2535585	0.17064
        30	0.1366105	0.1063408
        40	0.121379	0.0892029
        50	0.0803721	0.0639201
        60	0.0686294	0.0540531
        80	0.0423648	0.0398506
        100	0.0714574	0.059531
        110	0.0447575	0.0452627
        120	0.0446406	0.0541608
    };
    \addlegendentry{Random map}
    \addplot+[
      purple, 
      mark=square,
      mark options={purple, scale=0.8},
      smooth, 
      error bars/.cd, 
        y fixed,
        y dir=both, 
        y explicit
    ] 
    table [x=x, y=y,y error=error] {
        x  y       error
       5	6.067924062	4.243256328
        10	5.540682305	2.889229023
        15	3.427708126	2.172047964
        20	2.421275545	1.398563594
        30	1.985099139	1.322958029
        40	1.498467805	1.146891171
        50	1.448198181	0.907004649
        60	1.580900877	0.920386929
        80	1.225168965	0.698486432
        100	1.052965248	0.720228305
        110	0.806330795	0.503911677
        120	0.814885818	0.530813554
    };
    \addlegendentry{Clustered map}
    \addplot+[
      teal, 
      mark = *,
      mark options={teal, scale=0.8},
      smooth, 
      error bars/.cd, 
        y fixed,
        y dir=both, 
        y explicit
    ] 
    table [x=x, y=y,y error=error] {
        x  y       error
        5	4.186301508	3.416461228
        10	4.367998591	2.603629376
        15	3.197181614	1.956824877
        20	2.4767585	1.53920021
        30	2.244547732	1.318236871
        40	1.743591363	1.056856218
        50	1.577025349	1.180819495
        60	1.229930288	0.730441647
        80	1.221970883	0.760904112
        100	1.081998349	0.857802642
        110	1.46241008	0.849435435
        120	1.234080354	0.699857172
    };
    \addlegendentry{Mixed map}
    \end{axis}
    }
\caption{Dataset 1}
\label{fig:Out-of-sample-Test-10}
\end{subfigure}
\hfill
\begin{subfigure}{0.32\textwidth}
\captionsetup{justification=centering}
\footnotesize
\tikz{
\begin{axis}[
  width=\linewidth,
  ytick={0,5,...,50}, ytick align=outside, ytick style={draw=none},
  xtick={0,20,...,120}, xtick align=outside, xtick style={draw=none},
  xlabel=Scenarios,
  ylabel={Out-of-sample deviations (\%)},
  height=\axisdefaultheight,
  grid=both
  ]
    \addplot+[
      BlueViolet, 
      mark=triangle,
      mark options={BlueViolet, scale=0.8},
      smooth,
      error bars/.cd, 
        y fixed,
        y dir=both, 
        y explicit
    ] 
    table [x=x, y=y,y error=error] {
        x  y       error
        5	9.114420759	5.183449092
        10	2.262996122	1.696703307
        15	3.518212808	2.482370685
        20	1.277747152	1.000607006
        30	1.570933022	0.962052665
        40	0.831363877	0.675805317
        50	1.484347595	0.822934014
        60	1.233274836	0.542896197
        80	0.57680221	0.412391108
        100	0.417600448	0.301720095
        110	0.394037997	0.280611527
        120	0.946310858	0.569297824
    };
    \addplot+[
      purple, 
      mark=square,
      mark options={purple, scale=0.8},
      smooth, 
      error bars/.cd, 
        y fixed,
        y dir=both, 
        y explicit
    ] 
    table [x=x, y=y,y error=error] {
        x  y       error
        5	4.398329612	2.998012386
        10	6.266506361	4.907523682
        15	3.999860139	2.872170674
        20	1.339500378	1.031781873
        30	1.061879319	0.838017806
        40	0.937263075	0.724804683
        50	0.83214763	0.560523347
        60	0.777718424	0.447437618
        80	0.443560634	0.270029374
        100	0.535527683	0.319247051
        110	0.446894137	0.28408831
        120	0.445292862	0.289892677
    };
    \addplot+[
      teal, 
      mark = *,
      mark options={teal, scale=0.8},
      smooth, 
      error bars/.cd, 
        y fixed,
        y dir=both, 
        y explicit
    ] 
    table [x=x, y=y,y error=error] {
        x  y       error
        5	5.706144571	4.376936665
        10	4.356916088	3.395388327
        15	4.458029108	4.391755897
        20	4.522805987	3.752806991
        30	3.036274518	2.366658857
        40	3.163787069	2.315333047
        50	2.915919226	2.515850911
        60	2.538522393	1.701340053
        80	2.283567932	1.327476046
        100	1.756121773	1.082435282
        110	1.629940598	0.887883976
        120	1.69786123	1.14527838
    };
    \end{axis}
    }
\caption{Dataset 2}
\label{fig:Out-of-sample-Test-10}
\end{subfigure}
\hfill
\begin{subfigure}{0.32\textwidth}
\captionsetup{justification=centering}
\footnotesize
\tikz{
\begin{axis}[
  width=\linewidth,
  ytick={0,5,...,65}, ytick align=outside, ytick style={draw=none},
  xtick={0,20,...,120}, xtick align=outside, xtick style={draw=none},
  xlabel=Scenarios,
  ylabel={Out-of-sample deviations (\%)},
  height=\axisdefaultheight,
  grid=both
  ]
    \addplot+[
      BlueViolet, 
      mark=triangle,
      mark options={BlueViolet, scale=0.8},
      smooth,
      error bars/.cd, 
        y fixed,
        y dir=both, 
        y explicit
    ] 
    table [x=x, y=y,y error=error] {
        x  y       error
        5	29.64160672	21.76563642
        10	9.9790095	5.991730756
        15	3.535632549	2.813963069
        20	9.0303377	4.016645081
        30	1.62284016	1.524696926
        40	0.497564458	0.601780776
        50	1.01398492	0.842948379
        60	0.123698526	0.900164549
        80	1.170750849	1.717381934
        100	0.91422301	0.962751461
        110	2.806260275	2.523947497
        120	1.577739325	1.513237871
    };
    \addplot+[
      purple, 
      mark=square,
      mark options={purple, scale=0.8},
      smooth, 
      error bars/.cd, 
        y fixed,
        y dir=both, 
        y explicit
    ] 
    table [x=x, y=y,y error=error] {
        x  y       error
        5	29.37419702	21.26457054
        10	26.38604791	21.5582386
        15	3.123293236	1.955651535
        20	3.414167112	2.603535479
        30	2.245543433	1.462101392
        40	1.224911293	0.794216477
        50	1.093341662	0.649050139
        60	0.83731662	0.640311717
        80	0.701735552	0.559936284
        100	0.604635451	0.458768554
        110	0.49934169	0.381779126
        120	0.496067921	0.320273465
    };
    \addplot+[
      teal, 
      mark = *,
      mark options={teal, scale=0.8},
      smooth, 
      error bars/.cd, 
        y fixed,
        y dir=both, 
        y explicit
    ] 
    table [x=x, y=y,y error=error] {
        x  y       error
        5	8.38892527	6.844826533
        10	5.320879007	4.114942303
        15	5.66262626	5.108597753
        20	4.307482559	4.78345487
        30	3.027234461	2.76222913
        40	3.109151969	2.419139717
        50	2.542972735	2.093377814
        60	2.119500932	1.617372171
        80	1.541064922	1.476211586
        100	1.301379553	1.033388566
        110	1.502057004	1.162280745
        120	1.668945779	1.157735229
    };
    \end{axis}
    }
\caption{Dataset 3}
\label{fig:Out-of-sample-Test-10}
\end{subfigure}
\hfill
\begin{subfigure}{0.32\textwidth}
\captionsetup{justification=centering}
\footnotesize
\tikz{
\begin{axis}[
  width=\linewidth,
  ytick={0,5,...,30}, ytick align=outside, ytick style={draw=none},
  xtick={0,20,...,120}, xtick align=outside, xtick style={draw=none},
  xlabel=Scenarios,
  ylabel={Out-of-sample deviations (\%)},
  height=\axisdefaultheight,
  grid=both
  ]
    \addplot+[
      BlueViolet, 
      mark=triangle,
      mark options={BlueViolet, scale=0.8},
      smooth,
      error bars/.cd, 
        y fixed,
        y dir=both, 
        y explicit
    ] 
    table [x=x, y=y,y error=error] {
        x  y       error
        5	9.445931045	6.82587437
        10	4.085677105	2.974421875
        15	3.86680703	2.782176255
        20	4.276201065	3.007355175
        30	3.16883435	2.41212476
        40	2.04963251	1.620106675
        50	1.77868044	1.57752831
        60	2.47193092	1.685749875
        80	2.379754235	1.141922865
        100	2.11487529	1.143661115
        110	2.057431245	1.175885555
        120	1.53204177	0.89055169
    };
    \addplot+[
      purple, 
      mark=square,
      mark options={purple, scale=0.8},
      smooth, 
      error bars/.cd, 
        y fixed,
        y dir=both, 
        y explicit
    ] 
    table [x=x, y=y,y error=error] {
        x  y       error
        5	8.845449	5.227205603
        10	7.330357988	5.896879884
        15	3.18358652	2.150307104
        20	3.344311936	2.047057892
        30	2.241026436	1.484623656
        40	2.360051332	1.51224142
        50	1.399817312	0.824654308
        60	1.417758528	1.069555824
        80	1.088774196	0.660114056
        100	1.181843412	0.873277472
        110	1.131529724	0.796762548
        120	0.825258612	0.57595688
    };
    \addplot+[
      teal, 
      mark = *,
      mark options={teal, scale=0.8},
      smooth, 
      error bars/.cd, 
        y fixed,
        y dir=both, 
        y explicit
    ] 
    table [x=x, y=y,y error=error] {
        x  y       error
        5	10.09134028	7.70706658
        10	5.73434122	3.76021356
        15	6.939662336	5.513552452
        20	4.116095412	2.986554492
        30	4.827449732	3.209989664
        40	2.75816604	2.113339416
        50	1.444662876	1.317070076
        60	0.565216406	0.603928191
        80	1.112849798	0.762102656
        100	2.012487566	1.68361554
        110	1.287683584	1.254592063
        120	1.666735288	1.280770334
    };
    \end{axis}
    }
\caption{Dataset 4}
\label{fig:Out-of-sample-Test-10}
\end{subfigure}
\hfill
\begin{subfigure}{0.32\textwidth}
\captionsetup{justification=centering}
\footnotesize
\tikz{
\begin{axis}[
  width=\linewidth,
  ytick={0,5,...,55}, ytick align=outside, ytick style={draw=none},
  xtick={0,20,...,120}, xtick align=outside, xtick style={draw=none},
  xlabel=Scenarios,
  ylabel={Out-of-sample deviations (\%)},
  height=\axisdefaultheight,
  grid=both
  ]
    \addplot+[
      BlueViolet, 
      mark=triangle,
      mark options={BlueViolet, scale=0.8},
      smooth,
      error bars/.cd, 
        y fixed,
        y dir=both, 
        y explicit
    ] 
    table [x=x, y=y,y error=error] {
        x  y       error
        5	14.22047804	10.64743709
        10	15.58443708	8.203903584
        15	7.598447407	5.328921924
        20	9.951151307	4.852038469
        30	5.44639592	4.092206089
        40	3.00405581	2.390466326
        50	4.604455951	2.95480345
        60	5.208797855	3.068227506
        80	3.071090988	2.236238247
        100	2.793339757	2.013328437
        110	2.798586099	1.776854669
        120	2.65107939	1.690074413
    };
    \addplot+[
      purple, 
      mark=square,
      mark options={purple, scale=0.8},
      smooth, 
      error bars/.cd, 
        y fixed,
        y dir=both, 
        y explicit
    ] 
    table [x=x, y=y,y error=error] {
        x  y       error
        5	22.27818987	16.64591194
        10	11.00868761	8.287289888
        15	6.418877121	4.238878855
        20	6.944136478	6.007042202
        30	2.591192565	1.697511762
        40	3.491222312	2.143729683
        50	2.634175646	1.52809306
        60	2.125395208	1.176152223
        80	1.983774859	1.380620951
        100	1.158198226	0.942980594
        110	1.469411707	1.057738052
        120	1.323459201	1.105659459
    };
    \addplot+[
      teal, 
      mark = *,
      mark options={teal, scale=0.8},
      smooth, 
      error bars/.cd, 
        y fixed,
        y dir=both, 
        y explicit
    ] 
    table [x=x, y=y,y error=error] {
        x  y       error
        5	12.7618001	9.442415
        10	6.162972181	4.344198762
        15	7.105721591	3.715210267
        20	4.499876511	2.29738458
        30	2.513245128	0.938707593
        40	2.103092102	1.257173393
        50	3.222397863	1.614890866
        60	3.465838571	1.076162885
        80	2.375338367	1.053314609
        100	1.080933729	0.437432339
        110	1.67398491	0.330567491
        120	1.38583917	0.448263238
    };
    \end{axis}
    }
\caption{Dataset 5}
\label{fig:Out-of-sample-Test-10}
\end{subfigure}
\hfill
\begin{subfigure}{0.32\textwidth}
\centering
\vspace{-2cm} 
\ref{sharedlegend}
\end{subfigure}

\caption{Out-of-sample test (weak version): It measures the deviation between the stochastic best found  objective value and the best found  objective value of 10-customer instances obtained by re-evaluating the solution with different $\xi$ vector realizations.}
\label{fig:OutSample10}
\end{figure}
\section{Routing heuristics} \label{subsec:route}

Within the sALNS method, a routing heuristic substitutes an exact method for efficiency. We tested three alternatives:

\begin{itemize}
    \item[-] {Clarke and Wright’s savings heuristic (CW)} \citep{Clarke1964SchedulingPoints}: Starts with one route per customer and iteratively merges routes based on the highest calculated savings, computed for endpoint-customer pairs, while ensuring capacity and time window feasibility.
    \item[-] {Improved Clarke and Wright (ICW)} \citep{Braysy2005VehicleAlgorithms}: Builds on CW by incorporating roulette-wheel–selected local search moves (one-point, two-opt, and two-point exchange) to further refine routes and improve solution quality.
    \item[-] {Cluster First, Route Second (CFRS)}: Clusters customers initially and then constructs each route by greedily adding feasible customers. Customers deemed infeasible in the first pass are reassessed after route completion to determine if they can be inserted or require new routes.
\end{itemize}

We evaluated these heuristics across Solomon’s VRPTW instances \citep{Solomon1987AlgorithmsConstraints}, varying the number of customers (denoted \textit{I\#-C\#}, where \textit{I} is the instance map and \textit{C} the customer count) and adapting original time windows to align with our three predefined delivery slots. Figure~\ref{fig:routing_comparison} presents the  results.

\newcommand{\D}{26} 
\newcommand{\U}{150} 
\newdimen\R 
\R=3.2cm 
\newdimen\L 
\L=3.85cm
\newdimen\l 
\l=.35cm
\newcommand{\A}{360/\D} 
\begin{figure}[h!]
\centering
    \begin{subfigure}{0.45\textwidth}    
        \begin{tikzpicture}[scale=.8]
          \path (0:0cm) coordinate (O); 
        
          \foreach \X in {1,...,\D}{
            \draw [opacity=0.3] (\X*\A:0) -- (\X*\A:\R);
          }
        
          \foreach \Y in {10,30,...,\U}{
            \draw [opacity=0.3] (0:\Y*\R/\U) \foreach \X in {1,...,\D}{
                -- (\X*\A:\Y*\R/\U)
            } -- cycle;
            
            \pgfmathsetmacro{\Label}{\Y/100}
            \draw (90:\Y*\R/\U) node[anchor=south] 
            {\microsize \pgfmathprintnumber[precision=2, fixed, fixed zerofill] \Label}; 
          }
          \foreach \Y in {0,...,\U}{
            \foreach \X in {1,...,\D}{
              \path [opacity=0.3] (\X*\A:\Y*\R/\U) coordinate (D\X-\Y);
              \fill [opacity=0.1] (D\X-\Y) circle (0.5pt);
            }
          }

            \path (1*\A:\L) node (L1) {\microsize I1-C5};
              \path (2*\A:\L) node (L2) {\microsize I1-C10};
              \path (3*\A:\L) node (L3) {\microsize I1-C15};
              \path (4*\A:\L) node (L4) {\microsize I1-C20};
              \path (5*\A:\L) node (L5) {\microsize I1-C30};
              \path (6*\A:\L) node (L6) {\microsize I1-C40};
              \path (7*\A:\L) node (L7) {\microsize I1-C50};
              \path (8*\A:\L) node (L8) {\microsize I1-C60};
              \path (9*\A:\L) node (L9) {\microsize I2-C5};
              \path (10*\A:\L) node (L10) {\microsize I2-C10};
              \path (11*\A:\L) node (L11) {\microsize I2-C15};
              \path (12*\A:\L) node (L12) {\microsize I2-C20};
              \path (13*\A:\L) node (L13) {\microsize I2-C30};
              \path (14*\A:\L) node (L14) {\microsize I2-C40};
              \path (15*\A:\L) node (L15) {\microsize I2-C50};
              \path (16*\A:\L) node (L16) {\microsize I2-C60};
              \path (17*\A:\L) node (L17) {\microsize I2-C80};
              \path (18*\A:\L) node (L18) {\microsize I2-C100};
              \path (19*\A:\L) node (L19) {\microsize I3-C5};
              \path (20*\A:\L) node (L20) {\microsize I3-C10};
              \path (21*\A:\L) node (L21) {\microsize I3-C15};
              \path (22*\A:\L) node (L22) {\microsize I3-C20};
              \path (23*\A:\L) node (L23) {\microsize I3-C30};
              \path (24*\A:\L) node (L24) {\microsize I3-C40};
              \path (25*\A:\L) node (L25) {\microsize I3-C50};
              \path (26*\A:\L) node (L26) {\microsize I3-C60};

          \draw [color=BlueViolet,line width=1.0pt,opacity=0.5]
                (D1-0.00) --
                (D2-8.04) --
                (D3-4.88) --
                (D4-13.47) --
                (D5-14.75) --
                (D6-14.14) --
                (D7-13.17) --
                (D8-14.06) --
                (D9-0.00) --
                (D10-0.00) --
                (D11-0.00) --
                (D12-0.00) --
                (D13-0.00) --
                (D14-0.42) --
                (D15-4.39) --
                (D16-7.20) --
                (D17-8.43) --
                (D18-7.58) --
                (D19-0.00) --
                (D20-10.83) --
                (D21-6.42) --
                (D22-11.64) --
                (D23-11.46) --
                (D24-22.25) --
                (D25-11.94) --
                (D26-17.50) -- cycle;
        
          \draw [color=purple,line width=1.0pt,opacity=0.5, dotted]
                (D1-0.00) --
                (D2-8.04) --
                (D3-4.88) --
                (D4-13.47) --
                (D5-14.75) --
                (D6-14.14) --
                (D7-13.17) --
                (D8-14.06) --
                (D9-0.00) --
                (D10-0.00) --
                (D11-0.00) --
                (D12-0.00) --
                (D13-0.00) --
                (D14-0.42) --
                (D15-4.39) --
                (D16-7.20) --
                (D17-8.43) --
                (D18-7.58) --
                (D19-0.00) --
                (D20-10.83) --
                (D21-6.42) --
                (D22-11.64) --
                (D23-11.46) --
                (D24-22.25) --
                (D25-11.94) --
                (D26-17.50) -- cycle;
        
          \draw [color=PineGreen,line width=1.0pt,opacity=0.8, dashed]
                (D1-9.02) --
                (D2-11.23) --
                (D3-8.15) --
                (D4-6.91) --
                (D5-35.92) --
                (D6-59.02) --
                (D7-59.84) --
                (D8-64.65) --
                (D9-0.00) --
                (D10-7.84) --
                (D11-11.12) --
                (D12-20.42) --
                (D13-38.69) --
                (D14-75.93) --
                (D15-116.89) --
                (D16-103.79) --
                (D17-119.61) --
                (D18-136.46) --
                (D19-37.18) --
                (D20-2.14) --
                (D21-12.94) --
                (D22-15.16) --
                (D23-40.07) --
                (D24-41.47) --
                (D25-60.78) --
                (D26-64.45) -- cycle;

      \draw[very thick, BlueViolet] ($(4.5-0.3,0*.32cm-85)$) -- ++(0.6,0);
      \node[anchor=west] at ($(4.82cm, 0*0.32cm-0.15cm-80)$) {\scriptsize CW};
      \draw[very thick, purple, dotted] ($(4.5-0.3,-1*.32cm-85)$) -- ++(0.6,0);
      \node[anchor=west] at ($(4.82cm, -1*0.32cm-0.15cm-80)$) {\scriptsize ICW};
      \draw[very thick, teal, dashed] ($(4.5-0.3,-2*.32cm-85)$) -- ++(0.6,0);
      \node[anchor=west] at ($(4.82cm, -2*0.32cm-0.15cm-80)$) {\scriptsize CFRS};

        \end{tikzpicture}
        \caption{Comparison of the gap to the optimal objective value.}
        \label{fig:routing_comparison_obj}
    \end{subfigure}
\hfill
    \begin{subfigure}{0.45\textwidth}
        \begin{tikzpicture}[scale=.85]
    	\pgfplotstableread[col sep=space]{data.csv}\csvdata
    	\pgfplotstabletranspose\datatransposed{\csvdata} 
    	\begin{axis}[
    		boxplot/draw direction = y,
    		x axis line style = {opacity=1},
    		axis x line = bottom,
    		axis y line = left,
    		xtick = {1, 2, 3},
    		xticklabel style = {align=center, font=\scriptsize},
    		yticklabel style = {align=center, font=\scriptsize},
    		ylabel style = {align=center, font=\scriptsize},
    		xticklabels = {CW, ICW, CFRS},
    		xtick style = {draw=none}, 
    		ytick style = {draw=none}, 
    		ylabel = {Time (s)},
            cycle list={{BlueViolet},{purple},{teal}},
            every axis plot/.append style={fill,fill opacity=.5},
            height=7.2cm
    	]
    		\foreach \n in {1,...,3} {
    			\addplot+[boxplot] table[y index=\n] {\datatransposed};
    		}
    	\end{axis}
        \end{tikzpicture}
    \caption{Time comparison.}
    \label{fig:routing_comparison_time}
    \end{subfigure}        
\caption{Comparison of three different routing heuristics with exact solution results. CFRS: Cluster First, Route Second heuristic; CW: Clark and Wright's savings heuristic; ICW: Improved Clark and Wright's savings heuristic.}
\label{fig:routing_comparison}
\end{figure}

CW and ICW consistently return identical objective values across all instances, resulting in completely overlapping lines in Figure~\ref{fig:routing_comparison_obj}. However, in terms of computation time (Figure~\ref{fig:routing_comparison_time}), CW is both faster and more stable, while ICW introduces additional variability and longer runtimes without offering any improvement in accuracy. CFRS, on the contrary, shows greater variability and generally performs poorly in both areas. Overall, CW offers the best balance of optimality gap and computational efficiency, which led to its integration into our sALNS method.

\section{Full t-test results for all combination} \label{appSec:t-test}

\begin{landscape}
\begin{table}[]
\centering
\caption{Two-sided comparison of $t$-test statistics and p-values across random, clustered, and mixed map configurations, and datasets (\textit{i.e.} behavioral setting) 1 and 2, for 60, 80, and 100 customers, evaluated over different datasets and map instances. Significant differences (p-values $\leq 0.05$) indicate performance variations across the three spatial distributions.}
\label{tab:my-table}
\resizebox{\linewidth}{!}{%
\begin{tabular}{cccccccccccccccc}
\toprule
Customer \textbackslash Dataset &
  \multicolumn{7}{c}{1} &
   &
  \multicolumn{7}{c}{2} \\ \cline{1-8} \cline{10-16} 
\multirow{7}{*}{60} &
  Map &
  \multicolumn{2}{c}{Random} &
  \multicolumn{2}{c}{Clustered} &
  \multicolumn{2}{c}{Mixed} &
   &
  Map &
  \multicolumn{2}{c}{Random} &
  \multicolumn{2}{c}{Clustered} &
  \multicolumn{2}{c}{Mixed} \\ \cline{2-8} \cline{10-16} 
 &
  1 &
  T-value &
  P-value &
  T-value &
  P-value &
  T-value &
  P-value &
   &
  Instance &
  T-value &
  P-value &
  T-value &
  P-value &
  T-value &
  P-value \\ \cline{2-8} \cline{10-16} 
 &
  2 &
  9.4721 &
  0.00 &
  21.8762 &
  0.00 &
  2.7435 &
  0.01 &
   &
  1 &
  11.9467 &
  0.00 &
  7.5998 &
  0.00 &
  8.712 &
  0.00 \\
 &
  3 &
  12.1426 &
  0.00 &
  16.3073 &
  0.00 &
  5.6824 &
  0.00 &
   &
  2 &
  3.1454 &
  0.00 &
  8.537 &
  0.00 &
  11.5698 &
  0.00 \\
 &
  4 &
  19.0515 &
  0.00 &
  12.4182 &
  0.00 &
  15.9714 &
  0.00 &
   &
  3 &
  9.4631 &
  0.00 &
  6.7611 &
  0.00 &
  1.8968 &
  0.06 \\
 &
  5 &
  16.8987 &
  0.00 &
  7.5774 &
  0.00 &
  12.9098 &
  0.00 &
   &
  4 &
  6.6468 &
  0.00 &
  2.0369 &
  0.04 &
  7.5023 &
  0.00 \\
 &
  5 &
  24.8768 &
  0.00 &
  20.3031 &
  0.00 &
  22.0354 &
  0.00 &
   &
  5 &
  3.3437 &
  0.00 &
  6.7261 &
  0.00 &
  6.5822 &
  0.00 \\ \cline{1-8} \cline{10-16} 
\multirow{7}{*}{80} &
  Map&
  \multicolumn{2}{c}{Random} &
  \multicolumn{2}{c}{Clustered} &
  \multicolumn{2}{c}{Mixed} &
   &
  Map&
  \multicolumn{2}{c}{Random} &
  \multicolumn{2}{c}{Clustered} &
  \multicolumn{2}{c}{Mixed} \\ \cline{2-8} \cline{10-16} 
 &
  Instance &
  T-value &
  P-value &
  T-value &
  P-value &
  T-value &
  P-value &
   &
  Instance &
  T-value &
  P-value &
  T-value &
  P-value &
  T-value &
  P-value \\ \cline{2-8} \cline{10-16} 
 &
  1 &
  9.0698 &
  0.00 &
  24.7011 &
  0.00 &
  16.1647 &
  0.00 &
   &
  1 &
  2.7178 &
  0.01 &
  2.8273 &
  0.01 &
  8.3488 &
  0.00 \\
 &
  2 &
  16.6392 &
  0.00 &
  22.8292 &
  0.00 &
  18.9605 &
  0.00 &
   &
  2 &
  4.4907 &
  0.00 &
  3.5284 &
  0.00 &
  6.2942 &
  0.00 \\
 &
  3 &
  4.2411 &
  0.00 &
  15.6548 &
  0.00 &
  11.0525 &
  0.00 &
   &
  3 &
  4.472 &
  0.00 &
  8.5949 &
  0.00 &
  7.2261 &
  0.00 \\
 &
  4 &
  4.0953 &
  0.00 &
  14.8483 &
  0.00 &
  8.3731 &
  0.00 &
   &
  4 &
  4.9059 &
  0.00 &
  9.5353 &
  0.00 &
  2.7702 &
  0.01 \\
 &
  5 &
  20.7136 &
  0.00 &
  15.2051 &
  0.00 &
  10.7707 &
  0.00 &
   &
  5 &
  8.3606 &
  0.00 &
  7.4901 &
  0.00 &
  11.0601 &
  0.00 \\ \cline{1-8} \cline{10-16}
\multirow{7}{*}{100} &
  Map &
  \multicolumn{2}{c}{Random} &
  \multicolumn{2}{c}{Clustered} &
  \multicolumn{2}{c}{Mixed} &
   &
  Map &
  \multicolumn{2}{c}{Random} &
  \multicolumn{2}{c}{Clustered} &
  \multicolumn{2}{c}{Mixed} \\ \cline{1-8} \cline{10-16} 
 &
  Instance &
  T-value &
  P-value &
  T-value &
  P-value &
  T-value &
  P-value &
   &
  Instance &
  T-value &
  P-value &
  T-value &
  P-value &
  T-value &
  P-value \\ \cline{2-8} \cline{10-16} 
 &
  1 &
  8.7634 &
  0.00 &
  25.0762 &
  0.00 &
  15.6504 &
  0.00 &
   &
  1 &
  10.1256 &
  0.00 &
  2.2215 &
  0.03 &
  9.4065 &
  0.00 \\
 &
  2 &
  14.5533 &
  0.00 &
  5.2355 &
  0.00 &
  17.7989 &
  0.00 &
   &
  2 &
  3.4135 &
  0.00 &
  0.5771 &
  0.57 &
  3.0247 &
  0.00 \\
 &
  3 &
  17.1028 &
  0.00 &
  8.7009 &
  0.00 &
  19.3778 &
  0.00 &
   &
  3 &
  4.3461 &
  0.00 &
  -1.6803 &
  0.10 &
  12.3589 &
  0.00 \\
 &
  4 &
  31.5988 &
  0.00 &
  24.7608 &
  0.00 &
  5.3353 &
  0.00 &
   &
  4 &
  5.5519 &
  0.00 &
  1.7134 &
  0.09 &
  11.966 &
  0.00 \\
 &
  5 &
  13.8561 &
  0.00 &
  17.8784 &
  0.00 &
  9.5472 &
  0.00 &
   &
  5 &
  -0.3099 &
  0.76 &
  4.2756 &
  0.00 &
  7.56 &
  0.00 \\ \bottomrule 
\end{tabular}%
}
\end{table}
\end{landscape}

\begin{landscape}
\begin{table}[]
\centering
\caption{Two-sided comparison of $t$-test statistics and p-values across random, clustered, and mixed map configurations, and datasets (\textit{i.e.} behavioral setting) 3 and 4, for 60, 80, and 100 customers, evaluated over different datasets and map instances. Significant differences (p-values $\leq 0.05$) indicate performance variations across the three spatial distributions.}
\label{tab:my-table}
\resizebox{\linewidth}{!}{%
\begin{tabular}{cccccccccccccccc}
\toprule
Customer \textbackslash Dataset &
  \multicolumn{7}{c}{3} &
   &
  \multicolumn{7}{c}{4} \\ \cline{1-8} \cline{10-16} 
\multirow{7}{*}{60} &
  Map &
  \multicolumn{2}{c}{Random} &
  \multicolumn{2}{c}{Clustered} &
  \multicolumn{2}{c}{Mixed} &
   &
  Map &
  \multicolumn{2}{c}{Random} &
  \multicolumn{2}{c}{Clustered} &
  \multicolumn{2}{c}{Mixed} \\ \cline{2-8} \cline{10-16} 
 &
  Instance &
  T-value &
  P-value &
  T-value &
  P-value &
  T-value &
  P-value &
   &
  Instance &
  T-value &
  P-value &
  T-value &
  P-value &
  T-value &
  P-value \\ \cline{2-8} \cline{10-16} 
 &
  1 &
  10.6861 &
  0.00 &
  2.3497 &
  0.02 &
  15.4567 &
  0.00 &
   &
  1 &
  8.4083 &
  0.00 &
  2.3238 &
  0.02 &
  13.0496 &
  0.00 \\
 &
  2 &
  4.7334 &
  0.00 &
  5.9474 &
  0.00 &
  7.6912 &
  0.00 &
   &
  2 &
  5.9088 &
  0.00 &
  5.8853 &
  0.00 &
  7.4096 &
  0.00 \\
 &
  3 &
  9.0742 &
  0.00 &
  8.9495 &
  0.00 &
  2.0080 &
  0.05 &
   &
  3 &
  7.9070 &
  0.00 &
  8.9312 &
  0.00 &
  1.6414 &
  0.10 \\
 &
  4 &
  10.8969 &
  0.00 &
  5.3191 &
  0.00 &
  10.4703 &
  0.00 &
   &
  4 &
  9.9127 &
  0.00 &
  5.0501 &
  0.00 &
  8.7450 &
  0.00 \\
 &
  5 &
  1.7464 &
  0.08 &
  9.0870 &
  0.00 &
  10.6946 &
  0.00 &
   &
  5 &
  2.2973 &
  0.02 &
  8.3249 &
  0.00 &
  9.8057 &
  0.00 \\ \cline{1-8} \cline{10-16} 
\multirow{7}{*}{80} &
  Map &
  \multicolumn{2}{c}{Random} &
  \multicolumn{2}{c}{Clustered} &
  \multicolumn{2}{c}{Mixed} &
   &
  Map &
  \multicolumn{2}{c}{Random} &
  \multicolumn{2}{c}{Clustered} &
  \multicolumn{2}{c}{Mixed} \\ \cline{2-8} \cline{10-16} 
 &
  Instance &
  T-value &
  P-value &
  T-value &
  P-value &
  T-value &
  P-value &
   &
  Instance &
  T-value &
  P-value &
  T-value &
  P-value &
  T-value &
  P-value \\ \cline{2-8} \cline{10-16} 
 &
  1 &
  5.2963 &
  0.00 &
  -1.3443 &
  0.18 &
  6.5017 &
  0.00 &
   &
  1 &
  5.0342 &
  0.00 &
  -0.0433 &
  0.97 &
  6.2781 &
  0.00 \\
 &
  2 &
  -0.9187 &
  0.36 &
  1.2481 &
  0.21 &
  3.4390 &
  0.00 &
   &
  2 &
  -0.2257 &
  0.82 &
  -0.5766 &
  0.56 &
  3.5192 &
  0.00 \\
 &
  3 &
  1.0689 &
  0.29 &
  3.8798 &
  0.00 &
  5.7931 &
  0.00 &
   &
  3 &
  1.9745 &
  0.05 &
  4.1524 &
  0.00 &
  6.1043 &
  0.00 \\
 &
  4 &
  6.1905 &
  0.00 &
  3.7588 &
  0.00 &
  3.2435 &
  0.00 &
   &
  4 &
  6.5417 &
  0.00 &
  3.2674 &
  0.00 &
  4.4620 &
  0.00 \\
 &
  5 &
  8.6933 &
  0.00 &
  8.9453 &
  0.00 &
  12.9862 &
  0.00 &
   &
  5 &
  8.4489 &
  0.00 &
  9.0538 &
  0.00 &
  11.4594 &
  0.00 \\ \cline{1-8} \cline{10-16} 
\multirow{7}{*}{100} &
  Map &
  \multicolumn{2}{c}{Random} &
  \multicolumn{2}{c}{Clustered} &
  \multicolumn{2}{c}{Mixed} &
   &
  Map &
  \multicolumn{2}{c}{Random} &
  \multicolumn{2}{c}{Clustered} &
  \multicolumn{2}{c}{Mixed} \\ \cline{2-8} \cline{10-16} 
 &
  Instance &
  T-value &
  P-value &
  T-value &
  P-value &
  T-value &
  P-value &
   &
  Instance &
  T-value &
  P-value &
  T-value &
  P-value &
  T-value &
  P-value \\ \cline{2-8} \cline{10-16} 
 &
  1 &
  4.0165 &
  0.00 &
  3.5451 &
  0.00 &
  5.9080 &
  0.00 &
   &
  1 &
  4.9361 &
  0.00 &
  4.6293 &
  0.00 &
  5.5756 &
  0.00 \\
 &
  2 &
  -0.3394 &
  0.73 &
  6.8573 &
  0.00 &
  4.6164 &
  0.00 &
   &
  2 &
  -0.1033 &
  0.92 &
  6.0590 &
  0.00 &
  4.5326 &
  0.00 \\
 &
  3 &
  9.4074 &
  0.00 &
  4.4829 &
  0.00 &
  12.1499 &
  0.00 &
   &
  3 &
  9.9017 &
  0.00 &
  4.0735 &
  0.00 &
  10.8478 &
  0.00 \\
 &
  4 &
  3.0024 &
  0.00 &
  -0.1096 &
  0.91 &
  8.9149 &
  0.00 &
   &
  4 &
  3.5769 &
  0.00 &
  0.9127 &
  0.36 &
  8.6171 &
  0.00 \\
 &
  5 &
  3.1360 &
  0.00 &
  2.9425 &
  0.00 &
  7.9481 &
  0.00 &
   &
  5 &
  3.5328 &
  0.00 &
  3.7970 &
  0.00 &
  8.5296 &
  0.00 \\ \bottomrule
\end{tabular}%
}
\end{table}
\end{landscape}

\begin{table}[]
\tiny
\centering
\caption{Two-sided comparison of $t$-test statistics and p-values across random, clustered, and mixed map configurations, and dataset (\textit{i.e.} behavioral setting) 5, for 60, 80, and 100 customers, evaluated over different datasets and map instances. Significant differences (p-values $\leq 0.05$) indicate performance variations across the three spatial distributions.}
\label{tab:my-table}
\resizebox{\linewidth}{!}{%
\begin{tabular}{cccccccc}
\toprule
Customer \textbackslash Dataset & \multicolumn{7}{c}{5}                                                                             \\ \midrule
\multirow{7}{*}{60}             & Map      & \multicolumn{2}{c}{Random} & \multicolumn{2}{c}{Clustered} & \multicolumn{2}{c}{Mixed} \\ \cline{2-8} 
                                & Instance & T-value      & P-value     & T-value       & P-value       & T-value     & P-value     \\ \cline{2-8} 
                                & 1        & 17.2792      & 0.00        & 15.15         & 0.00          & 5.4641      & 0.00        \\
                                & 2        & 9.1588       & 0.00        & 7.9244        & 0.00          & 14.9602     & 0.00        \\
                                & 3        & 18.1849      & 0.00        & 18.9486       & 0.00          & 12.3058     & 0.00        \\
                                & 4        & 20.3378      & 0.00        & 15.604        & 0.00          & 16.9418     & 0.00        \\
                                & 5        & 14.9309      & 0.00        & 17.6602       & 0.00          & 13.3051     & 0.00        \\ \midrule
\multirow{7}{*}{80}             & Map      & \multicolumn{2}{c}{Random} & \multicolumn{2}{c}{Clustered} & \multicolumn{2}{c}{Mixed} \\ \cline{2-8} 
                                & Instance & T-value      & P-value     & T-value       & P-value       & T-value     & P-value     \\ \cline{2-8} 
                                & 1        & 8.6953       & 0.00        & 4.7777        & 0.00          & 15.3163     & 0.00        \\
                                & 2        & 11.7637      & 0.00        & 10.359        & 0.00          & 15.2971     & 0.00        \\
                                & 3        & 13.9386      & 0.00        & 16.4205       & 0.00          & 15.7669     & 0.00        \\
                                & 4        & 15.9811      & 0.00        & 12.1274       & 0.00          & 11.8162     & 0.00        \\
                                & 5        & 12.0941      & 0.00        & 8.3681        & 0.00          & 20.6263     & 0.00        \\ \midrule
\multirow{7}{*}{100} & Map & \multicolumn{2}{c}{Random} & \multicolumn{2}{c}{Clustered} & \multicolumn{2}{c}{Mixed} \\ \cline{2-8} 
                                & Instance & T-value      & P-value     & T-value       & P-value       & T-value     & P-value     \\  \cline{2-8} 
                                & 1        & 7.3317       & 0.00        & 14.5709       & 0.00          & 10.4639     & 0.00        \\
                                & 2        & 5.8636       & 0.00        & 6.9392        & 0.00          & 12.0040     & 0.00        \\
                                & 3        & 19.9409      & 0.00        & 8.6319        & 0.00          & 19.7068     & 0.00        \\
                                & 4        & 14.2066      & 0.00        & 15.4699       & 0.00          & 12.3752     & 0.00        \\
                                & 5        & 16.4049      & 0.00        & 0.3306        & 0.74          & 17.9086     & 0.00        \\ \bottomrule
\end{tabular}%
}
\end{table}

\end{document}